\pdfoutput=1
\RequirePackage{ifpdf}
\ifpdf 
\documentclass[pdftex]{sigma}
\else
\documentclass{sigma}
\fi

\newcommand{\kk}{\Bbbk}
\newcommand{\Z}{\mathbb{Z}}
\newcommand{\N}{\mathbb{N}}
\newcommand{\Q}{\mathbb{Q}}
\newcommand{\C}{\mathbb{C}}
\newcommand{\R}{\mathbb{R}}

\renewcommand{\hat}[1]{\widehat{#1}}
\renewcommand{\tilde}[1]{\widetilde{#1}}

\newcommand{\opname}[1]{\operatorname{\mathsf{#1}}}

\newcommand{\Spec}{\operatorname{\mathsf{Spec}}}

\newcommand{\pr}{\opname{pr}}

\newcommand{\Gr}{\opname{Gr}}

\renewcommand{\deg}{\opname{deg}}

\newcommand{\Rm}[1]{{\longmapsto}}
\newcommand{\Lm}[1]{{\longmapsfrom}}

\newcommand{\cA}{{\mathcal A}}

\newcommand{\cF}{{\mathcal F}}

\newcommand{\cS}{{\mathcal S}}
\newcommand{\cT}{{\mathcal T}}
\newcommand{\cU}{{\mathcal U}}

\newcommand{\cZ}{{\mathcal Z}}

\newcommand{\uk}{{\underline{k}}}

\newcommand{\uu}{{\underline{u}}}

\newcommand{\tB}{{\widetilde{B}}}

\newcommand{\gen}{\mathbb{L}}

\newcommand{\torus}{\cT}

\newcommand{\clAlg}{{\cA}}

\newcommand{\diag}{{\delta}}

\newcommand{\row}{{\mathrm{row}}}
\newcommand{\tw}{{\tilde{w}}}

\renewcommand{\diag}{{d}}

\newcommand{\var}{\opname{var}}

\newcommand{\Diag}{\opname{Diag}}

\newcommand{\fv}{\opname{f}}
\newcommand{\ufv}{\opname{uf}}

\newcommand{\Li}{\opname{Li}}
\newcommand{\isom}{\stackrel{\sim}{\rightarrow}}

\newcommand{\dCan}{\opname{B}^*}

\renewcommand{\tw}{\opname{tw}}

\newcommand{\Mc}{M^\circ}

\newcommand{\Nufv}{N_{\opname{uf}}}

\newcommand{\Iufv}{I_{\ufv}}
\newcommand{\Ifv}{I_{\fv}}

\newcommand{\upClAlg}{\mathcal{U}^A}
\newcommand{\XClAlg}{\mathcal{U}^X}

\newcommand{\AVar}{\mathbb{A}}
\newcommand{\XVar}{\mathbb{X}}

\newcommand{\LP}{{\mathcal{LP}}}

\newcommand{\bClAlg}{{\overline{\clAlg}}}

\newcommand{\fFac}{\mathcal{Z}}
\newcommand{\fFacRoot}{\fFac^{\Q}}

\renewcommand{\diag}{\opname{d}'}
\renewcommand{\Diag}{D}

\newcommand{\Id}{{\opname{Id}}}
\newcommand{\prin}{{\opname{prin}}}

\newcommand{\col}{\opname{col}}
\renewcommand{\row}{\opname{row}}

\newcommand{\seq}{\boldsymbol{\mu}}

\newcommand{\hLP}{\widehat{\mathcal{LP}}}

\newcommand{\trans}{\opname{P}}

\renewcommand{\tw}{\opname{tw}}
\newcommand{\qO}{{\opname{A_q}}}

\newcommand{\Proj}{\mathrm{\Proj}}
\newcommand{\ee}{\mathbf{e}}
\newcommand{\ff}{\mathbf{f}}
\renewcommand{\uu}{\mathbf{u}}
\newcommand{\uee}{\underline{\ee}}
\newcommand{\uff}{\underline{\ff}}

\newcommand{\high}{\opname{\ufv,\fv}}
\newcommand{\low}{\opname{\fv,\ufv}}

\numberwithin{equation}{section}

\newtheorem{Theorem}{Theorem}[section]
\newtheorem*{Theorem*}{Theorem}
\newtheorem{Corollary}[Theorem]{Corollary}
\newtheorem{Lemma}[Theorem]{Lemma}
\newtheorem{Proposition}[Theorem]{Proposition}
 { \theoremstyle{definition}
\newtheorem{Definition}[Theorem]{Definition}

\newtheorem{Example}[Theorem]{Example}
\newtheorem{Remark}[Theorem]{Remark}
\newtheorem{Assumption}[Theorem]{Assumption}
}

\begin{document}
\allowdisplaybreaks

\newcommand{\arXivNumber}{2201.10284}

\renewcommand{\PaperNumber}{105}

\FirstPageHeading

\ShortArticleName{Twist Automorphisms and Poisson Structures}

\ArticleName{Twist Automorphisms and Poisson Structures}

\Author{Yoshiyuki KIMURA~$^{\rm a}$, Fan QIN~$^{\rm b}$ and Qiaoling WEI~$^{\rm c}$}

\AuthorNameForHeading{Y.~Kimura, F.~Qin and Q.~Wei}

\Address{$^{\rm a)}$~Faculty of Liberal Arts, Sciences and Global Education, Osaka Metropolitan University, Japan}
\EmailD{\href{mailto:ysykimura@omu.ac.jp}{ysykimura@omu.ac.jp}}

\Address{$^{\rm b)}$~School of Mathematical Sciences, Beijing Normal University, P.R.~China}
\EmailD{\href{mailto:qin.fan.math@gmail.com}{qin.fan.math@gmail.com}}

\Address{$^{\rm c)}$~School of Mathematical Sciences, Capital Normal University, P.R.~China}
\EmailD{\href{mailto:wql03@cnu.edu.cn}{wql03@cnu.edu.cn}}

\ArticleDates{Received April 03, 2023, in final form December 01, 2023; Published online December 23, 2023}

\Abstract{We introduce (quantum) twist automorphisms for upper cluster algebras and cluster Poisson algebras with coefficients. Our constructions generalize the twist automorphisms for quantum unipotent cells. We study their existence and their compatibility with Poisson structures and quantization. The twist automorphisms always permute well-behaved bases for cluster algebras. We explicitly construct (quantum) twist automorphisms of Donaldson--Thomas type and for principal coefficients.}


\Classification{13F60; 17B63}

\section{Introduction}\label{sec:intro}

\subsection{Background}

\subsubsection*{Cluster algebras}

The theory of cluster algebras was introduced by Fomin and Zelevinsky
\cite{FominZelevinsky02} as a combinatorial framework to study the
dual canonical bases of quantum groups \cite{Kas:crystal,Lusztig90,Lusztig91}. In this theory, one has the cluster $A$-variety
$\AVar$ (also called the cluster $K_{2}$-variety). It is a scheme
equipped with a cluster structure: $\AVar$ is the union of many tori
which are glued by birational maps called mutations \cite{gross2013birational}.
Let $\upClAlg$ denote the ``function ring'' of $\AVar$ (by which
we mean the ring of global sections of its structure sheaf; it is
the coordinate ring when $\AVar$ is affine). Then $\upClAlg$ is
called the upper cluster algebra (or upper cluster $A$-algebra).
Under a mild assumption (\emph{full rank assumption}), one can endow
$\AVar$ with a Poisson structure \cite{GekhtmanShapiroVainshtein03,GekhtmanShapiroVainshtein05}.
Correspondingly, the upper cluster algebra $\upClAlg$ becomes a Poisson
algebra, which can be naturally quantized \cite{BerensteinZelevinsky05}.

One also has the cluster $X$-variety $\XVar$ (also called the cluster
Poisson variety). It is a scheme equipped with the same cluster structure:
$\XVar$ is the union of many dual tori which are glued by mutations.
It has a canonical Poisson structure. Let $\XClAlg$ denote its function
ring, which is called the cluster Poisson algebra (or upper cluster
$X$-algebra). Then $\XClAlg$ is a Poisson algebra, and one can quantize
it naturally.

Fock and Goncharov \cite{FockGoncharov06a,FockGoncharov09} found
that the cluster varieties $\AVar$ and $\XVar$ naturally arise in
the study of the (higher) Teichm\"uller theory of a surface. They further
conjectured that the upper cluster algebra $\upClAlg$ should possess
a basis naturally parametrized by the tropical points of the cluster
$X$-variety associated to the \emph{Langlands dual }cluster structure
and, conversely, the cluster Poisson algebra $\XClAlg$ should possess
a basis naturally parametrized by the tropical points of the cluster
$A$-variety associated to the Langlands dual cluster structure. (By
\cite{gross2018canonical}, Fock--Goncharov conjecture is true for
many cases, but might not be true in general.)

In view of Fock--Goncharov conjecture, it is important to understand
the upper cluster $A$-algebras, upper cluster $X$-algebras and their
bases.

On the one hand, many well-known cluster $A$-varieties (strictly
speaking, the sets of their rational points) are smooth manifolds.
Examples include the unipotent cells $N_{-}^{w}$ \cite{GeissLeclercSchroeer10},
double Bruhat cells \cite{BerensteinFominZelevinsky05}, and top dimensional
cells of the Grassmannians $\Gr(k,n)$, $k\leq n\in\N$ \cite{scott2006grassmannians}.
On the other hand, there exists few literature on cluster $X$-varieties.
\cite{ip2018cluster,schrader2019cluster,shen2022cluster} embedded the Drinfeld
double quantum groups of Dynkin types to quantized $\XClAlg$. \cite{schrader2019k}
embedded a subalgebra of a $K$-theoretic Coulomb branch to quantized
$\XClAlg$.

The bases for (quantized) $\upClAlg$ have been extensively studied
and they have been related to representation theory and geometry (see
the survey \cite{qin2021cluster}). Not much is known for the bases
of (quantized) $\XClAlg$ (see \cite{allegretti2019categorified,FockGoncharov06a,gross2018canonical} for some results).

\subsubsection*{Twist automorphisms}

\cite{BerensteinFominZelevinsky96,berenstein1997total} introduced
an automorphism $\tilde{\eta}_{w}$ on the unipotent cell $N_{-}^{w}$.
Its gives rise to an automorphism $\eta_{w}$ on the coordinate ring
$\C[N_{-}^{w}]$. $\eta_{w}$ was called a twist automorphism and
has been studied via cluster algebras \cite{GeissLeclercSchroeer10b}.
A quantum unipotent cell $\qO[N_{-}^{w}]$ is a quantum analog of
the coordinate ring $\C[N_{-}^{w}]$, which is defined for all Kac--Moody
types. In a joint work \cite{kimura2017twist}, the first author introduced
the quantum twist automorphism $\eta_{w}$ on $\qO[N_{-}^{w}]$ and
further showed that the dual canonical basis $\dCan$ of $\qO[N_{-}^{w}]$
is permuted by $\eta_{w}$.

By \cite{GeissLeclercSchroeer10,GeissLeclercSchroeer11,GY13,goodearl2020integral},
$\qO[N_{-}^{w}]$ is a (quantum) upper cluster algebra $\upClAlg$.
By a motivational conjecture of Fomin and Zelevinsky \cite{FominZelevinsky02}
and its natural generalization \cite{Kimura10}, the dual canonical
basis of $\qO[N_{-}^{w}]$ should contain all cluster monomials, i.e.,
the monomials of coordinate functions in some toric chart which are globally regular
on $\AVar$. With the help of the existence of the twist automorphism
$\eta_{w}$ and the fact that it permutes the dual canonical basis
$\dCan$, the second author gave a proof of the generalized conjecture
for all cases \cite{qin2020dual}. (See \cite{Kang2018,mcnamara2021cluster,qin2017triangular} for other approaches
and the corresponding results.)

In view of the successful application of the twist automorphism $\eta_{w}$
to the (quantum) upper cluster algebra $\qO[N_{-}^{w}]$, it is natural
to ask for a twist automorphism in general cases.

\subsection{Main results}

\subsubsection*{Preliminaries}

Let $I$ denote a finite set of vertices and $I=I_{\ufv}\sqcup I_{\fv}$
a partition of $I$. Choose skew-symmetrizers $d_{i}\in\Z_{>0}$,
$i\in I$. The vertices in $I_{\ufv}$ and $\Ifv$ are called unfrozen
and frozen, respectively. A~seed $t$ is a~collection $((A_{i})_{i\in I},(X_{i})_{i\in I},B)$
whose matrix $B=(b_{ij})_{i,j\in I}$ is $\Z$-valued such that $\frac{1}{d_{i}}b_{ij}=-\frac{1}{d_{j}}b_{ji}$,
and $A_{i}$, $X_{i}$ are indeterminates called cluster $A$-variables
and cluster $X$-variables, respectively.

Choose a base ring $\kk=\Z\text{ or }\Z\big[v^{\pm\frac{1}{\diag}}\big]$ for the classical cases or the quantum cases, respectively, where $\diag$ is a sufficiently divisible positive integer. We associate to $t$ two Laurent polynomial rings $\LP^{A}=\kk\big[A_{i}^{\pm}\big]_{i\in I}$,
$\LP^{X}=\kk\big[X_{i}^{\pm}\big]_{i\in I}$. Then, if we take $\kk=\C$, $\LP^{\bullet}$ is the
coordinate ring of the corresponding torus $\cT^{\bullet}:=\Spec\LP^{\bullet}$,
where $\bullet$ stands for $A$ or $X$. We denote the Laurent monomials by
$A^{m}$ for any $m=(m_{i})_{i\in I}\in\Z^{I}$
and, similarly, $X^{n}$ for any $n=(n_{i})_{i\in I}\in\Z^{I}$.

Unless otherwise specified, we will use a $v$-twisted product on $\LP^\bullet$ (see Section \ref{sec:quantization}). We will set $v=1$ if we work at the classical level. Let $\cF^{\bullet}$ denote the skew-field of fractions of~$\LP^{\bullet}$. A~monomial map between $\LP^\bullet$ is an algebra
homomorphism sending Laurent monomials to Laurent monomials. We have a canonical
a monomial map $p^{*}\colon\LP^{X}\rightarrow\LP^{A}$. Following~\mbox{\cite{Qin12, qin2017triangular}}, we consider
the pointed elements, i.e., the elements of the form $S_{m}^{A}=A^{m}\cdot F_{m}|_{Z_{k}\mapsto p^{*}X_{k}}$
in~$\LP^{A}$ and~$S_{n}^{X}=X^{n}\cdot F_{n}|_{Z_{k}\mapsto X_{k}}$
in $\LP^{X}$, where $F_{m}$ and $F_{n}$ are polynomials in indeterminates~$Z_{k}$, ${k\in\Iufv}$, with constant term $1$ (called $F$-polynomials).
Here, $Z_{k}\mapsto p^{*}X_{k}$ and $Z_{k}\mapsto X_{k}$ denote
the evaluation of $Z_{k}$, respectively, and $\cdot$ denotes the commutative product.

For any unfrozen vertex $k\in\Iufv$, we have an algorithm called
mutation which generates a~new seed \[\mu_{k}t=((A_{i}(\mu_{k}t))_{i\in I},(X_{i}(\mu_{k}t))_{i\in I},B(\mu_{k}t)).\]
In particular, we have an isomorphism $\mu_{k}^{\bullet}\colon \cF^{\bullet}(\mu_{k}t)\isom\cF^{\bullet}(t)$
called a mutation birational map.\looseness=1

Let $\Delta_{t}^{+}$ denote the set of seeds $t'=\seq t$ obtained
from $t$ by any finite sequence of mutations~$\seq$, where $\seq$ takes the form $\mu_{k_s}\cdots\mu_{k_2}\mu_{k_1}$ and starts from the seed $t$. Note that the mutation sequence $\seq$ depends on $t$, but we omit the symbol $t$ as in standard literature. Recall that
the cluster $\bullet$-variety is the union \smash{$\cup_{t'\in\Delta_{t}^{+}}\cT^{\bullet}(t')$}
such that the tori (local charts) are glued by mutation birational
maps $(\seq^{\bullet})^{*}$ (coordinate change). The upper cluster
$\bullet$-algebra $\mathcal{U}^{\bullet}$ is defined to be its function
ring, which turns out to be \smash{$\mathcal{U}^{\bullet}(t)=\cap_{t'\in\Delta_{t}^{+}}\seq^{*}\LP^{\bullet}(t')$}
by choosing any initial seed $t$ (local chart). It is known that,
up to identification by mutation maps, these objects are independent
of the choice of the initial seed $t$. So we can simply write $\cU^{\bullet}$.

\subsubsection*{Constructions and results}

Let $t=((A_{i}),(X_{i}),B)$ and $t'=\seq t=((A_{i}'),(X_{i}'),B')$
denote two similar seeds, i.e., there exists a~permutation $\sigma$
on $\Iufv$ such that $b_{i,j}=b_{\sigma i,\sigma j}'$, $d_{i}=d_{\sigma i}$,
$\forall i,j\in I_{\ufv}$. We can relabel the vertices~$\Iufv$ when
working with $t'$, so that we can assume $\sigma=\Id$ from now on.

Let $\pr_{\Iufv}$ denote the natural projection from $\Z^{I}$ to
$\Z^{\Iufv}$. Following \cite{Qin12,qin2017triangular}, two pointed
elements \smash{$S_{m}^{A}\in\LP^{A}(t)$} and \smash{$S_{m'}^{A}\in\LP^{A}(t')$}
are said to be similar if $\pr_{\Iufv}m=\pr_{\Iufv}m'$ and they have
the same $F$-polynomial. Similarly, we define two pointed elements
\smash{$S_{n}^{X}\in\LP^{X}(t)$} and \smash{${S_{n'}^{X}\in\LP^{X}(t')}$} to be similar
if $\pr_{\Iufv}B(t)n=\pr_{\Iufv}B(t')n'$ and they have the same $F$-polynomial.

Following \cite{Qin12,qin2017triangular}, a variation map (or a correction
map, a coefficient twist map) $\var_{t}^{\bullet}$ is a~map sending
pointed elements in $\LP^{\bullet}(t)$ to similar pointed elements
in $\LP^{\bullet}(t')$, where $\bullet$ stands for~$A$ or $X$, see Remarks \ref{rem:variation_map_pointed_elements} and \ref{rem:X-variation_map_pointed_elements}.
For the purpose of this paper, we require that $\var_{t}^{\bullet}$
is an algebra homomorphism from $\LP^{\bullet}(t)$ to $\LP^{\bullet}(t')$.

We define a twist endomorphism $\tw_{t}^{\bullet}$ on $\cU^{\bullet}(t)$
to be the composition $\seq^{*}\var_{t}^{\bullet}$ of the mutation
map $\seq^{*}$ with a variation map $\var_{t}^{\bullet}$. In the classical case, it is
called Poisson if it preserves the Poisson structure.
We can show that the construction is independent of the choice of
the initial seed, so we can simply say that $\tw^{\bullet}$ acts
on $\cU^{\bullet}$ (Propositions \ref{prop:change_seed_A_variation} and \ref{prop:change_seed_X_variation}). Note that different choices of variation maps still give different twist endomorphisms.

Let $\fFacRoot$ denote the multiplicative group generated by the
roots of coefficients (frozen cluster $A$-variables $A_{j}$, $j\in\Ifv$).
Let \smash{$\upClAlg_{\fFacRoot}$} denote the ring \smash{$\upClAlg\otimes \kk\big[\fFacRoot\big]$}.
The notion of the twist endomorphism~$\tw^{A}$ can be naturally generalized
for \smash{$\upClAlg_{\fFacRoot}$}.

In the classical case, the existence of the twist endomorphisms $\tw^{A}$ is given by Theorem
\ref{thm:twist_A} \big(see Remark \ref{rem:exist_X_twist} for the existence
of $\tw^{X}$\big):
\begin{itemize}\itemsep=0pt
\item Assume that the full rank assumption holds. Then the twist endomorphisms
$\tw^{A}$ on $\upClAlg_{\fFacRoot}$ exist, and they are in bijection
with the solutions of an inhomogeneous linear equation system.
\end{itemize}
Let $\var_{t'}^{X}$ denote a monomial map from $\LP^{X}(t')$ to
$\LP^{X}(t)$ and $\var_{t}^{A}$ the corresponding monomial map from
$\LP^{A}(t)_{\fFacRoot}$ to $\LP^{A}(t')_{\fFacRoot}$ constructed
using the pullback (see Theorem \ref{thm:twist_X}). Let~$\tw^{\bullet}$
denote their compositions with the mutation maps, respectively. Then
$\tw^{\bullet}$ are related by Theorem~\ref{thm:twist_X}:
\begin{itemize}\itemsep=0pt
\item $\tw^{X}$ is a twist endomorphism (resp.\ twist automorphism) on $\XClAlg$
if and only if $\tw^{A}$ is a~twist endomorphism (resp.\ twist automorphism\footnote{While we only consider upper cluster algebras in this paper, Proposition \ref{prop:change_seed_A_variation} implies that our twist automorphism $\tw^{A}$ restricts to an automorphism on the $\kk[\fFacRoot]$-subalgebra (cluster algebra) generated by the cluster variables.})
on $\upClAlg_{\fFacRoot}$.
\item Assume that $\tw^{X}$ is a twist automorphism on $\XClAlg$.
Then it is Poisson if and only if $\tw^{A}$ and $\bigl(\tw^{X}\bigr)^{-1}$
commute with the natural homomorphism $p^{*}\colon\XClAlg\!\!\rightarrow\upClAlg$,
i.e., ${p^{*}\bigl(\tw^{X}\bigr)^{-1}\!=\!\tw^{A}p^{*}}$.
\end{itemize}
Note that a variation map $\var_{t'}^{X}$ is Poisson if and only if the corresponding linear map satisfies the quadratic equation in Lemma \ref{lem:commute_variation_p}.

We also explicitly construct Poisson (or quantum) twist automorphisms $\tw^{\bullet}$
on $\cU^{\bullet}$ in the following cases:
\begin{itemize}\itemsep=0pt
\item The case $t'=t[1]$ (see Definition \ref{def:Injective-reachable} and Theorem
\ref{thm:twist_DT}): the corresponding twist automorphisms are said
to be of Donaldson--Thomas type (DT-type for short). The original twist
automorphism $\eta_{w}$ is of this type \cite{qin2020dual}.
\item The case when there is a seed of principal coefficients (see Section~\ref{subsec:Principal-coefficients} and
Theorem \ref{thm:twist_principal_coeff}).
\end{itemize}
We prove general results that well-behaved bases for $\cU^{\bullet}$ satisfying Assumptions \ref{assumption:A_basis} or \ref{assumption:X_basis}
are permuted by twist automorphisms $\tw^{\bullet}$ (see Theorems \ref{thm:twist_A_basis}
and \ref{thm:twist_X_basis}). In addition, we proposed a~method for constructing bases of $\XClAlg$ in Theorem \ref{thm:construct_X_bases}.

Recall that we can naturally quantize $\cU^{\bullet}$ if a compatible Poisson
structure is given. Then we can lift a Poisson twist endomorphism
to a quantum twist endomorphism. Conversely, the classical limit of a quantum twist endomorphism is a Poisson twist endomorphism. See Remarks~\ref{rem:Poisson_quantum_variation} and \ref{rem:X-Poisson_quantum_variation}.

\begin{Remark}[a comparison with the previous literature]
The second author introduced the correction technique to compare similar
pointed elements for similar seeds $t$ and $t'$ in~\cite{Qin12}.
In order to facilitate the comparison, in \cite{qin2017triangular},
he introduced a variation map $\var_{t}^{A}$ sending pointed elements
in $\LP^{A}(t)$ to similar pointed elements in $\LP^{A}(t')$. The
variation map $\var^{A}\colon\LP^{A}(t)_{\fFacRoot}\rightarrow\LP^{A}(t')_{\fFacRoot}$
considered in this paper is defined in the same spirit, but chosen
slightly differently, such that it becomes an algebra homomorphism.
The variation map $\var^{X}\colon\LP^{X}(t')\rightarrow\LP^{X}(t)$ has
not been considered before.

In the classical case $\kk=\Z$, it is easy to see that a variation map $\var^{A}\colon\LP^{A}(t)\rightarrow\LP^{A}(t')$
in this paper is the same as a quasi-homomorphism introduced by Fraser
\cite{fraser2016quasi}. When one identifies~$\cF^{A}(t)$
and $\cF^{A}(t')$ using the mutation map $\seq^{*}$, a twist endomorphism
$\tw^{A}$ on $\upClAlg$ becomes the same as a quasi-homomorphism
for a normalized seed pattern in the sense of \cite{fraser2016quasi}. \cite{chang2019cluster} studied the group of the twist automorphisms
on $\upClAlg$ for principal coefficient cases and several other special
cases.

To the best of the authors' knowledge, the Poisson or quantum twist endomorphisms
as well as the twist endomorphisms $\tw^{X}$ have not been introduced
in the previous literature (although specific examples have risen
from Lie theory and from higher Teichm\"uller theory, see Section~\ref{sec:Examples}).
\end{Remark}

\begin{Remark}[morphisms for schemes]
Consider the classical cases.

On the one hand, our twist endomorphism $\tw^{\bullet}$ provides
an endomorphism for the affine scheme $\Spec\cU^{\bullet}$.

On the other hand, a cluster variety might not be affine, and we do
not know if $\tw^{\bullet}$ provide an endomorphism for it. While
the variation map is always a scheme morphism, the mutations might
be not.

See \cite{gross2013birational} for the comparison between the schemes.
\end{Remark}

\subsection{Contents}

We provide preliminaries for the theory of cluster algebras in Section
\ref{sec:Preliminaries}.

In Sections \ref{sec:Twist-endomorphisms-A} and \ref{sec:Twist-endomorphisms-X},
we introduce (Poisson, quantum) twist endomorphisms for upper cluster $A$-algebras
and upper cluster $X$-algebras, respectively. We discuss their existence.

In Section \ref{sec:Twist-automorphisms-special-cases}, we explicitly
construct (quantum) twist automorphisms for two special cases: the Donaldson--Thomas
type and the principal coefficients.

In Section \ref{sec:Permutation-on-bases}, we prove that a basis
with nice properties is permuted by a twist automorphism. We also construct bases for cluster Poisson algebras.

In Section \ref{sec:Examples}, we give some explicit examples for
Poisson twist automorphisms.

\subsection{Convention}

We fix a finite set of vertices $I$ together with a partition $I=I_{\ufv}\sqcup I_{\fv}$.
The elements in $I_{\ufv}$ and~$I_{\fv}$ are said to be unfrozen
and frozen, respectively. We also fix skew-symmetrizers $d_{i}\in\Z_{>0}$
for~$i\in I$. Let
$\Diag$ denote the diagonal matrix whose diagonal entries are $\frac{1}{d_{i}}$, $i\in I$.

We choose a base ring $\kk$ of characteristic
$0$ and a unit $v\in \kk$. For classical cluster algebras, we choose $\kk=\Z$ and $v=1$. For quantum cluster algebras, we choose \smash{$\kk=\Z\big[v^{\pm\frac{1}{\diag}}\big]$} where $v$ is an indeterminate and $\diag\in \N_{>0}$ is sufficiently divisible. The ring multiplication for cluster algebras will be the $v$-twisted product $*$. We will also use the commutative product $\cdot$.

An $I\times I$ matrix $H$ is said to be skew-symmetrizable by $\Diag$
if $\Diag H$ is skew-symmetric.

For any $I\times I$-matrix $H$ and any $J_{1},J_{2}\subset I$,
let $H_{J_{1},J_{2}}$ denote the $J_{1}\times J_{2}$-submatrix of
$H$. Then we can denote $H$ as a block matrix:
\begin{align*}
H & =\left(\begin{matrix}
H_{\Iufv,\Iufv} & H_{\Iufv,\Ifv}\\
H_{\Ifv,\Iufv} & H_{\Ifv,\Ifv}
\end{matrix}\right)=:\left(\begin{matrix}
H_{\ufv} & H_{\high}\\
H_{\low} & H_{\fv}
\end{matrix}\right).
\end{align*}

Let $\col_{i}H$ and $\row_{i}H$ denote the $i$-th column and the $i$-th
row of $H$, respectively. For any permutation $\sigma$ of $I$ (resp.\
of $\Iufv$), let $\trans_{\sigma}$ denote the $I\times I$-matrix
(resp.\ the $\Iufv\times\Iufv$-matrix) such that $\col_{i}\trans_{\sigma}$
is the $\sigma i$-th unit vector. Then we have $\col_{k}(HP_{\sigma})=\col_{\sigma k}H$ as column vectors and~$\row_{i}(P_{\sigma}H)=\row_{\sigma^{-1}i}H$ as row vectors. We work with column vectors unless otherwise specified.

Assume that $\sigma$ is a permutation on $\Iufv$. If $w=(w_k)_{k\in \Iufv}$ is a vector in $\R^{\Iufv}$, we define~$\sigma w\in \R^{\Iufv}$ such that $(\sigma w)_{\sigma k}= w_k$. Similarly, if $\seq=\mu_{k_{r}}\cdots\mu_{k_{2}}\mu_{k_{1}}$ is a sequence of
mutations with $k_i\in \Iufv$, then $\sigma\seq$
denotes $\mu_{\sigma k_{r}}\cdots\mu_{\sigma k_{2}}\mu_{\sigma k_{1}}$.

For a $\Z$-lattice $L$, we denote $L_{\Q}=L\otimes_{\Z}\Q$. If
$L$ has an $I$-labeled basis $\{u_{i}\mid i\in I\}$, let $L_{\ufv}$
and~$L_{\fv}$ denote its sublattices spanned by $\{u_{i}\mid i\in I_{\ufv}\}$
and $\{u_{i}\mid i\in I_{\fv}\}$, respectively.

Let $\pr_{\Iufv}$ denote the natural projection from $\Z^{I}$ to
$\Z^{\Iufv}$.

\section{Preliminaries}\label{sec:Preliminaries}

\subsection{Seeds and tori}

\begin{Definition}[seeds]
A seed $t$ is a collection $((A_{i})_{i\in I},(X_{i})_{i\in I},B)$,
where $A_{i}$ and $X_{i}$ are indeterminates called cluster $A$-variables
and cluster $X$-variables, respectively, and $B=(b_{ij})_{i,j\in I}$
is a $\Z$-valued matrix such that $\Diag B$ is skew-symmetric.
\end{Definition}

Define Laurent polynomial rings \smash{$\LP^{A}=\kk\big[A_{i}^{\pm}\big]_{i\in I}$}
and \smash{$\LP^{X}=\kk\big[X_{i}^{\pm}\big]_{i\in I}$}. For any vector $m,n\in\Z^{I}$,
denote the Laurent monomials $A^{m}=\prod A_{i}^{m_{i}}$ and $X^{n}=\prod X_{i}^{n_{i}}$. We use $\cdot$ to denote their commutative multiplication, which is often omitted for simplicity.

We associate to $t$ a lattice $N=\Z^{I}$ with the natural basis
$\ee_{i}$, $i\in I$, whose elements $n=\sum n_{i}\ee_{i}=(n_{i})\in\Z^{I}$
are viewed as the Laurent degrees for $X^{n}\in\LP^{X}$. Let $M$
denote its dual lattice with the dual basis $\ee_{i}^{*}$. Let $\langle\ ,\ \rangle$
denote the natural pairing between $N_{\Q}$ and $M_{\Q}$. Define
$\ff_{i}=\frac{1}{d_{i}}\ee_{i}^{*}$ and the sublattice $\Mc=\oplus_{i\in I}\Z\ff_{i}\subset M_{\Q}$.
We identify $\Mc$ with $\Z^{I}$ such that~$\ff_{i}$ become the
unit vectors. View its elements \smash{$m=\sum m_{i}\ff_{i}=(m_{i})\in\Z^{I}$}
as the Laurent degrees for~$A^{m}\in\LP^{A}$.

Let $\uee$ and $\uff$ denote the $I$-labeled bases $(\ee_{i})_{i\in I}$
and $(\ff_{i})_{i\in I}$, respectively. We can view them as matrices
whose columns are the basis elements. Define the linear map $p^{*}\colon N\rightarrow\Mc$
such that~$p^{*}\ee_{j}=\sum_{i}b_{ij}\ff_{i}$, $\forall j$. It
has the following matrix presentation:
\begin{align*}
p^{*}(\uee) & =\uff B.
\end{align*}

In the classical case, extend $\kk$ to a field containing $\Z$. Denote the affine schemes $\torus^{A}=\Spec\LP^{A}$ and $\torus^{X}=\Spec\LP^{X}$.
We will call $\torus^{A}$ and $\torus^{X}$ the tori associated to
$t$, since $\torus^{A}(\kk)$ and $\torus^{X}(\kk)$ coincide with
the split algebraic torus $(\kk^{\times})^{I}$.

\subsection{Poisson structures}

Recall that $d$ is the least common multiplier of all $d_i$. We can associate to $N$ a $\frac{1}{d}\Z$-valued canonical skew-symmetric bilinear
form $-\omega$ such that $\omega(\ee_{i},\ee_{j})=\frac{1}{d_{j}}b_{ji}$,
$\forall i,j$. The corresponding canonical Poisson structure on \smash{$\cT^{X}$}
is given by
\begin{align}
\bigl\{X^{n},X^{n'}\bigr\} =-\omega(n,n')X^{n+n'}.\label{eq:log-Poisson}
\end{align}

Let $W$ denote the $\Q$-valued matrix $(\omega(\ee_{i},\ee_{j}))_{i,j}$. Recall that the diagonal entries of $D$ are~$\frac{1}{d_i}$. Then we have \smash{$WD^{-1}=B^{\rm T}$}.

\begin{Remark}
For any $i\in I$, we have $p^{*}\ee_{i}(\ )=\omega(\ee_{i},\ )$.
\end{Remark}

\begin{Remark}
Following \cite{GekhtmanShapiroVainshtein03,GekhtmanShapiroVainshtein05},
the Poisson structure of the form \eqref{eq:log-Poisson} is usually
called log-canonical in the sense that the following holds (for $\kk=\R$):
\begin{align*}
\big\{\ln X^{n},\ln X^{n'}\big\} =-\omega(n,n').
\end{align*}
\end{Remark}

We often impose the following assumption. Denote the matrix $\tB:=(b_{ik})_{i\in I,k\in\Iufv}$.

\begin{Assumption}[full rank assumption]\label{fullrank}
We assume that the linear map $p^{*}$ restricts to an injective map
on $\Nufv$. Equivalently, the matrix $\tB$ is of full rank.
\end{Assumption}

\begin{Definition}[{\cite{BerensteinZelevinsky05}}]\label{def:compatible_poisson_A}
By a (compatible) Poisson structure on $\cT^{A}$, we mean a collection
of strictly positive numbers $\diag_{k} \in \Q_{>0}$, $k\in\Iufv$, and a $\Q$-valued
skew-symmetric bilinear form $\lambda$ on~$\Mc$ such that
\begin{gather*}
\lambda(\ff_{i},p^{*}\ee_{k}) = -\delta_{ik}\diag_{k}
\end{gather*}
for any $i\in I$, $k\in\Iufv$.

The corresponding Poisson bracket on $\cT^{A}$ is defined as
\begin{gather*}
\bigl\{A^{m},A^{m'}\bigr\} = \lambda(m,m')A^{m+m'}.
\end{gather*}

The bilinear form $\lambda$ is represented by the $\Lambda$-matrix
$\Lambda=(\Lambda_{ij})_{i,j\in I}:=(\lambda(\ff_{i},\ff_{j}))_{i,j\in I}$.
\end{Definition}

Note that the existence of a compatible Poisson structure implies the full rank assumption (see Assumption~\ref{fullrank}). Conversely, by \cite{GekhtmanShapiroVainshtein03,GekhtmanShapiroVainshtein05},
when the full rank assumption is satisfied, we can always choose a
(not necessarily unique) compatible Poisson structure.

From now on, we will make
the choice such that \smash{$\diag_{k}=\frac{1}{d_{k}}$}, $\forall k\in\Iufv$
(see \cite[Proposition 3.3]{BerensteinZelevinsky05}). Correspondingly, $\lambda$ is $\frac{1}{\diag}\Z$-valued for some sufficiently divisible $\diag\in \N_{>0}$. We choose $\diag$ such that~$d| \diag$.

Note that the bilinear form $\lambda$ on $\Mc(t)$ naturally induces
a bilinear form $\lambda$ on $\Nufv$ such that
\begin{gather*}
\lambda(n,n') :=\lambda(p^{*}n,p^{*}n').
\end{gather*}
It is easy to check that $\lambda$ on $\Nufv$ coincides with $-\omega$
(see \cite[Lemma 2.1.11\,(2)]{qin2020dual}).

\begin{Definition}[connected matrix]\label{def:connected_matrix}
An $I\times I$ matrix $B$ is said to be connected, if for any~$i,j\in I$,
there exists finitely many vertices $i_{s}\in I$, $0\leq s\leq l$,
such that $i_{0}=i$, $i_{l}=j$, and $b_{i_{s}i_{s+1}}\neq0$ for
any $0\leq s\leq l-1$.
\end{Definition}

It is straightforward to check the following result.

\begin{Lemma}\label{lem:connected_symmetrizer}
Assume that $B$ is a connected $I\times I$ matrix. If $DB$
and $D'B$ are both skew-symmetric for some invertible diagonal matrices
$D$, $D'$, then there exists some $0\neq\alpha\in\Q$ such that
$D=\alpha D'$.
\end{Lemma}

\subsection{Quantum torus algebras}\label{sec:quantization}

Recall that we have chosen $d|\diag$ and $\omega$ on $N$ is $\frac{1}{d}\Z$-valued. Using the
canonical Poisson structure on $\torus^{X}$, we equip $\LP^{X}$ with an extra multiplication $*$, called the $v$-twisted product, such
that
\begin{align*}
X^{n}*X^{n'} & =v^{-\omega(n,n')}X^{n+n'}.
\end{align*}
Then $\LP^{X}$ is called a quantum torus algebra for the quantum case $\kk=\Z\big[v^{\pm\frac{1}{\diag}}\big]$. The canonical Poisson bracket
\eqref{eq:log-Poisson} can be recovered from the twisted product
by
\begin{align*}
\bigl\{X^{n},X^{n'}\bigr\} =\lim_{v\rightarrow1}\frac{1}{2(v-1)}(X^{n}*X^{n'}-X^{n'}*X^{n}).
\end{align*}

If there exists a $\frac{1}{\diag}\Z$-valued compatible Poisson structure
$\lambda$ on $\torus^{A}$ (see Definition~\ref{def:compatible_poisson_A}), we can equip $\LP^{A}$ with an extra multiplication $*$, called the $v$-twisted product, such that
\begin{align*}
 A^{m}*A^{m'} =v^{\lambda(m,m')}A^{m+m'}.
 \end{align*}
Then $\LP^{A}$ is called a quantum torus algebra for the quantum case $\kk=\Z\big[v^{\pm\frac{1}{\diag}}\big]$. As before, the compatible Poisson structure can be recovered from the twisted product
by
\begin{align*}
\bigl\{A^{m},A^{m'}\bigr\} =\lim_{v\rightarrow1}\frac{1}{2(v-1)}(A^{m}*A^{m'}-A^{m'}*A^{m}).
\end{align*}

For the classical case $\kk=\Z$, we have $v=1$ and we simply define the $v$-twisted product $*$ of~$\LP^\bullet$ to be the commutative products $\cdot$. From now on, we always view $\LP^\bullet$ as a $\kk$-algebra whose multiplication is the $v$-twisted product.

Let $\LP=\kk[L]=\oplus \kk_{u\in L} \chi^u$ and $\LP'=\kk[L']=\oplus \kk_{u'\in L'} \chi^{u'}$ be two quantum torus algebras as above, viewed as $\kk$-modules. Here $L$ and $L'$ denote the lattices of Laurent degrees, respectively, and $\chi$ can denote $X$ or $A$. By a monomial map $\varPhi$ from $\LP$ to $\LP'$ we mean a $\kk$-linear map such that, there exists some linear map $\varPsi\colon L\rightarrow L'$ satisfying $\varPhi(\chi^u)=\chi^{\varPsi(u)}$ $\forall u\in L$. In this case, $\varPhi$~is called the monomial map associated to $\varPsi$ and $\varPsi$ the linear map associated to $\varPhi$.

In particular, the linear map $p^*\colon N\rightarrow \Mc$ determines a monomial map $\LP^X\rightarrow \LP^A$ sending~$X^n$ to $A^{p^* n}$, which is still denoted by $p^*$ for simplicity. Then $p^*$ is an algebra homomorphism preserving the $v$-twisted products.

Let $T$ denote the subalgebra $\kk[X_k]_{k\in \Iufv}$ of the quantum torus algebra $\LP^X$. We introduce the $\kk$-algebra \smash{$\hLP^X=\LP^X\otimes_T \hat{T}$}, where \smash{$\hat{T}$} is the completion of $T$ with respect to its maximal ideal generated by $X_k$, $k\in \Iufv$. Similarly, define \smash{$\hLP^A:=\LP^A\otimes_{p^*T}p^*\hat{T}$}. The elements in \smash{$\hLP^\bullet$} will be called formal Laurent series. Note that $p^*$ extends to a homomorphism from \smash{$\hLP^X$} to~\smash{$\hLP^A$}.

\subsection{Mutation maps}\label{sec:mutation}

Denote $[\ ]_{+}=\max(\ ,0)$. For any vector $(g_{i}),$ denote $[(g_{i})]_{+}=([g_{i}]_{+})$.

Let $t$ denote a given seed and $k$ any chosen unfrozen vertex.
Choose any sign $\varepsilon\in\{1,-1\}$. We define a new seed $t'=\mu_{k}t=((A_{i}'),(X_{i}'),B')$,
such that the matrix $B'=(b_{ij}')$ is given by
\begin{align*}
b_{ij}' & =\begin{cases}
-b_{ij}, & k\in\{i,j\},\\
b_{ij}+b_{ik}[\varepsilon b_{kj}]_{+}+[-\varepsilon b_{ik}]_{+}b_{kj}, & k\notin\{i,j\}.
\end{cases}
\end{align*}

Let $\cF^\bullet$ denote the fraction fields of the Laurent polynomial rings $\LP^\bullet$ for classical cases and the skew fields of fractions of the quantum torus algebras $\LP^\bullet$ for quantum cases. We always call $\cF^\bullet$ fraction fields for simplicity.

We further relate the cluster variables for $t'$ with those for $t$,
by introducing isomorphisms for the fraction fields $\big(\mu_{k}^{X}\big)^{*}\colon\big(\cF^{X}\big)'\simeq\cF^{X}$
and $\big(\mu_{k}^{A}\big)^{*}\colon\big(\cF^{A}\big)'\simeq\cF^{A}$. For classical cases, we define
\begin{gather*}
\big(\mu_{k}^{X}\big)^{*}X_{i}' =\begin{cases}
X_{i}X_{k}^{[\varepsilon b_{ki}]_{+}}(1+X_{k}^{\varepsilon})^{-b_{ki}},& i\neq k,\\
X_{k}^{-1}, & i=k
\end{cases}\\
\hphantom{\big(\mu_{k}^{X}\big)^{*}X_{i}'}{} =\begin{cases}
X_{i}(1+X_{k})^{-b_{ki}}, & b_{ki}\leq0,\\
\displaystyle X_{i}X_{k}^{b_{ki}}\frac{1}{(1+X_{k})^{b_{ki}}},& b_{ki}>0,\\
X_{k}^{-1}, & i=k,
\end{cases}\\
\big(\mu_{k}^{A}\big)^{*}A_{i}' =\begin{cases}
 A_{i}, & i\neq k,\\
 \displaystyle A_{k}^{-1}\prod_{j}A_{j}^{[-\varepsilon b_{jk}]_{+}}(1+p^{*}X_{k}^\varepsilon), & i=k
\end{cases} \\
\hphantom{\big(\mu_{k}^{A}\big)^{*}A_{i}'}{}
=\begin{cases}
A_{i}, & i\neq k,\\
\displaystyle A_{k}^{-1}\prod_{j}A_{j}^{[-b_{jk}]_{+}}(1+p^{*}X_{k}), & i=k.
\end{cases}
\end{gather*}
Following \cite[equation (2.6)]{qin2019bases} and \cite[equations (2.2) and (2.4)]{qin2020dual}, for quantum cases, we define $\big(\mu_{k}^{A}\big)^{*}$ such that
 \begin{align*}
 \big(\mu_{k}^{A}\big)^{*}A_{i}' =\begin{cases}
 A_{i}, & i\neq k,\\
 A^{-\ff_k+\sum_j [-b_{jk}]_{+}\ff_j}+A^{-\ff_k+\sum_i [b_{ik}]_{+}\ff_i}, & i=k.
 \end{cases}
 \end{align*}
We define $\big(\mu_{k}^{X}\big)^{*}$ such that
\begin{align*}
 \big(\mu_{k}^{X}\big)^{*}X_{i}' =\begin{cases}
 X_{i}\cdot \displaystyle\sum_{s=0}^{-b_{ki}}\begin{pmatrix}
 -b_{ki}\\
 s\\
 \end{pmatrix}_{v_k} X_{k}^{s}, & b_{ki}\leq0,\\
 \Biggl(\displaystyle\sum_{s=0}^{b_{ki}}\begin{pmatrix}
 b_{ki}\\
 s\\
 \end{pmatrix}_{v_k} X_{i}^{-1}\cdot X_{k}^{-b_{ki}} \cdot X_{k}^{s} \Biggr)^{-1},& b_{ki}>0,\\
 X_{k}^{-1}, & i=k,
 \end{cases}
 \end{align*}
where $\cdot$ denotes the commutative product. Here, we denote $v_k:=v^{\diag_{k}}=v^{\frac{1}{d_k}}$. The quantum numbers for $b\leq a \in \N$ are defined as
\[[a]_v:=\frac{v^a-v^{-a}}{v-v^{-1}},\qquad [a]_v!:=[a]_v[a-1]_{v}\cdots [1]_v, \qquad \begin{pmatrix}
 a\\
 b \end{pmatrix}_{v_k}:=\frac{[a]_v!}{[b]_v! [a-b]_v !}.\]

The maps $(\mu_{k}^{\bullet})^{*}$ are called the \emph{mutation
maps}, where $\bullet$ stands for $A$ or $X$\emph{.} For classical cases, they induce
birational maps $\mu_{k}^{\bullet}\colon\torus^{\bullet}(t)\dashrightarrow\torus^{\bullet}(t')$.\footnote{We use dashed arrows to denote rational maps.}

It is straightforward to check the well-known fact that $B'$ does not depend on the choice of the sign $\varepsilon$. Moreover,
mutation is an involution on the seeds, and the compositions \[((\mu_{k}^{\bullet})^{*})^{2}\colon \ (\cF^{\bullet})'\simeq\cF^{\bullet}\simeq(\cF^{\bullet})'\]
are the identity. One can also check that mutations commute with the monomial map $p^*$:
\begin{align}\label{eq:X_A_mutation_commutativity}
 \big(\mu_k^A\big)^*\circ p^*(t')(X_i(t'))=p^*(t)\circ \big(\mu_k^X\big)^*(X_i(t')),\qquad \forall i\in I.
\end{align}

For simplicity, we might omit the symbol $\bullet$ when there is
no confusion.

\subsubsection*{Hamiltonian formalism}

For classical cases, let us recall the Hamiltonian formalism for the mutation maps following
\cite{gekhtman2017hamiltonian} (see also \cite{FockGoncharov09a,gross2018canonical}). Recall that the Euler dilogarithm function
is given by
\begin{align*}
\Li_{2}(z) =\sum_{n\geq1}\frac{z^{n}}{n^{2}}.
\end{align*}

Let \smash{$\big(\rho_{k,\varepsilon}^{X}\big)^{*}$} denote the automorphism $\exp\{\varepsilon d_{k}\Li_{2}(-X_{k}^{\varepsilon}),\ \}$
on $\cF^{X}$. When $\cF^{A}$ possesses a~compatible Poisson structure,
let \smash{$\big(\rho_{k,\varepsilon}^{A}\big)^{*}$} denote its automorphism $\exp\big\{\varepsilon\Li_{2}\bigl(-A^{\varepsilon p^{*}\ee_{k}}\bigr),\ \big\}$.
One can check that
\begin{gather*}
\big(\rho_{k,\varepsilon}^{X})^{*}(X^{n}\big) = X^{n}\big(1+X_{k}^{\varepsilon}\big)^{\omega(n,-d_{k}\ee_{k})},\qquad
\big(\rho_{k,\varepsilon}^{A}\big)^{*}(A^{m}) = A^{m}\big(1+A^{\varepsilon p^{*}\ee_{k}}\big)^{-m_{k}}.
\end{gather*}
In particular, we have
\begin{gather*}
\big(\rho_{k,\varepsilon}^{X}\big)^{*}(X_{i}) = \begin{cases}
X_{i}(1+X_{k}^{\varepsilon})^{-b_{ki}}, & i\neq k,\\
X_{k}, & i=k,
\end{cases}\\
\big(\rho_{k,\varepsilon}^{A}\big)^{*}(A_{i}) = \begin{cases}
A_{i}, & i\neq k,\\
A_{k}(1+(p^{*}X_{k})^{\varepsilon})^{-1}, & i=k.
\end{cases}
\end{gather*}

We define monomial maps $\psi_{k,\varepsilon}^{X}\colon\big(\LP^{X}\big)'\simeq\LP^{X}$
and $\psi_{k,\varepsilon}^{A}\colon\big(\LP^{A}\big)'\simeq\LP^{A}$ such that
\begin{gather*}
\psi_{k,\varepsilon}^{X}X_{i}' = \begin{cases}
X_{i}\cdot X_{k}^{[\varepsilon b_{ki}]_{+}}, & i\neq k,\\
X_{k}^{-1}, & i=k,
\end{cases}\qquad
\psi_{k,\varepsilon}^{A}A_{i}' = \begin{cases}
A_{i}, & i\neq k,\\
\displaystyle A_{k}^{-1}\cdot \prod_{j}A_{j}^{[-\varepsilon b_{jk}]_{+}}, & i=k,
\end{cases}
\end{gather*}
where the commutative products are used. Then the mutation map $(\mu_{k}^{\bullet})^{*}$ decomposes as the
composition $(\rho_{k,\varepsilon}^{\bullet})^{*}\circ\psi_{k,\varepsilon}^{\bullet}$.
The factor $\psi_{k,\varepsilon}^{\bullet}$ is called its monomial
part and $(\rho_{k,\varepsilon}^{\bullet})^{*}$ its Hamiltonian part. Note that
we have $(\rho_{k,\varepsilon}^{\bullet})^{*}\circ\psi_{k,\varepsilon}^{\bullet}=\psi_{k,\varepsilon}^{\bullet}\circ(\rho_{k,-\varepsilon}^{\bullet})^{*}$.

By \cite{BerensteinZelevinsky05}, if $t$ is endowed with a compatible Poisson structure, then we have a unique compatible Poisson structure for $t'$ such
that the homomorphism $\psi^{A}_{k,\varepsilon}$ is a Poisson homomorphism. In addition, it is independent of $\varepsilon$. We still use $\lambda$ to denote the corresponding skew-symmetric
bilinear form on $(\Mc)'$, and denote the corresponding matrix by
$\Lambda'=(\lambda(f_{i}',f_{j}'))_{i,j\in I}$.

\subsubsection*{Mutation sequences}
Let $\uk=(k_{1},\ldots,k_{l})$ denote a finite sequence of unfrozen vertices.
A sequence of mutations ${\seq=\seq_{\uk}}$ is the composition of mutations
$\mu_{k_{l}}\cdots\mu_{k_{2}}\mu_{k_{1}}$ (read from right to left).
The corresponding mutation map is $(\seq^{\bullet})^{*}=(\mu_{k_{1}}^{\bullet})^{*}(\mu_{k_{2}}^{\bullet})^{*}\cdots(\mu_{k_{l}}^{\bullet})^{*}$,
which we often write $\seq^{*}$ for simplicity. Note that $\seq^{-1}=\mu_{k_{1}}\mu_{k_{2}}\cdots \mu_{k_{l}}$ and \smash{$(\seq^*)^{-1}=\big(\seq^{-1}\big)^{*}$}.

We deduce the following equality in $\cF^A(t)$ from \eqref{eq:X_A_mutation_commutativity}:
\begin{align}\label{eq:X_A_mutation_seq_comm}
 \big(\seq^A\big)^*\circ p^*(t')(X(t')^n)=p^*(t)\circ \big(\seq^X\big)^*(X(t')^n),\qquad \forall n\in \Z^I.
\end{align}

\subsubsection*{Formal Laurent series expansions}
We recall the maps $\iota^\bullet$ taking formal Laurent series expansions following \cite[Section 3.3]{qin2019bases}. They will only be used in Lemma \ref{lem:X_set_univ_Laurent}.

Take $t'=\seq t$ for any mutation sequence $\seq$. The mutation map $(\mu^\bullet)^*\colon \cF^\bullet(t')\rightarrow \cF^\bullet(t)$ induces an algebra homomorphism $\iota^\bullet(t')$ from $\LP^\bullet(t')$ to \smash{$\hLP^\bullet(t)$}, sending the Laurent monomials $(X')^n$ to $\big(\seq^X\big)^*(X')^n$ \big(or $(A')^m$ to $\big(\mu^A\big)^* (A')^m$\big). It is injective, see \cite[Lemma 3.3.7\,(1)]{qin2019bases}. In addition, the image of the Laurent monomials are pointed elements in the sense of Section \ref{sec:pointedness}, see \cite[Lemma 3.3.6]{qin2019bases}.
\begin{Lemma}
 For any $Z\in \LP^\bullet(t')$, the following statements are true.
\begin{itemize}\itemsep=0pt
\item[$(1)$] If $(\mu^\bullet)^* Z\in \LP^\bullet(t)$, then $\iota(Z)=(\mu^\bullet)^* Z$.
\item[$(2)$] If $\iota Z\in \LP^\bullet(t)$, then $\iota(Z)=(\mu^\bullet)^* Z$.
\end{itemize}
\end{Lemma}
\begin{proof}
 (1) See \cite[Lemma 3.3.7\,(2)]{qin2019bases}.

 (2) We prove it for the $X$-side, and the proof for the $A$-side is the same. Our proof is similar to that of \cite[Lemma 3.3.7\,(2)]{qin2019bases}.

 We can write $(X')^n * Z= F$ for some $F\in \kk[X'_i]_{i\in I}$, $n\in \N^I$. On the one hand, we have $\iota(X')^n* \iota(Z)=\iota(F)$ in~$\LP^X(s)$. On the other hand, we have $\seq^*(X')^n* \seq^*(Z)=\seq^*(F)$ in~$\cF^X(s)$. By the definition of $\iota$, we have $\iota(X')^n=\seq^*(X')^n$ and $\iota(F)=\seq^*(F)$. It follows that $\iota(Z)=\seq^*(Z)$.
\end{proof}

\subsection{Cluster algebras}

Let there be any given initial seed $t_{0}$. We use $\Delta^{+}=\Delta_{t_{0}}^{+}$
to denote the set of seeds $t=((A_{i}(t))_{i\in I},(X_{i}(t))_{i\in I},B(t))$
obtained from $t_{0}$ by sequences of mutations. If we work with the quantum cases, we also choose a compatible Poisson structure and consider the associated quantization as in Section~\ref{sec:quantization}. Recall that the
cluster variables $A_{j}(t_{0})$, $j\in I_{\fv}$, are unchanged
by mutations, which are denoted by $A_{j}$ and are called the frozen
variables. We use $\fFac$ to denote the multiplicative group generated
by $A_{j}^{\pm}$, $j\in\Ifv$.

Using the $v$-twisted product $*$, the partially compactified (quantum) cluster algebra $\bClAlg(t_{0})$ with
the initial seed $t_{0}$ is defined to be the $\kk$-subalgebra of
$\cF^{A}(t_{0})$ generated by all the cluster variables $\seq^{*}A_{i}(t)$,
$i\in I$, $t=\seq t_{0}\in\Delta_{t_{0}}^{+}$. The (localized) cluster
algebra $\clAlg(t_{0})$ is defined to be its localization $\bClAlg(t_{0})\big[A_{j}^{-1}\big]_{j\in I_{\fv}}$.

The upper cluster algebra (or upper cluster $A$-algebra) $\upClAlg(t_{0})$
with the initial seed $t_{0}$ is defined as the intersection $\cap_{t=\seq t_{0}\in\Delta^{+}}\seq^{*}\LP^{A}(t)$
inside $\cF^{A}(t_{0})$. By the Laurent phenomenon~\cite{BerensteinZelevinsky05,FominZelevinsky02},
it contains the cluster algebra $\clAlg(t_{0})$. For classical cases, if a compatible
Poisson structure on~$\LP^{A}(t_{0})$ is given, then $\upClAlg(t_{0})$
inherits the Poisson structure.

The cluster Poisson algebra (or upper cluster $X$-algebra) $\XClAlg(t_{0})$
with the initial seed $t_{0}$ is defined as the intersection $\cap_{t=\seq t_{0}\in\Delta^{+}}\seq^{*}\LP^{X}(t)$
inside $\cF^{X}(t_{0})$. For classical cases, it inherits the canonical Poisson structure
from that of $\LP^{X}(t_{0})$.

We often identify fraction fields $\cF^{\bullet}(t)$ and $\cF^{\bullet}(t_{0})$
via the mutation map $\seq^{*}$ for simplicity. Correspondingly,
we omit the symbol $t_{0}$ and $\seq^{*}$ in the above notations,
and we can write~${\upClAlg=\cap_{t}\LP^{A}(t)}$ and $\XClAlg=\cap_{t}\LP^{X}(t)$.

\subsection{Cluster varieties}

Let us work at the classical cases. Given two seeds $t'=\seq t$. Using the mutation birational maps $\seq\colon\torus^{A}(t)\dashrightarrow\torus^{A}(t')$,
we can glue the tori $\cT^{A}(t)$, $t\in\Delta^{+}$, into a scheme
$\AVar$ (see \cite[Proposition~2.4]{gross2013birational}), which
is called the cluster $A$-variety or cluster K2 variety. Note that
the upper cluster algebra $\upClAlg$ is the ring of global sections
of its structure sheaf. A choice of compatible Poisson structure on
$\LP^{A}(t_{0})$ gives rise to a Poisson structure on $\AVar$. It
often happens that $\AVar$ is a smooth manifold, for example, for
many well-known the cluster algebra arising from Lie theory (unipotent
cells \cite{GeissLeclercSchroeer10}, double Bruhat cells \cite{BerensteinFominZelevinsky05,goodearl2016berenstein}).

Similarly, using the mutation birational maps $\seq\colon\cT^{X}(t)\dashrightarrow\cT^{X}(t')$,
we can glue the tori $\cT^{X}(t)$, $t\in\Delta^{+}$, into a scheme
$\XVar$ called the cluster $X$-variety or the cluster Poisson variety. The cluster Poisson algebra
$\XClAlg$ is the ring of global sections of its structure sheaf.
Note that $\XVar$ has the canonical Poisson structure.

\subsection{Transition matrices}

We use $(\ )^{\rm T}$ to denote matrix transpose.

Given seeds $t'=\mu_{k}t$ and a mutation sign $\varepsilon\in\{+,-\}$.
Let us describe the monomial part of mutation using transition matrices.\footnote{In \cite{BerensteinZelevinsky05}, our matrices $\trans_{\varepsilon}^{N}$
and $\trans_{\varepsilon}^{M}$ are denoted by $F_{\varepsilon}$
and $E_{\varepsilon}$, respectively.}

Define the following $I\times I$-matrix $\trans_{k,\varepsilon}^{N}(t)$:
\begin{gather*}
\bigl(\trans_{k,\varepsilon}^{N}(t)\bigr)_{ij} = \begin{cases}
-1, & i=j=k,\\{}
[\varepsilon b_{kj}]_{+}, & i=k,\ j\neq k,\\
\delta_{ij}, & \mathrm{else}.
\end{cases}
\end{gather*}
Then it represents the linear map
\[\psi_{k,\varepsilon}^{N}\colon\ N(t')\rightarrow N(t),\qquad \psi_{k,\varepsilon}^{N}(\uee') =\uee\trans_{k,\varepsilon}^{N}(t).\]

We see that $\psi_{k,\varepsilon}^{N}$ induces the monomial part
$\psi_{k,\varepsilon}^{X}$ of the mutation from $\cF^{X}(t')$ to
$\cF^{X}(t)$.

Similarly, define the following $I\times I$-matrix $\trans_{k,\varepsilon}^{M}(t)$:
\begin{gather*}
\bigl(\trans_{k,\varepsilon}^{M}(t)\bigr)_{ij} = \begin{cases}
-1, & i=j=k,\\{}
[-\varepsilon b_{ik}]_{+}, & i\neq k,\ j=k,\\
\delta_{ij}, & \mathrm{else}.
\end{cases}
\end{gather*}
Then it represents the linear map
\[\psi_{k,\varepsilon}^{M}\colon\ \Mc(t')\simeq\Mc(t),\qquad \psi_{k,\varepsilon}^{M}(\uff') =\uff\trans_{k,\varepsilon}^{M}(t).\]
We see that $\psi_{k,\varepsilon}^{M}$ induces the monomial part
$\psi_{k,\varepsilon}^{A}$ of the mutation from~$\cF^{A}(t')$ to~$\cF^{A}(t)$.

The following results were known by \cite{BerensteinZelevinsky05,nakanishi2012tropical},
see \cite[Section~5.6]{Keller12} for a summary.

\begin{Proposition}\label{prop:matrix_properties} We have the following
equalities:
\begin{enumerate}\itemsep=0pt
 \item[$(1)$] $\trans_{k,\varepsilon}^{N}(t)^{2}=\Id$, $\trans_{k,-\varepsilon}^{N}(t')=\trans_{k,\varepsilon}^{N}(t)$,
 $\trans_{k,-\varepsilon}^{N}(t')\trans_{k,\varepsilon}^{N}(t)=\Id$.

 \item[$(2)$] $\trans_{k,\varepsilon}^{M}(t)^{2}=\Id$, $\trans_{k,-\varepsilon}^{M}(t')=\trans_{k,\varepsilon}^{M}(t)$,
 $\trans_{k,-\varepsilon}^{M}(t')\trans_{k,\varepsilon}^{M}(t)=\Id$.

 \item[$(3)$] $D\trans_{k,-\varepsilon}^{M}(t')D^{-1}=D\trans_{k,\varepsilon}^{M}(t)D^{-1}=\bigl(\trans_{k,\varepsilon}^{N}(t)\bigr)^{\rm -T}$.

 \item[$(4)$] $B'=\trans_{k,\varepsilon}^{M}(t)\cdot B\cdot\trans_{k,\varepsilon}^{N}(t)$.

 \item[$(5)$] $\Lambda'=\trans_{k,\varepsilon}^{M}(t)^{\rm T}\cdot\Lambda\cdot\trans_{k,\varepsilon}^{M}(t)$, where $\Lambda$ is the matrix of the Poisson structure.
\end{enumerate}
\end{Proposition}

\begin{Remark}\label{rem:signed_mutation_non_identity}
 It is straightforward to check that we have
 \[\bigl(\trans^M_{k,\varepsilon}(t')\trans^M_{k,\varepsilon}(t)\bigr)_{ij}=\begin{cases}
 -\varepsilon b_{ik}, & i\neq k,\ j=k,\\
 \delta_{ij}, & \mathrm{else}.
 \end{cases}\]
 and similarly
 \[
 \bigl(\trans^N_{k,\varepsilon}(t')\trans^N_{k,\varepsilon}(t)\bigr)_{ij} = \begin{cases}
 -\varepsilon b_{kj}, & i=k,\ j\neq k,\\
 \delta_{ij}, & \mathrm{else}.
 \end{cases}
 \]
\end{Remark}

\subsection[Cluster expansions, g-vectors, c-vectors]{Cluster expansions, $\boldsymbol{g}$-vectors, $\boldsymbol{c}$-vectors}

Let there be given an initial seed $t_{0}=((A_{i}),(X_{i}),B)$. Take
any sequence of unfrozen vertices $\uk=(k_{0},k_{1},\ldots,k_{r})$.
Denote the corresponding mutation sequence by $\seq=\seq_{\uk}=\mu_{k_{r}}\cdots\mu_{k_{0}}$
and the resulting seed by $t=\seq t_{0}$.

Recall that a vector is said to be sign-coherent if its coordinates
are all non-negative or all non-positive.

\begin{Theorem}[{\cite{DerksenWeymanZelevinsky09,gross2018canonical,Tran09}}]\label{thm:cluster_expansion}
There exist $I\times I$ invertible $\Z$-matrices \[E(t)=\left(\begin{matrix}
E(t)_{\ufv} & E(t)_{\high}\\
0 & \Id_{\fv}
\end{matrix}\right)\qquad \text{and}\qquad F(t)=\left(\begin{matrix}
F(t)_{\ufv} & 0\\
F(t)_{\low} & \Id_{\fv}
\end{matrix}\right),\] such that the following statements hold.
\begin{itemize}\itemsep=0pt
\item[$(1)$] Any cluster $A$-variable $A_{i}(t)$, $i\in I$, has the following
Laurent expansion in $\LP^{A}(t_{0})$:
\begin{align*}
\seq^{*}A_{i}(t) & =A^{\col_{i}F(t)}*\sum_{n\in\N^{I_{\ufv}}}c_{n}A^{p^{*}n},
\end{align*}
such that $c_{0}=1$.

\item[$(2)$] Any cluster $X$-variable $X_{i}(t)$, $i\in I$, has the following
expression in $\cF^{X}(t_{0})$:
\begin{align*}
\seq^{*}X_{i}(t) & =X^{\col_{i}E(t)}*P*Q^{-1},
\end{align*}
where $P$, $Q$ are polynomials in $\kk[X_{k}]_{k\in I_{\ufv}}$ with
constant term $1$.

\item[$(3)$] The row vectors of $F(t)$ are sign-coherent.

\item[$(4)$] The column vectors of $E(t)$ are sign-coherent.
\end{itemize}
\end{Theorem}

The $I_{\ufv}\times I_{\ufv}$-submatrices $C(t):=E(t)_{\ufv}$ and
$G(t):=F(t)_{\ufv}$ are usually called the $C$-matrix and the $G$-matrix
of $t$ (with respect to the initial seed $t_{0}$), respectively.

\begin{Definition}[\cite{qin2019bases}]
For any pair of seeds $t,t'\in\Delta^{+}$, let $t$ be the
initial seed and $E(t')$ (resp.~$F(t')$) denote the corresponding
$E$-matrix (resp.\ $F$-matrix) of $t'$. We define the linear map:
\begin{alignat*}{3}
& \psi_{t,t'}^{N}\colon\ N(t')\rightarrow N(t),\qquad&& \psi_{t,t'}^{N}(\uee(t'))=\uee(t)E(t');&\\
& \psi_{t,t'}^{M}\colon\ \Mc(t')\rightarrow\Mc(t),\qquad&& \psi_{t,t'}^{M}(\uff(t'))=\uff(t)F(t').&
\end{alignat*}
\end{Definition}

\begin{Example}\label{eg:A_1_eg}
Take the index set $I=\{1,2\}$ such that $1$ is the only unfrozen
vertex. Choose $d_{1}=2$, $d_{2}=1$. The initial seed $t=((X_{1},X_{2}),B)$
is given such that $B=\left(\begin{smallmatrix}
\hphantom{-}0 & 2\\
-1 & 0
\end{smallmatrix}\right)$.

Then \smash{$W=B^T \left(\begin{smallmatrix}
 \frac{1}{2} & 0\\
 0 & 1
 \end{smallmatrix}\right)=\left(\begin{smallmatrix}
 0 & -1\\
 1 & \hphantom{-}0
 \end{smallmatrix}\right)$}. We choose $\Lambda=-\frac{1}{2}W$, so that $(\Lambda B)_{i1}=-\delta_{i,1}\frac{1}{d_1}$ for $i=1,2$.

The seed $t'=\mu_{1}t$ is the only non-initial seed. Using the commutative product $\cdot$, we can write $A_{1}'=A^{-\ff_1+\ff_2}+A^{-\ff_1}=A_{1}^{-1}\cdot A_{2}\cdot \big(1+A_{2}^{-1}\big)=A_{1}^{-1}\cdot A_{2}\cdot (1+p^{*}X_{1})$
and $A_{2}'=A_{2}.$ It follows that $\psi_{t,t'}^{M}\colon \Mc(t')\rightarrow\Mc(t)$
is represented by $F(t')=\left(\begin{smallmatrix}
-1 & 0\\
\hphantom{-}1 & 1
\end{smallmatrix}\right)$. We also have \[X_{1}'=X_{1}^{-1}, \qquad (X_{2}')^{-1}=X_{2}^{-1}\cdot X_{1}^{-2}\cdot \big(1+\big(v^{\frac{1}{2}}+v^{-\frac{1}{2}}\big)X_1+X_{1}^2\big)\] or, equivalently, \[X_2'=\big(X^{-2\ee_1-\ee_2}+\big(v^{\frac{1}{2}}+v^{-\frac{1}{2}}\big)X^{-\ee_1-\ee_2}+X^{-\ee_2}\big)^{-1}.\]
So $\psi_{t,t'}^{N}\colon N(t')\rightarrow N(t)$ is represented by $E(t')=\left(\begin{smallmatrix}
-1 & 2\\
\hphantom{-}0 & 1
\end{smallmatrix}\right)$.

\end{Example}

Note that, for any $k\in\Iufv$, $\psi_{t,\mu_{k}t}^{\bullet}$ is
represented by the matrix $\trans_{k,+}^{\bullet}(t)$, where $\bullet$
stands for $M$ or $N$. In particular, $\psi_{t,\mu_{k}t}^{\bullet}\psi_{\mu_{k}t,t}^{\bullet}$
might not be the identity by Remark \ref{rem:signed_mutation_non_identity}.

Identify $N(t)$ with $N(t_{0})$ by using the linear map $\psi_{t_{0},t}^{N}$.
Then the basis vector $\ee_{i}(t)$ has the coordinate vector $\col_{i}E(t)$
in $N(t_{0})\simeq\Z^{I}$. Similarly, identify $M(t)$ with $M(t_{0})$
by using the linear map $\psi_{t_{0},t}^{M}$. Then the basis vector
$\ff_{i}(t)$ has the coordinate vector $\col_{i}F(t)$ in $\Mc(t_{0})\simeq\Z^{I}$.

We will often work with vectors, linear maps and bilinear forms in
the fixed lattices $N(t_{0})$ and $M(t_{0})$. We refer the reader
to \cite{gross2013birational} for more details on the fixed data.

\subsection{Canonical mutation signs}

For any $k\in I_{\ufv}$, recall that the $k$-th $c$-vector $\col_{k}C(t)$
is sign coherent. Then we define the canonical mutation sign $\varepsilon$
for the mutation of $t$ at the vertex $k$ to be the sign of $\col_{k}C(t)$.
From now on, we always choose the canonical mutation sign unless otherwise
specified.

For $\uk=(k_{0},k_{1},\ldots,k_{r})$, denote $\uk_{\leq s}=(k_{0},\ldots,k_{s})$
for any $0\leq s\leq r-1$ and $t_{s}=\seq_{\uk_{\leq s-1}}t_{0}$.
Let $\varepsilon_{s}$ denote the canonical sign for the mutation
of $t_{s}$ at the direction $k_{s}$.

\begin{Proposition}[{\cite{gross2018canonical,nakanishi2012tropical}}]\label{prop:recursive_bases}
We have the following
\begin{gather*}
E(t) =\trans_{k_0,\varepsilon_{0}}^{N}(t_{0})\trans_{k_1,\varepsilon_{1}}^{N}(t_{1})\cdots\trans_{k_{r-1},\varepsilon_{r-1}}^{N}(t_{r-1}),\\
F(t) =\trans_{k_0,\varepsilon_{0}}^{M}(t_{0})\trans_{k_1,\varepsilon_{1}}^{M}(t_{1})\cdots\trans_{k_{r-1},\varepsilon_{r-1}}^{M}(t_{r-1}).
\end{gather*}
\end{Proposition}
By Proposition \ref{prop:recursive_bases}, the matrices $C(t)$ and $G(t)$ only depend on $B_{\ufv}$
and $\uk$.

Using Proposition \ref{prop:matrix_properties}, we obtain
\begin{align}
E(t)^{\rm T}=DF(t)^{-1}D^{-1} &.\label{eq:E_D_F_relation}
\end{align}
Its restriction gives $C(t)^{\rm T}=D_{\ufv}G(t)^{-1}D_{\ufv}^{-1}$.

\begin{Lemma}\label{lem:bilinear_form_basis_change}

The following statements are true.
\begin{itemize}\itemsep=0pt
\item[$(1)$] $F(t)B(t)(E(t))^{-1}=B(t_{0})$ or, equivalently, $\psi_{t_{0},t}^{M}(p^{*}\uee(t))=p^{*}\big(\psi_{t_{0},t}^{N}\uee(t)\big)$.
\item[$(2)$] $F(t)^{\rm T}\Lambda(t_{0})F(t)=\Lambda(t)$ or, equivalently, $\lambda\big(\psi_{t_{0},t}^{M}(\uff(t)),\psi_{t_{0},t}^{M}(\uff(t))\big)=\lambda(\uff(t),\uff(t))$.
\item[$(3)$] $E(t)^{\rm T}W(t_{0})E(t)=W(t)$ or, equivalently, $\omega\big(\psi_{t_{0},t}^{N}(\uee(t)),\psi_{t_{0},t}^{N}(\uee(t))\big)=\omega(\uee(t),\uee(t)).$
\end{itemize}
\end{Lemma}

\begin{proof}

(1)--(2) follow from Propositions \ref{prop:matrix_properties} and
\ref{prop:recursive_bases}.

(3) follow from $(1)$, $W(t)D^{-1}=B(t)^{\rm T}$ and $E(t)^{\rm T}=DF(t)^{-1}D^{-1}$. \hfill $\qed$
 \renewcommand{\qed}{}
\end{proof}

\subsection{Principal coefficients\label{subsec:Principal-coefficients}}

Let us recall the seeds with principal coefficients and their relation
with the $C$-matrix and the $G$-matrix. This part will only be used in Sections
\ref{subsec:Cluster-twist-automorphism-principal} and \ref{sec:construction_X_bases}.

Denote a copy of $\Iufv$ by $\Iufv'=\{k'|k'\in\Iufv\}$. We extend
the corresponding principal $B$-matrix $B_{\ufv}$ to the $(\Iufv\sqcup\Iufv')\times(\Iufv\sqcup\Iufv')$
matrix \smash{$B^{\prin}=\left(\begin{smallmatrix}
B_{\ufv} & -\Id_{\Iufv}\\
\Id_{\Iufv} & 0
\end{smallmatrix}\right)$}, which is called the principal coefficient $B$-matrix, where $\Id_{\Iufv}$
represents the natural isomorphism $\Iufv\simeq\Iufv'$. The corresponding
diagonal matrices are \[D^{\prin}=\left(\begin{matrix}
D_{\ufv}\\
 & D_{\ufv}
\end{matrix}\right)\qquad \text{and}\qquad W^{\prin}=\left(\begin{matrix}
W_{\ufv} & D_{\ufv}\\
-D_{\ufv} & 0
\end{matrix}\right).\] We denote by $t_{0}^{\prin}$ the seed obtained from $t_{0}$ by
changing the fixed data as above.

Let $\uk$ denote a sequence of vertices and $t=\seq_{\uk}t_{0}$.
Then it is known that the $C$-matrix and the $G$-matrix can be computed
using principle coefficients:
\begin{align*}
F\big(\seq_{\uk}t_{0}^{\prin}\big) & =\left(\begin{matrix}
G(t) & 0\\
0 & \Id_{\Iufv}
\end{matrix}\right),\\
E\big(\seq_{\uk}t_{0}^{\prin}\big) & =\big(D^{\prin}\big)^{-1}F(\seq_{\uk}t_{0})^{\rm -T}D^{\prin}=\left(\begin{matrix}
C(t) & 0\\
0 & \Id_{\Iufv}
\end{matrix}\right),
\end{align*}
where $C(t)=D_{\ufv}^{-1}G(t)^{\rm -T}D_{\ufv}$. Then, using Lemma \ref{lem:bilinear_form_basis_change},
we obtain
\begin{align*}
B\big(\seq_{\uk}t_{0}^{\prin}\big) & =F\big(\seq_{\uk}t_{0}^{\prin}\big)^{-1}B\big(t_{0}^{\prin}\big)E(\seq_{\uk}t_{0}) =\left(\begin{matrix}
G(t)^{-1}B_{\ufv}C(t) & -G(t)^{-1}\\
C(t) & 0
\end{matrix}\right)\\
 & =\left(\begin{matrix}
B(t)_{\ufv} & -G(t)^{-1}\\
C(t) & 0
\end{matrix}\right)=\left(\begin{matrix}
B(t)_{\ufv} & -D_{\ufv}^{-1}C(t)^{\rm T}D_{\ufv}\\
C(t) & 0
\end{matrix}\right).
\end{align*}

\subsection{Degrees and pointedness}\label{sec:pointedness}

We endow $\Nufv(t)=\oplus\Z \ee_{k}\simeq\Z^{\Iufv}$ with the natural
partial order such that $n\geq n'$ if $n-n'\geq0$. Denote $\Nufv^{\geq0}(t)=\oplus\N \ee_{k}\simeq\N^{\Iufv}$.
Recall that we have the linear map $p^{*}\colon\Nufv(t)\rightarrow\Mc(t)$
represented by the matrix \smash{$\tB=\left(\begin{smallmatrix}
B_{\ufv}\\
B_{\low}
\end{smallmatrix}\right)$}.

\begin{Definition}[pointedness]
A formal Laurent series $Z\in\hLP^{X}(t)$ is said to have degree $n$,
denoted by $\deg_{t}Z=n$, if it takes the form \smash{$Z=X(t)^{n}\cdot\big(\sum_{n'\in\Nufv^{\geq0}(t)}c_{n'}X(t)^{n'}\big)$}
for some $c_{0}\neq0$, $c_{n'}\in\kk$. Its $F$-function is defined as \smash{$\sum_{n'\in\Nufv^{\geq0}(t)}c_{n'}X(t)^{n'}$}. It is further said to be pointed at $n$ or $n$-pointed if $c_{0}=1$.
\end{Definition}

Now assume that $t$ is a seed such that $\ker p^{*}\cap\Nufv^{\geq0}(t)=0$.
This condition is satisfied when the seed $t$ satisfies the full
rank assumption (see Assumption~\ref{fullrank}). We recall the degrees and pointedness introduced
in \cite{qin2017triangular,qin2019bases}.

\begin{Definition}[dominance order \cite{qin2017triangular}]
For any $m,m'\in\Mc(t)$, we say $m'$ is dominated by~$m$, denoted
by $m'\preceq_{t}m$, if $m'=m+p^{*}n$ for some $n\in\Nufv^{\geq0}(t)$.
\end{Definition}

\begin{Definition}[pointedness \cite{qin2017triangular}]\label{def:A-pointedness}
A formal Laurent series $Z\in\hLP^{A}(t)$ is said to have degree
$m$, denoted by $\deg_{t}Z=m$, if it takes the form \smash{$Z=A(t)^{m}\cdot\big(\sum_{n''\in\Nufv^{\geq0}(t)}c_{n}A(t)^{p^{*}n''}\big)$}
for some $c_{0}\neq0$, $c_{n''}\in\kk$, i.e., its $\prec_{t}$-maximal
Laurent degree is $m$. Its $F$-function is defined as $\sum_{n''\in\Nufv^{\geq0}(t)}c_{n''}X(t)^{n''}$. It is further said to be pointed at $m$ or
$m$-pointed if $c_{0}=1$.

\end{Definition}

Theorem \ref{thm:cluster_expansion} implies that $\deg_{t_{0}}X_{i}(t)=\col_{i}E(t)$
and $\deg_{t_{0}}A_{i}(t)=\col_{i}F(t)$. By definition, $\psi_{t_{0},t}^{N}$
is the linear map sending $\ee_{i}(t)=\deg_{t}X_{i}(t)$ to $\deg_{t_{0}}\seq^{*}X_{i}(t)$. Similarly, $\psi_{t_{0},t}^{M}$ is the linear map sending $\ff_{i}(t)=\deg_{t}A_{i}(t)$
to $\deg_{t_{0}}\seq^{*}A_{i}(t)$.

\subsection{Injective-reachable}

\begin{Definition}[{injective-reachable \cite[Section 2.3]{qin2017triangular}}]\label{def:Injective-reachable}
A seed $t_{0}$ is said to be injective-reach\-able, if there exists
another seed $t_{0}[1]=\seq t_{0}\in\Delta_{t_{0}}^{+}$ and a permutation
$\sigma$ of $I_{\ufv}$, such that for any~$k\in\Iufv$, we have
\begin{align}
\psi_{t_{0},t_{0}[1]}^{N}\ee_{\sigma k}(t_{0}[1])=-\ee_{k}(t_{0}).\label{eq:injective_reachable_def}
\end{align}
\end{Definition}

The sequence $\seq$ is also called a green to red sequence from $t_{0}$
to $t_{0}[1]$ in the sense of Keller~\cite{keller2011cluster}.

Assume that $t_{0}$ is injective-reachable. Then all $t\in\Delta_{t_{0}}^{+}$
are injective-reachable, i.e., we can always find a seed $t[1]$, see
\cite[Proposition~5.1.4]{qin2017triangular} or \cite{muller2015existence}.

By \cite[Proposition~2.3.3]{qin2017triangular}, for any $k\in\Iufv$,
we have \[\psi_{t_{0},t_{0}[1]}^{N}\ff_{\sigma k}(t_{0}[1])\in-\ff_{k_{0}}(t_{0})+\oplus_{j\in\Ifv}\Z\ff_{j},\qquad d_{k}=d_{\sigma k}.\]
 We will see that $t_{0}[1]\in\Delta^{+}$ is similar
to $t_{0}$ up to $\sigma$ in the sense of Definition \ref{def:similar_seeds}.

\section[Twist endomorphisms for upper cluster A-algebras]{Twist endomorphisms for upper cluster $\boldsymbol{A}$-algebras}\label{sec:Twist-endomorphisms-A}

In this section, we introduce the notion of a twist endomorphisms
for a pair of similar seeds $t$, $t'$, which is defined as the composition
of a mutation map with a monomial map called a variation map.

\subsection{Similar seeds and variation maps}

We first recall the definition of similar seeds.

\begin{Definition}[{similar seeds \cite{Qin12,qin2017triangular}}]\label{def:similar_seeds}
Two seeds $t$, $t'$ are called similar, if there is a permutation~$\sigma$ of $\Iufv$ such that for any $i,j\in\Iufv$, we have $b_{ij}(t)=b_{\sigma i,\sigma j}(t')$
and skew-symmetrizers $d_{i}=d_{\sigma i}$, $\forall i,j\in I_{\ufv}$.
\end{Definition}
If $t$ and $t'$ are similar, in our choices of compatible Poisson structures, we automatically have $\diag_{k}(t)=\diag_{\sigma k}(t')$ for $k\in\Iufv$
(see Definition~\ref{def:compatible_poisson_A}).

Note that we can always trivially extends $\sigma$ to a permutation
on $I$.

\begin{Proposition}
Assume that $B_{\ufv}$ is connected Definition~{\rm \ref{def:connected_matrix}}.
Let $t=\seq t_{0}$ be a given seed. If there is a permutation $\sigma$
of $\Iufv$ such that $b_{ij}(t)=b_{\sigma i,\sigma j}(t_0)$ for any
$i,j\in\Iufv$, then $d_{k}=d_{\sigma k}$ for any $k\in\Iufv$.
\end{Proposition}

\begin{proof}
Let $P_{\sigma}$ be the permutation matrix associated to $\sigma$ of rank $\Iufv$, i.e., the $i$-th column of~$P_{\sigma}$ is the $\sigma i$-th unit vector. Notice that $P_{\sigma}^{\rm T}=P_{\sigma}^{-1}=P_{\sigma^{-1}}$.

Rewrite the equality $b_{\sigma i,\sigma j}(t)=b_{ij}(t_{0})$ as $P_{\sigma}^{-1}B(t)_{\ufv}P_{\sigma}=B(t_0)_{\ufv}$.
Since $D_{\ufv}B(t_0)_{\ufv}$ is skew-symmetric,
\[\big(P_{\sigma}D_{\ufv}P_{\sigma}^{-1}\big)B(t)_{\ufv}=P_{\sigma}D_{\ufv}B(t_0)_{\ufv}P_{\sigma}^{-1}\]
is skew-symmetric.

On the other hand, we have $B(t)=F(t)^{-1}B(t_{0})E(t)$ by Lemma~\ref{lem:bilinear_form_basis_change}. Using $DF(t)^{-1}=E(t)^{\rm T}D$
by \eqref{eq:E_D_F_relation}, we obtain
\[DB(t)=E(t)^{\rm T}DB(t_{0})E(t).\]
Since $DB(t_{0})$ is skew-symmetric, so is $DB(t)$. It follows that
$D_{\ufv}B(t)_{\ufv}$ is also skew-symmetric.

Then Lemma \ref{lem:connected_symmetrizer} implies that $P_{\sigma}^{-1}D_{\ufv}P_{\sigma}=\alpha D_{\ufv}$
for some $\alpha\neq0$. Since $D_{\ufv}$ is of full rank,
by taking the determinant, we see that $\alpha=1$. The claim follows.
\end{proof}

Assume $t$, $t'$ are similar up to a permutation $\sigma$. In view
of the correction technique \cite{Qin12}, it is often useful to compare
pointed elements in $\LP^A(t)$ with those in $\LP^A(t')$. As in \cite{qin2017triangular},
we define a variation map sending pointed elements in $\LP^A(t)$ to
those in $\LP^A(t')$, which differ only by the frozen variables.

In this paper, we will allow roots of frozen variables. For any given
integer number $r>0$, we define \smash{$\fFac^{\frac{1}{r}}$} to be the
multiplicative group generated by the $r$-th roots \smash{$A_{j}^{\pm\frac{1}{r}}$},
$j\in I_{\fv}$. Define $\fFacRoot$ to be the multiplicative group
generated by any \smash{$A_{j}^{\pm\frac{1}{r}}$}, $j\in I_{\fv}$, $r>0$
.

For any $\kk[\fFac]$-algebra $R$, we use $R_{\fFacRoot}$ to denote
$R\otimes_{\kk[\fFac]}\kk[\fFacRoot]$. Define $R_{\fFac^{\frac{1}{r}}}$
similarly.

\begin{Definition}\label{def:variation_homomorphism} Let $t$, $t'$ be two seeds similar up to a permutation $\sigma$. A $\kk$-algebra homomorphism $\var_{t}^{A}\colon \LP(t)_{\fFacRoot}\rightarrow\LP(t')_{\fFacRoot}$
is called a (monomial) variation map from $t$ to $t'$ if, $\forall k\in\Iufv$,
$j\in\Ifv$, it takes the form
\begin{gather*}
\var_{t}^{A}(A_{k}) = p_{\sigma k}\cdot A_{\sigma k}',\qquad
\var_{t}^{A}\big(A^{\col_{k}B}\big) = (A')^{\col_{\sigma k}B'},\qquad
\var_{t}^{A}(A_{j}) = p_{j},
\end{gather*}
for some $p_{\sigma k},p_{j}\in\fFacRoot$.

For classical cases, if $t$ and $t'$ are similar seeds with compatible Poisson structures, then
$\var_{t}^{A}$ is called a Poisson variation map if it further preserves
the Poisson structures:
\begin{align*}
\var_{t}^{A}\{A_{i},A_{j}\} & =\bigl\{\var_{t}^{A}A_{i},\var_{t}^{A}A_{j}\bigr\},\qquad \forall i,j\in I.
\end{align*}

\end{Definition}

\begin{Definition}
 In $\hLP^A(t)$, take any pointed formal Laurent series \smash{$Z=A^m \cdot F|_{X^{n}\mapsto A^{p^*n}}$}, where $m\in \Z^I$, \smash{$F=\sum_{n\in \N_{\ufv}^{\geq 0}} c_n X^n$}, $c_n\in \kk$, $c_0= 1$. We say $Z$ and \smash{$Z'\in \hLP^A(t')_{\fFacRoot}$} are similar if $Z'=(A')^{m'}\cdot F|_{X^{n}\mapsto (A')^{p^*(t')(\sigma n)}}$ with $\pr_{\Iufv}m'=\sigma \pr_{\Iufv}m$.
\end{Definition}

\begin{Remark}\label{rem:variation_map_pointed_elements}
 The variation map $\var_t^A$ is a $\kk$-algebra homomorphism such that it sends a pointed element $Z\in\LP^A(t)$ to a similar element in $\LP^A(t')_{\fFacRoot}$.
\end{Remark}

The variation map in Definition \ref{def:variation_homomorphism}
is the monomial map associated to the following linear variation map.

Denote $\Mc_{\ufv}(t)=\oplus_{k\in\Iufv}\Z\ff_{k}$, $\Mc_{\fv}(t)=\oplus_{j\in\Ifv}\Z\ff_{j}$.

\begin{Definition}\label{def:variation_linear}
Let $t$, $t'$ be two seeds similar up to a permutation $\sigma$.
A linear map $\var_{t}^{M}\colon\Mc(t)_{\Q}\allowbreak\rightarrow\Mc(t')_{\Q}$ is
called a (linear) variation map from $t$ to $t'$, for any $k\in\Iufv$,
$j\in\Ifv$, if it takes the following form
\begin{gather}
\var_{t}^{M}(\ff_{k}) = \ff'_{\sigma k}+\uu_{\sigma k},\qquad
\var_{t}^{M}(\ff_{j}) = \uu_{j},\qquad
\var_{t}^{M}\Bigl(\sum b_{ik}\ff_{i}\Bigr) = \sum b'_{i,\sigma k}\ff'_{i},\label{eq:preserve_p_matrix}
\end{gather}
for some $\uu_{k},\uu_{j}\in(\Mc_{\fv}(t'))_{\Q}$.

If $t$ and $t'$ are similar seeds with compatible Poisson structures $\lambda$, then
it is called a Poisson variation map if we have
\begin{gather*}
\lambda(\ff_{i},\ff_{j}) = \lambda\bigl(\var_{t}^{M}\ff_{i},\var_{t}^{M}\ff_{j}\bigr),\qquad \forall i,j\in I.
\end{gather*}
\end{Definition}

Note that \eqref{eq:preserve_p_matrix} can be written as
\begin{equation}\label{eq:preserve_p_map}
\var_{t}^{M}(p^{*}\ee_{k}) =p^{*}\ee'_{\sigma k}.
\end{equation}

\begin{Remark}\label{rem:Poisson_quantum_variation}
Assume the existence of compatible Poisson structures, i.e, Assumption \ref{fullrank} holds. Then the following statements are equivalent:
\begin{itemize}\itemsep=0pt
 \item A linear map $\varPsi\colon\Mc(t)\rightarrow \Mc(t')$ is a Poisson variation map.
 \item For the quantum case $\kk=\Z\big[v^{\pm\frac{1}{\diag}}\big]$, the monomial map $\varPhi_v$ associated to $\varPsi$ is a variation map. In particular, it preserves the $v$-twisted products.
 \item For the classical case $\kk=\Z$, the monomial map $\varPhi_1$ associated to $\varPsi$ is a Poisson variation map.
\end{itemize}
Therefore, a Poisson variation map is the classical limit of a quantum variation map at $v=1$, see Section \ref{sec:quantization}. Conversely, a Poisson variation map gives rise to a quantum variation map by the above equivalent statements.
\end{Remark}

Let $P_{\sigma}$ denote the permutation matrix associated to $\sigma$
of rank $\Iufv$.

\begin{Lemma}\label{lem:matrix_eq_A_variation}\quad
\begin{itemize}\itemsep=0pt
\item[$(1)$] The variation map $\var_{t}^{M}$ has the following matrix representation
in the bases $\uff$, $\uff'$:
\begin{align*}
\var_{t}^{M}(\uff) & =\uff'\cdot\begin{pmatrix}
\Id_{\ufv} & 0\\
U_{\low} & U_{\fv}
\end{pmatrix}\begin{pmatrix}
\trans_{\sigma}\\
 & \Id_{\fv}
\end{pmatrix},
\end{align*}
where $\col_{k}U_{\low}$ $($resp.\ $\col_{j}U_{\fv})$ are the coordinates of $\uu_{k}$ $($resp.\ $\uu_{j})$ in the basis $\{\ff_{j'},\,j'\in\Ifv\}$, for $k\in \Iufv $
$($resp.\ $j\in I_{\fv})$.

\item[$(2)$] Moreover, equation \eqref{eq:preserve_p_matrix} is equivalent
to the following
\begin{gather}
\begin{pmatrix}
U_{\low}\trans_{\sigma} & U_{\fv}\end{pmatrix}\begin{pmatrix}
B_{\ufv}\\
B_{\low}
\end{pmatrix} = B'_{\low}\trans_{\sigma}.\label{eq:solve_cluster_A_twist}
\end{gather}
\end{itemize}
\end{Lemma}

\begin{proof}

(1) Recall that $\col_{k}(HP_{\sigma})=\col_{\sigma k}H$ for any
matrix $H$. The first statement follows.

(2) For any $k\in\Iufv$, we have
\begin{align*}
\var_{t}^{M}(p^{*}\ee_{k}) & =\var_{t}^{M}(\uff\col_{k}B) =\uff'\cdot\begin{pmatrix}
\Id_{\ufv} & 0\\
U_{\low} & U_{\fv}
\end{pmatrix}\begin{pmatrix}
\trans_{\sigma} & 0\\
0 & \Id_{\fv}
\end{pmatrix}\col_{k}B\\
 & =\uff'\cdot\begin{pmatrix}
\trans_{\sigma} & 0\\
U_{\low}\trans_{\sigma} & U_{\fv}
\end{pmatrix}\col_{k}\begin{pmatrix}
B_{\ufv}\\
B_{\low}
\end{pmatrix} =\uff'\cdot\col_{k}\left(\begin{pmatrix}
\trans_{\sigma} & 0\\
U_{\low}\trans_{\sigma} & U_{\fv}
\end{pmatrix}\begin{pmatrix}
B_{\ufv}\\
B_{\low}
\end{pmatrix}\right).
\end{align*}
We have
\begin{align*}
p^{*}\ee_{\sigma k}' & =\uff'\col_{\sigma k}B' =\uff'\col_{k}\left(B'\begin{pmatrix}
\trans_{\sigma} & 0\\
0 & \Id_{\fv}
\end{pmatrix}\right) =\uff'\col_{k}\begin{pmatrix}
B'_{\ufv}\trans_{\sigma}\\
B'_{\low}\trans_{\sigma}
\end{pmatrix}.
\end{align*}
Note that $\trans_{\sigma}^{-1}B'_{\ufv}\trans_{\sigma}=B_{\ufv}$
since $b'_{\sigma i,\sigma j}=b_{i,j}$ for $i,j\in\Iufv$. So $\trans_{\sigma}B_{\ufv}=B'_{\ufv}\trans_{\sigma}$.
The claim follows.
\end{proof}

\begin{Lemma}\label{lem:inverse_M_variation}

If a $($Poisson$)$ variation map $\var_{t}^{M}$ is invertible, its inverse
is still a $($Poisson$)$ variation map.

\end{Lemma}

\begin{proof}

Since $\var_{t}^{M}$ is represented by an invertible matrix
\[
\begin{pmatrix}
\Id_{\ufv} & 0\\
U_{\low} & U_{\fv}
\end{pmatrix}\begin{pmatrix}
\trans_{\sigma}\\
 & \Id_{\fv}
\end{pmatrix},\]
its inverse $\bigl(\var_{t}^{M}\bigr)^{-1}$ is represented by
\[\begin{pmatrix}
\trans_{\sigma^{-1}}\\
 & \Id_{\fv}
\end{pmatrix}\begin{pmatrix}
\Id_{\ufv} & 0\\
-U_{\fv}^{-1}U_{\low} & U_{\fv}^{-1}
\end{pmatrix}=
\begin{pmatrix}
\Id_{\ufv} & 0\\
-U_{\fv}^{-1}U_{\low}\trans_{\sigma} & U_{\fv}^{-1}
\end{pmatrix}\begin{pmatrix}
\trans_{\sigma^{-1}}\\
 & \Id_{\fv}
\end{pmatrix}.
\]
Thus $\bigl(\var_{t}^{M}\bigr)^{-1}$ is a variation map with permutation ${\sigma}^{-1}$.

Since $\var_{t}^{M}$ is a variation map, we have $\var_{t}^{M}p^{*}\ee_{k}=p^{*}\ee_{\sigma k}'$
for any $k\in\Iufv$, it follows that \smash{$\bigl(\var_{t}^{M}\bigr)^{-1}p^{*}\ee_{k}'=p^{*}\ee_{{\sigma}^{-1}k}$}.
So~\eqref{eq:preserve_p_map} holds. Finally, if $\var_{t}^{M}$ preserves
the compatible Poisson structures, so does its inverse.
\end{proof}

\begin{Lemma}\label{lem:A_variation_restrict_cluster_alg}
 A variation map $\var_{t}^{A}\colon\LP^A(t)\rightarrow\LP^A(t')$ restricts to an algebra homomorphism from $\upClAlg(t)$
 to $\upClAlg(t')$.
 \end{Lemma}

 \begin{proof}
 Take any mutation sequence $\seq=\mu_{k_{r}}\cdots\mu_{k_{1}}$, which is the identify if the sequence $(k_r,\ldots,k_1)$ is empty. Denote
 $\sigma\seq=\mu_{\sigma k_{r}}\cdots\mu_{\sigma k_{1}}$. Then $s=\seq t$
 and $s'=(\sigma\seq)t'$ are similar up to $\sigma$. To prove the claim, it suffices to show that, for any $z\in \LP^A(t)\cap \seq^*\LP^{A}(s)$, we have~${\var^A_t(z)\in (\sigma\seq)^*\LP^{A}(s')}$. If so, we obtain that, for any \begin{gather*}
 z\in \upClAlg(t)=\LP^A(t)\bigcap (\cap_{s=\seq t} \seq^*\LP^{A}(s)),\\
 \var^A_t(z)\in \cap_{s'=(\sigma \seq) t'} (\sigma\seq)^*\LP^{A}(s')=\upClAlg(t').\end{gather*}

 By the Laurent expansion of $z$ in the seed $s$, we have $z* \big(\seq^* A(s)^d\big)=\sum c_m (\seq^* A(s)^m)$ for some $c_m\in \kk$, $d,m\in \N^I$. By \cite[Lemma 4.2.2\,(ii)]{qin2017triangular}, the variation map $\var^A_t$ sends cluster monomials of $s$ to the cluster monomials of $s'$ up to frozen factors: we have ${\var^A_t\big(\seq^* A(s)^d\big)=(\sigma \seq)^* A(s')^{d'}}$ and \smash{$\var^A_t(\seq^* A(s)^m)=(\sigma \seq)^* A(s')^{m'}$}, $\pr_{\Iufv}d',\pr_{\Iufv}m'\geq 0$. Therefore, $\var^A_t(z)$ is contained in $(\sigma\seq)^*\LP^{A}(s')$.
 \end{proof}

\subsection{Existence of variation maps}

Assume that there is a linear variation map $\var_{t}^{M}\colon\Mc(t)_{\Q}\rightarrow\Mc(t')_{\Q}$.
Then it restricts to a linear map from $\Mc(t)$ to $\Mc_{\ufv}(t')\oplus\frac{1}{r}\Mc_{\fv}(t')$
for some $r>0$. Similarly, the corresponding variation map $\var_{t}^{A}\colon\LP^{A}(t)_{\fFacRoot}\rightarrow\LP^{A}(t')_{\fFacRoot}$
restricts to an algebra homomorphism \smash{$\var_{t}^{A}\colon\LP^{A}(t)\rightarrow\LP^{A}(t')_{\fFac^{\frac{1}{r}}}$}.
The restriction of a linear variation map to sublattices (resp.\ the
restriction of a monomial variation maps to subalgebras) will be still
called a variation map.

\begin{Proposition}\label{prop:exist_linear_variation}
Let there be given two seeds $t$, $t'$ similar up to a permutation $\sigma$.
Assume that the full rank assumption holds.
\begin{itemize}\itemsep=0pt
\item[$(1)$] There exists a variation map $\var_{t}^{M}\colon\Mc(t)_{\Q}\rightarrow\Mc(t')_{\Q}$.
\item[$(2)$] The set of variation maps from $\Mc(t)_{\Q}$ to $\Mc(t')_{\Q}$
takes the form $\bigl\{\var_{t}^{M}+z\mid z\in Z_{\Q}\bigr\}$, where $Z_{\Q}$ is a $\Q$-vector space
of dimension $|I_{\fv}|\cdot|\Iufv|$.
\item[$(3)$] If an order $|\Iufv|$ minor of $\tB$ is $\pm1$, then there
exists a variation map $\var_{t}^{M}\colon\Mc(t)\rightarrow\Mc(t')$. Moreover,
the set of variation maps from $\Mc(t)$ to $\Mc(t')$ takes the form
$\bigl\{\var_{t}^{M}\allowbreak+z\mid z\in Z_{\Z}\bigr\}$, where $Z_{\Z}$ is a lattice of rank $|I_{\fv}|\cdot|\Iufv|$.
\end{itemize}
\end{Proposition}

\begin{proof}
Relabeling the vertices $\sigma k$ by $k$ for $t'$ if necessary,
we assume $\sigma=\Id$ in this proof.

(1) We need to construct a linear map $\var_{t}^{M}\colon\Mc(t)_{\Q}\rightarrow\Mc(t')_{\Q}$
represented by a matrix \smash{$\left(\begin{smallmatrix}
\Id_{\ufv} & 0\\
U_{\low} & U_{\fv}
\end{smallmatrix}\right)$} in the bases $\uff$, $\uff'$ such that \eqref{eq:solve_cluster_A_twist}
holds. The unknown matrix $U:=\begin{pmatrix}
U_{\low} & U_{\fv}\end{pmatrix}$ needs to satisfy the following inhomogeneous equation:
\begin{gather}
U\tB = B'_{\low} \label{eq:matrix_eq_inhom}.
\end{gather}

Notice that \smash{$\tB=\left(\begin{smallmatrix}
B_{\ufv}\\
B_{\low}
\end{smallmatrix}\right)$} is of full rank. Then we can always choose a (not necessarily unique)
size-$|\Iufv|$ subset $J\subset I$ such that the submatrix $B_{J,\Iufv}$
is full rank. Its inverse matrix $(B_{J,\Iufv})^{-1}$ has entries
in
\smash{$\frac{1}{\det B_{J,\Iufv}}\Z$}.
Let $u_i$, $i\in I$, denote the (unknown) columns of $U$ and $b_i$, $i\in I$, denote the rows of $\tilde{B}$, then
\[U\tilde{B}=\sum_{i\in I}u_ib_i.\]
Taking $u_i=0$ for $i\in I\backslash J$ and letting $u_i$, $i\in J$ equal the columns of the $I_{\fv}\times J$-matrix $B'_{\low}(B_{J,\Iufv})^{-1}$, we obtain a special solution $U_0$ for \eqref{eq:matrix_eq_inhom}, whose entries lie in $\frac{1}{\det B_{J,\Iufv}}\Z$.

(2) Let $U_{0}$ denote the special solution in (1). Let $Z_{\Q}$ denote
the set of solutions of the homogeneous linear system
\begin{gather}\label{eq:matrix_eq_homog}
U\tB=0, \qquad \text{or equivalently}\quad\sum_{i\in J}u_ib_i=-\sum_{i\in I\setminus J}u_ib_i.
\end{gather}
Equation \eqref{eq:matrix_eq_homog} has a unique solution $u_i$, $i\in I$, for any given $u_i$, $i\in I\backslash J$. Hence the set of $\Q$-solutions is a vector space of dimension $|I_{\fv}|\cdot|\Iufv|$.
It follows that the solutions for \eqref{eq:matrix_eq_inhom} take
the form $U+U_{0}$, $U\in Z_{\Q}$.

(3) It is a direct consequence of the above argument since $\det B_{J,\Iufv}$ is $\pm 1$.
\end{proof}

\subsection{Change of seeds}

We treat lattices $\Mc(t)$ and the corresponding quantum torus algebras
$\LP^A(t)$ in this subsection. Our arguments and results remain valid
after the extension to $\Mc(t)_{\Q}$ and $\LP^A(t)_{\fFacRoot}$.

Let $t',t\in\Delta^{+}$ denote seeds similar up to a permutation
$\sigma$. Let there be given a linear variation map $\var_{t}^{M}\colon\Mc(t)\rightarrow\Mc(t')$.
Choose any $k\in\Iufv$ and denote $s=\mu_{k}t$, ${s'=\mu_{\sigma k}t'}$.
Recall that we have linear isomorphisms $\psi_{t,s}^{M}\colon\Mc(s)\simeq\Mc(t)$
and $\psi_{t',s'}^{M}\colon\Mc(s')\simeq\Mc(t')$, which are represented
by the matrices $\trans_{k,+}^{M}(t)$ and $\trans_{\sigma k,+}^{M}(t')$,
respectively. Define \[\var_{s}^{M}=\big(\psi_{t',s'}^{M}\big)^{-1}\var_{t}^{M}\psi_{t,s}^{M}\colon\ \Mc(s)\simeq\Mc(s')\]
so that the following diagram commutes:
\begin{align}
\begin{array}{@{}ccc}
\Mc(t) & \xrightarrow{\var_{t}^{M}} & \Mc(t')\\
\uparrow\psi_{t,s}^{M} & & \uparrow\psi_{t',s'}^{M}\\
\Mc(s) & \xrightarrow{\var_{s}^{M}} & \,\Mc(s').
\end{array}\label{eq:commute_M_variation}
\end{align}

The following properties for $\var_{s}^{M}$ are analogous to those
in \cite[Proposition 5.1.6]{qin2020dual}.

\begin{Proposition}\label{prop:change_seed_M_variation}\quad
\begin{itemize}\itemsep=0pt
\item[$(1)$] The linear map $\var_{s}^{M}$ is a variation map.
\item[$(2)$] If $\var_{t}^{M}$ is invertible, so is $\var_{s}^{M}$.
\item[$(3)$] If $\var_{t}^{M}$ is a Poisson variation map, so is $\var_{s}^{M}$.
\end{itemize}
\end{Proposition}

\begin{proof}
As before, we can assume $\sigma=\Id$ by relabeling the vertices
for $t'$ and $s'$. Then $\trans_{k,+}^{M}(t)_{\ufv}$ coincides
with $\trans_{k,+}^{M}(t')_{\ufv}$ by the similarity between $t$
and $t'$. It is of the block matrix form
\[\begin{pmatrix}
P_{\ufv} & 0\\
P_{\fv,\ufv} & \Id_{\fv}
\end{pmatrix}.\]
(1) The statement can be translated from \cite[Proposition 5.1.6]{qin2020dual}
by using Lemma \ref{lem:deg_variation_change_A}. Let us give a direct
proof.

The linear map $\var_{s}^{M}$ has the following matrix representation:
\begin{align*}
\var_{s}^{M}(\uff(s)) & =\big(\psi_{t',s'}^{M}\big)^{-1}\var_{t}^{M}\psi_{t,s}^{M}(\uff(s))=\big(\psi_{t',s'}^{M}\big)^{-1}\var_{t}^{M}(\uff(t)\trans_{k,+}^{M}(t))\\
 & =\big(\psi_{t',s'}^{M}\big)^{-1}\left(\uff(t')\begin{pmatrix}
\Id_{\ufv} & 0\\
U_{\low} & U_{\fv}
\end{pmatrix}\trans_{k,+}^{M}(t)\right)=\uff(s')\trans_{k,+}^{M}(t')^{-1}\begin{pmatrix}
\Id_{\ufv} & 0\\
U_{\low} & U_{\fv}
\end{pmatrix}\trans_{k,+}^{M}(t).
\end{align*}
We see that the matrix \smash{$\trans_{k,+}^{M}(t')^{-1}\left(\begin{smallmatrix}
\Id_{\ufv} & 0\\
U_{\low} & U_{\fv}
\end{smallmatrix}\right)\trans_{k,+}^{M}(t)$} is block diagonal, such that its $\Iufv\times\Iufv$ submatrix is
$\Id_{\ufv}$ (by the similarity between $t$, $t'$) and its $I_{\ufv}\times I_{\fv}$
submatrix is $0$.

We need to check \eqref{eq:preserve_p_map}. Let $\uee_{\ufv}$ denote
the matrix $(\ee_{i})_{i\in\Iufv}$ where $\ee_{i}$ are column vectors.
Then $\trans_{k,+}^{N}(t)_{\ufv}$ coincides with $\trans_{k,+}^{N}(t')_{\ufv}$
by the similarity between $t$ and $t'$. Recall that, for $j\in I_{\ufv}$, we have \[\col_{j}\big(\uee(t)\trans_{k,+}^{N}(t)\big)=\begin{cases}
 \ee_j(t)+[b_{kj}]_+\ee_k(t), & j\neq k,\\
 -\ee_k(t), & j=k.
\end{cases}\]
We deduce that $\col_{j}\big(\uee(t)\trans_{k,+}^{N}(t)\big)=\col_{j}\big(\uee(t)_{\ufv}\trans_{k,+}^{N}(t)_{\ufv}\big)$. Using Lemma \ref{lem:bilinear_form_basis_change},
we have
\begin{align*}
\psi_{t',s'}^{M}\var_{s}^{M}(p^{*}\ee_{j}(s)) & =\var_{t}^{M}\psi_{t,s}^{M}(p^{*}\ee_{j}(s))=\var_{t}^{M}p^{*}\big(\psi_{t,s}^{N}\ee_{j}(s)\big)\\
 & =\var_{t}^{M}p^{*}\bigl(\col_{j}(\uee(t)\trans_{k,+}^{N}(t))\bigr)=\var_{t}^{M}p^{*}\bigl(\col_{j}(\uee(t)_{\ufv}\trans_{k,+}^{N}(t)_{\ufv})\bigr)\\
 & =\col_{j}\var_{t}^{M}p^{*}\big(\uee(t)_{\ufv}\cdot\trans_{k,+}^{N}(t)_{\ufv}\big).
\end{align*}
We also have
\begin{align*}
\psi_{t',s'}^{M}(p^{*}\ee_{j}(s')) & =p^{*}\big(\psi_{t',s'}^{N}\ee_{j}(s')\big)=p^{*}\bigl(\col_{j}\big(\uee(t')\trans_{k,+}^{N}(t')\big)\bigr)\\
 & =p^{*}\bigl(\col_{j}\big(\uee(t')_{\ufv}\trans_{k,+}^{N}(t')_{\ufv}\big)\bigr) =\col_{j}p^{*}\big(\uee(t')_{\ufv}\cdot\trans_{k,+}^{N}(t')_{\ufv}\big).
\end{align*}

Recall that $\var_{t}^{M}(p^{*}\ee_{j}(t))=p^{*}\ee_{j}(t')$ for
any $j\in I_{\ufv}$. Then
\[\var_{t}^{M}\big(p^{*}(\uee(t)_{\ufv}\cdot\trans_{k,+}^{N}(t)_{\ufv})\big)=p^{*}\big(\uee(t')_{\ufv}\cdot\trans_{k,+}^{N}(t')_{\ufv}\big).\]
Therefore $\psi_{t',s'}^{M}\var_{s}^{M}(p^{*}\ee_{j}(s))=\psi_{t',s'}^{M}(p^{*}\ee_{j}(s'))$,
$\forall j\in I_{\ufv}$. It follow that $\var_{s}^{M}$ is a variation
map.

(2) The statement is obvious.

(3) By Lemma \ref{lem:bilinear_form_basis_change}, $\psi_{t,s}^{M}$
and $\psi_{t',s'}^{M}$ preserve the bilinear form $\lambda$. The
claim follows.
\end{proof}

Let $\var_{t}^{A}\colon\LP^{A}(t)\rightarrow \LP^{A}(t')$ denote the monomial map associated to $\var_{t}^{M}$. We further assume that it is a variation map (equivalently, $\var_t^M$ needs to be Poisson if we work at the quantum level, see Remark \ref{rem:Poisson_quantum_variation}).

\begin{Lemma}\label{lem:deg_variation_change_A}
 We have $\var_{s}^{M}\ff_{i}(s)=\deg_{s'}\big(\big(\mu_{\sigma k}^{A}\big)^{*}\big)^{-1}\var_{t}^{A}\big(\mu_{k}^{A}\big)^{*}A_{i}(s)$.
 \end{Lemma}

 \begin{proof}
 On the one hand, the $i$-th cluster variable $\big(\mu_{k}^{A}\big)^{*}A_{i}(s)$
 of $s$, $i\in I$, is pointed at \smash{$\deg_{t}\big(\mu_{k}^{A}\big)^{*}A_{i}(s)=\psi_{t,s}^{M}\ff_{i}(s)$}
 in $\LP^{A}(t)$. Since $\var_{t}^{A}$ is a variation map, $\var_{t}^{A}\big(\mu_{k}^{A}\big)^{*}A_{i}(s)$
 is pointed at $\var_{t}^{M}\psi_{t,s}^{M}\ff_{i}(s)$ in $\LP(t')$.
 In particular, we have
 \begin{align*}
 \deg_{t'}\var_{t}^{A}\big(\mu_{k}^{A}\big)^{*}A_{i}(s) & =\var_{t}^{M}\psi_{t,s}^{M}\ff_{i}(s).
 \end{align*}

 On the other hand, by the similarity between $t$ and $t'$, the image
 \smash{$\var_{t}^{A}\big(\mu_{k}^{A}\big)^{*}A_{i}(s)$} of the $i$-th cluster variable
 \smash{$\big(\mu_{k}^{A}\big)^{*}A_{i}(s)$} agrees with the $\sigma i$-th cluster
 variable \smash{$\big(\mu_{\sigma k}^{A}\big)^{*}A_{\sigma i}(s')$} up to a frozen
 factor, see \cite[Lemma 4.2.2]{qin2017triangular}. It follows that,
 similar to the cluster variable $\big(\mu_{\sigma k}^{A}\big)^{*}A_{\sigma i}(s')$,
 the degree of \smash{$\var_{t}^{A}\big(\mu_{k}^{A}\big)^{*}A_{i}(s)$} in $t'$ can
 be computed from the degree of its Laurent expansion \smash{$\big(\big(\mu_{\sigma k}^{A}\big)^{*}\big)^{-1}\var_{t}^{A}\big(\mu_{k}^{A}\big)^{*}A_{i}(s)$}
 in $s'$ by using the linear map $\psi_{t',s'}^{M}$:
 \begin{align*}
 \deg_{t'}\var_{t}^{A}\big(\mu_{k}^{A}\big)^{*}A_{i}(s) & =\psi_{t',s'}^{M}\deg_{s'}\big(\big(\mu_{\sigma k}^{A}\big)^{*}\big)^{-1}\var_{t}^{A}\big(\mu_{k}^{A}\big)^{*}A_{i}(s).
 \end{align*}

 Therefore, we obtain $\var_{t}^{M}\psi_{t,s}^{M}\ff_{i}(s)=\psi_{t',s'}^{M}\deg_{s'}\big(\big(\mu_{\sigma k}^{A}\big)^{*}\big)^{-1}\var_{t}^{A}(\mu_{k}^{A})^{*}A_{i}(s)$.
 The desired claim follows from the definition \smash{$\var_{s}^{M}=\big(\psi_{t',s'}^{M}\big)^{-1}\var_{t}^{M}\psi_{t,s}^{M}$}.
 \end{proof}

 We define the monomial map $\var_{s}^{A}\colon\LP^{A}(s)\rightarrow\LP^{A}(s')$ to be the monomial map associated to \smash{$\var_{s}^{M}=\big(\psi_{t',s'}^{M}\big)^{-1}\var_{t}^{M}\psi_{t,s}^{M}$}. Since $\var_{t}^A$ is a variation map, it is also a variation map, see Remark~\ref{rem:Poisson_quantum_variation}.

\begin{Proposition}\label{prop:change_seed_A_variation}
The map $\var_{t}^{A}\big(\mu_{k}^{A}\big)^{*}$ agrees with $\big(\mu_{\sigma k}^{A}\big)^{*}\var_{s}^{A}$,
i.e., the following diagram commutes:
\begin{align*}
\begin{array}{@{}ccc}
 \upClAlg(t) & \xrightarrow{\var_{t}^{A}} & \upClAlg(t')\\
\uparrow\big(\mu_{k}^{A}\big)^{*} & & \uparrow\big(\mu_{\sigma k}^{A}\big)^{*}\\
\upClAlg(s) & \xrightarrow{\var_{s}^{A}} &\,\upClAlg(s').
\end{array}
\end{align*}
\end{Proposition}

\begin{proof}
It suffices to show $\big(\big(\mu_{\sigma k}^{A}\big)^{*}\big)^{-1}\var_{t}^{A}\big(\mu_{k}^{A}\big)^{*}A_{i}(s)=\var_{s}^{A}A_{i}(s)$, $\forall i$.

Both $\big(\big(\mu_{\sigma k}^{A}\big)^{*}\big)^{-1}\var_{t}^{A}\big(\mu_{k}^{A}\big)^{*}A_{i}(s)$
and $\var_{s}^{A}A_{i}(s)$ are cluster variables in $s'$ up to a
frozen factor by \cite[Lemma 4.2.2]{qin2017triangular}. By Lemma
\ref{lem:deg_variation_change_A}, they have the same degree.
The claim follows.
\end{proof}

We can given a different proof for Proposition \ref{prop:change_seed_A_variation}
by straightforward computation (see Proposition \ref{prop:change_seed_X_variation}).

Take any mutation sequence $\seq$ and denote $s=\seq t$, $s'=(\sigma\seq)t'$.
By using Proposition \ref{prop:change_seed_A_variation} recursively,
$\var_{t}^{A}$ and $\seq$ uniquely determine a variation map $\var_{s}^{A}$
such that the following diagram is commutative:
\begin{align*}
\begin{array}{@{}ccc}
\upClAlg(t) & \xrightarrow{\var_{t}^{A}} & \upClAlg(t')\\
\uparrow(\seq)^{*} & & \uparrow(\sigma\seq)^{*}\\
\upClAlg(s) & \xrightarrow{\var_{s}^{A}} & \,\upClAlg(s').
\end{array}
\end{align*}
By tracking the degrees of cluster variables, we see that the variation
map $\var_{s}^{M}$ associated to~$\var_{s}^{A}$ is determined
by $\var_{t}^{M}$ and $\seq$. It still makes \eqref{eq:commute_M_variation}
commutative.

\subsection{Twist endomorphisms}

We are now ready to define twist endomorphisms. Notice that the mutation
maps always preserves the compatible Poisson structures and the $v$-twisted products. Let $t$ and $t'=\seq t$ be two
seeds similar up to $\sigma$ as before. Let there be given a variation
map $\var_{t}^{A}$ from $\LP^{A}(t)_{\fFacRoot}$ to~$\LP^{A}(t')_{\fFacRoot}$.
Note that it induces a homomorphism $\var_{t}^{A}\colon\upClAlg(t)_{\fFacRoot}\rightarrow \upClAlg(t')_{\fFacRoot}$.

\begin{Definition}
The composition $\seq^{*}\var_{t}^{A}\colon\upClAlg(t)_{\fFacRoot}\rightarrow\upClAlg(t)_{\fFacRoot}$
is called a twist endomorphism passing through the seed $t'=\seq t$, which
is denoted by $\tw_{t}^{A}$.

For classical cases, if $t$ and $t'$ are equipped with compatible Poisson structures
which are preserved by $\tw_{t}^{A}$, then $\tw_{t}^{A}$ is called
a Poisson twist endomorphism.
\end{Definition}

Note that $\var_{t}^{A}$ restricts to a variation map \smash{$\var_{t}^{A}\colon\LP(t)\rightarrow\LP(t')_{\fFac^{\frac{1}{r}}}$}
for some integer $r>0$. So the twist endomorphism $\tw_{t}^{A}$
restricts to an algebra homomorphism from $\upClAlg(t)$ to $\upClAlg(t)_{\fFac^{\frac{1}{r}}}$.

Propositions \ref{prop:change_seed_M_variation} and \ref{prop:change_seed_A_variation}
have the following consequence.

\begin{Corollary}
A twist endomorphism $\tw_{t}^{A}$ on $\upClAlg(t)_{\fFacRoot}$ $($resp.\ on $\upClAlg(t))$ gives rise to twist endomorphisms $\tw_{s}^{A}$ on
$\upClAlg(s)_{\fFacRoot}$ $($resp.\ on $\upClAlg(s))$ for all seeds $s\in\Delta^{+}$
via the mutation maps. For classical cases, if $\tw_{t}^{A}$ is Poisson,
so are $\tw_{s}^{A}$.
\end{Corollary}

Correspondingly, we can simply denote the twist endomorphisms $\tw_{t}^{A}$
by $\tw^{A}$ when we do not want to choose a specific seed $t$.

\begin{Proposition}\label{prop:A_variation_map_restrict_cluster_alg}
The twist endomorphism $\tw^{A}$ restricts to an endomorphism
on $\clAlg_{\fFacRoot}$.
\end{Proposition}

\begin{proof}
Recall that $\var_{t}^{A}$ sends cluster variables to cluster
variables up to a frozen factor, see \cite[Lemma 4.2.2]{qin2017triangular}.
The claim follows.
\end{proof}

If $\var_{t}^{A}$ is an invertible variation map, then we have a
twist automorphism $\tw_{t}^{A}:=\seq^{*}\var_{t}^{A}$ on~$\upClAlg(t)_{\fFacRoot}$.
Lemma \ref{lem:inverse_M_variation} implies that \smash{$\bigl(\var_{t}^{A}\bigr)^{-1}$}
is still a variation map. So \smash{$\eta:=(\seq^{*})^{-1}\bigl(\var_{t}^{A}\bigr)^{-1}$}
is a twist automorphism on $\upClAlg(t')_{\fFacRoot}$.

\begin{Lemma}\label{lem:inverse_A_twist_automorphism}
We have $\seq^{*}\eta(\seq^{*})^{-1}=\bigl(\tw_{t}^{A}\bigr)^{-1}$. Namely,
by identifying the fraction fields by mutations, the twist automorphisms
associated to invertible variation maps $\smash{\var_{t}^{A}\colon\LP^{A}(t)_{\fFacRoot}\rightarrow}\allowbreak\LP^{A}(t')_{\fFacRoot}$
and \smash{$\bigl(\var_{t}^{A}\bigr)^{-1}\colon\LP^{A}(t')_{\fFacRoot}\rightarrow\LP^{A}(t)_{\fFacRoot}$} are inverse
to each other.
\end{Lemma}

\begin{proof}
We have
\begin{align*}
\seq^{*}\eta(\seq^{*})^{-1} & =\seq^{*}\circ(\seq^{*})^{-1}\bigl(\var_{t}^{A}\bigr)^{-1}\circ(\seq^{*})^{-1} =\bigl(\var_{t}^{A}\bigr)^{-1}\circ(\seq^{*})^{-1} =\bigl(\tw_{t}^{A}\bigr)^{-1}.\tag*{\qed}
\end{align*} \renewcommand{\qed}{}
\end{proof}

Proposition \ref{prop:exist_linear_variation} implies the following
result.

\begin{Theorem}\label{thm:twist_A}
Consider the classical case $\kk=\Z$. Assume that $t,t'\in\Delta^{+}$ are similar and the full rank assumption
holds, then the following statements are true.
\begin{itemize}\itemsep=0pt
\item[$(1)$] There exists a twist endomorphism $\tw_{t}^{A}$ on $\upClAlg_{\fFacRoot}$
passing through the seed $t'$. The set of such twist endomorphisms
is in bijection with a $\Q$-vector space of dimension $|\Iufv|\cdot|\Ifv|$.
\item[$(2)$] If an order $|\Iufv|$ minor of $\tB$ is $\pm1$, then there
exists a twist endomorphism $\tw_{t}^{A}$ on $\upClAlg$ passing
through the seed $t'$. The set of such twist endomorphisms is in
bijection with a~lattice of rank $|\Iufv|\cdot|\Ifv|$.
\end{itemize}
\end{Theorem}

\section[Twist endomorphisms for upper cluster X-algebras]{Twist endomorphisms for upper cluster $\boldsymbol{X}$-algebras}\label{sec:Twist-endomorphisms-X}

\subsection{Variation maps}\label{subsec:Variation-N-maps}

Let $t$, $t'$ be similar seeds as before. In the following calculation,
we choose $t$ as the initial seed. We want to construct twist endomorphism
on $\XClAlg(t')$ preserving the canonical Poisson structure~\eqref{eq:log-Poisson}.

Let us first investigate linear maps between lattices $N(t)$ which
arise as the pullback of those between $\Mc(t)$. Let $\var_{t}^{M}$
denote a $\Z$-linear map $\Mc(t)\rightarrow\Mc(t')_{\Q}$,
its pullback between the dual lattices gives rise to a $\Z$-linear
map $\bigl(\var_{t}^{M}\bigr)^{*}$ from $N(t')$ to $N(t)_{\Q}$. Recall that
the diagonal entries of $D$ are $\frac{1}{d_{i}}$, $i\in I$.

\begin{Lemma}\label{lem:pullback_matrix}

The pullback $\bigl(\var_{t}^{M}\bigr)^{*}\colon N(t')\rightarrow N(t)_{\Q}$ is represented
by the matrix \[D^{-1}\begin{pmatrix}
\trans_{\sigma^{-1}} & 0\\
0 & \Id_{\fv}
\end{pmatrix}\begin{pmatrix}
\Id_{\ufv} & U_{\low}^{\rm T}\\
0 & U_{\fv}^{\rm T}
\end{pmatrix}D\]
 with respect to the bases $\uee$ and $\uee'$.

\end{Lemma}

\begin{proof}

By assumption, we have \smash{$\var_{t}^{M}(\uff)=\uff'\cdot\left(\begin{smallmatrix}
\Id_{\ufv} & 0\\
U_{\low} & U_{\fv}
\end{smallmatrix}\right)\left(\begin{smallmatrix}
\trans_{\sigma} & 0\\
0 & \Id_{\fv}
\end{smallmatrix}\right)$}. Then its pullback with respect to the dual bases is represented
by the transpose:
\begin{gather*}
\bigl(\var_{t}^{M}\bigr)^{*}((\uff')^{*}) = (\uff^{*})\cdot\begin{pmatrix}
\trans_{\sigma^{-1}} & 0\\
0 & \Id_{\fv}
\end{pmatrix}\begin{pmatrix}
\Id_{\ufv} & U_{\low}^{\rm T}\\
0 & U_{\fv}^{\rm T}
\end{pmatrix}.
\end{gather*}
Since $\ff_{i}^{*}=d_{i}\ee_{i}$ and the diagonal entries of $D$
are $\frac{1}{d_{i}}$., we have $\uff^{*}=\uee D^{-1}$ and, similarly,
$(\uff')^{*}=(\uee')D^{-1}$. So we obtain
\begin{align*}
\bigl(\var_{t}^{M}\bigr)^{*}\bigl(\uee'D^{-1}\bigr) & =\uee D^{-1}\cdot\begin{pmatrix}
\trans_{\sigma^{-1}} & 0\\
0 & \Id_{\fv}
\end{pmatrix}\begin{pmatrix}
\Id_{\ufv} & U_{\low}^{\rm T}\\
0 & U_{\fv}^{\rm T}
\end{pmatrix},
\end{align*}
and, equivalently, $\bigl(\var_{t}^{M}\bigr)^{*}(\uee')=\uee D^{-1}\left(\begin{smallmatrix}
\trans_{\sigma^{-1}} & 0\\
0 & \Id_{\fv}
\end{smallmatrix}\right)\left(\begin{smallmatrix}
\Id_{\ufv} & U_{\low}^{\rm T}\\
0 & U_{\fv}^{\rm T}
\end{smallmatrix}\right)D$.
\end{proof}

Since $t$ and $t'$ are similar, we have $d_{k}=d_{\sigma k}$ by
definition. Then it follows that \[D^{-1}\begin{pmatrix}
\trans_{\sigma^{-1}} & 0\\
0 & \Id_{\fv}
\end{pmatrix}=\begin{pmatrix}
\trans_{\sigma^{-1}} & 0\\
0 & \Id_{\fv}
\end{pmatrix}D^{-1}.\] In view of Lemma \ref{lem:pullback_matrix}, we propose the following
definition.

\begin{Definition}\label{def:N_variation}

A linear map $\var_{t'}^{N}\colon N(t')_{\Q}\rightarrow N(t)_{\Q}$ is said
to be a variation map, if it has the following $\Q$-valued matrix
representation in the basis $\uee'=\{\ee'_{i}\mid i\in I\}$, $\uee=\{\ee_{i}\mid i\in I\}$:
\begin{align}
\var_{t'}^{N}(\uee') & =\uee\cdot\begin{pmatrix}
\trans_{\sigma^{-1}} & 0\\
0 & \Id_{\fv}
\end{pmatrix}\begin{pmatrix}
\Id_{\ufv} & V_{\high}\\
0 & V_{\fv}
\end{pmatrix},\label{eq:N_variation_triangular}
\end{align}
and it satisfies
\begin{align}\label{eq:N_variation_condition}
\omega\big(\var_{t'}^{N}\ee'_{i},\var_{t'}^{N}\ee'_{k}\big) & =\omega(\ee'_{i},\ee'_{k}),\qquad \forall i\in I,\quad k\in\Iufv.
\end{align}

It is said to be Poisson if it further preserves the Poisson structures:
\begin{align*}
\omega\big(\var_{t'}^{N}\ee'_{i},\var_{t'}^{N}\ee'_{j}\big) & =\omega(\ee_{i}',\ee_{j}'),\qquad \forall i,j\in I.
\end{align*}

\end{Definition}
We will see in Proposition \ref{prop:linear_N_variation} that
\eqref{eq:N_variation_condition} is natural.

Definition \ref{def:N_variation} will be naturally deduced from the definition of the following monomial variation map, see Remark \ref{rem:variation_X_N_translation}.
Note that we always have $\var_{t'}^{N}(\ee'_{\sigma k})=\ee_{k}$
for $k\in\Iufv$. In particular, \eqref{eq:N_variation_condition}
is an inhomogeneous linear system of equations on $\var_{t'}^{N}$.

\begin{Definition}\label{def:X_variation_map}
 A $\kk$-algebra homomorphism $\var_{t'}^X\colon\LP^{X}(t')\rightarrow \LP^{X}(t)$ is called a (monomial) variation map if, for any $k\in \Iufv$, $j\in \Ifv$, we have
 \begin{equation*}
 \var_{t'}^{X}((X')_k) = X_{\sigma^{-1}k},\qquad
 \var_{t'}^{X}((X')_{j}) = X^{n^{(j)}},
 \end{equation*}
 for some $n^{(j)}\in \Z^I$ satisfying, $\forall n'\in \Z^{I}$,
 \begin{align}\label{eq:var_X_B}
 \sigma \pr_{\Iufv} B(t) n&=\pr_{\Iufv}B(t')(n'),
 \end{align}
 where $X^n:=\var_{t'}^X (X')^{n'}$.

 For classical cases, $\var_{t'}^{X}$ is called a Poisson variation map if it further preserves
 the Poisson structures:
 \begin{align*}
 \var_{t'}^{X}\{(X')_{i},(X')_{j}\} =\bigl\{\var_{t'}^{X}(X')_{i},\var_{t'}^{X}(X')_{j}\bigr\},\qquad\forall i,j\in I.
 \end{align*}
\end{Definition}

\begin{Definition}
 In $\hLP^X(t')$, take any pointed formal Laurent series $Z'=(X')^{n'} \cdot F$, where
 \[n'\in \Z^I, \qquad F=\sum_{n''\in N_{\ufv}^{\geq 0}} c_{n''} (X')^{n''}, \qquad c_{n''}\in \kk, \qquad c_0=1.\] We say $Z'$ and \smash{$Z\in \hLP^X(t)$} are similar if $Z$ takes the form \[X^{n}\cdot F|_{(X')^{n''}\mapsto (X)^{\sigma^{-1} n''}}\qquad \text{with}\quad \sigma \pr_{\Iufv}B(t)n=\pr_{\Iufv}B(t')n'.\]
\end{Definition}

\begin{Remark}\label{rem:X-variation_map_pointed_elements}
 The variation map $\var_{t'}^X$ is a $\kk$-algebra homomorphism such that it sends a pointed element $Z'\in \LP^X(t')$ to a similar element in $\LP^X(t)$.
\end{Remark}

\begin{Remark}\label{rem:variation_X_N_translation}
 Let $\var^X_{t'}$ be a monomial variation map and $\var^N_{t'}$ the associated linear map. Then $\var^N_{t'}$ must satisfies \eqref{eq:N_variation_triangular}. We can study condition~\eqref{eq:var_X_B} by taking $n'$ to be $\ee_i'$ and consider the $k$-th rows of both sides for all $i\in I$, $k\in \Iufv$. Recall that \smash{$\omega(n,\ee_{\sigma^{-1}k})=\frac{1}{d_{\sigma^{-1}k}}(B n)_{\sigma^{-1}k}$}, $\omega(\ee'_i,\ee'_k)=\frac{1}{d_k}b'_{ki}$, and $d_k=d_{\sigma^{-1} k}$. Then condition \eqref{eq:var_X_B} is equivalent to
\begin{align*}
 \omega\big(\var^N_{t'}\ee_i',\ee_{\sigma^{-1} k}\big)=\omega(\ee_i',\ee_k').
\end{align*}
We see that this equation is equivalent to \eqref{eq:N_variation_condition} by using $\ee_{\sigma^{-1} k}=\var^N_{t'}\ee_k'$.
\end{Remark}

\begin{Remark}\label{rem:X-Poisson_quantum_variation}
By using Remark \ref{rem:variation_X_N_translation}, we obtain the following equivalent statements:
\begin{itemize}\itemsep=0pt
 \item A linear map $\varPsi\colon N(t')\rightarrow N(t)$ is a Poisson variation map.
 \item For the quantum case $\kk=\Z\big[v^{\pm\frac{1}{\diag}}\big]$, the monomial map $\varPhi_v$ associated to $\varPsi$ is a variation map.
 \item For the classical case $\kk=\Z$, the monomial map $\varPhi_1$ associated to $\varPsi$ is a Poisson variation map.
\end{itemize}
Therefore, a Poisson variation map is the classical limit of a quantum variation map at $v=1$, see Section \ref{sec:quantization}. Conversely, a Poisson variation map gives rise to a quantum variation map by the above equivalent statements.
\end{Remark}

\begin{Proposition}\label{prop:linear_N_variation}
A linear map $\var_{t'}^{N}\colon N(t')_{\Q}\rightarrow N(t)_{\Q}$ is a
variation map if and only if its pullback $\bigl(\var_{t'}^{N}\bigr)^{*}\colon\Mc(t)_{\Q}\rightarrow\Mc(t')_{\Q}$
is a variation map.
\end{Proposition}

\begin{proof}
Relabeling the vertices $\sigma k$ by $k$ for $t'$ if necessary,
we assume that $\sigma$ is the identity. Denote the pullback of $\var_{t'}^{N}$
by $\var_{t}^{M}$.

By Lemma \ref{lem:pullback_matrix}, $\var_{t'}^{N}$ takes the form
of \eqref{eq:N_variation_triangular} if and only if $\var_{t}^{M}$
takes the form of \eqref{eq:preserve_p_map}. It remains
to show that $\var_{t}^{M}p^{*}(\ee_{k})=p^{*}\ee_{k}'$ is equivalent
to $\omega\bigl(\var_{t'}^{N}\ee'_{i},\var_{t'}^{N}\ee'_{k}\bigr)=\omega(\ee'_{i},\ee'_{k})$, $\forall i\in I$, $k\in\Iufv$.

Recall that $p^{*}\ee_{i}(\ )=\omega(\ee_{i},\ )$, we
have
\begin{align*}
\omega\bigl(\var_{t'}^{N}\ee'_{i},\var_{t'}^{N}\ee'_{k}\bigr) =\omega\bigl(\var_{t'}^{N}\ee'_{i},\ee_{k}\bigr) =-p^{*}(\ee_{k})\bigl(\var_{t'}^{N}\ee'_{i}\bigr) =-\var_{t}^{M}p^{*}(\ee_{k})(\ee_{i}').
\end{align*}
On the other hand, $\omega(\ee'_{i},\ee'_{k})=-p^{*}(\ee'_{k})(\ee'_{i})$.

Therefore, the condition \[\omega\bigl(\var_{t'}^{N}\ee'_{i},\var_{t'}^{N}\ee'_{k}\bigr)=\omega(\ee'_{i},\ee'_{k}),\qquad\forall i\in I,\]
 is equivalent to $\var_{t}^{M}p^{*}(\ee_{k})=p^{*}(\ee'_{k})$.
\end{proof}

The following result is analogous to Lemma \ref{lem:inverse_M_variation}.

\begin{Lemma}\label{lem:inverse_N_variation}
An invertible linear map $\var_{t'}^{N}\colon N(t)\rightarrow N(t')$ is
a $($Poisson$)$ variation map if and only if its inverse is.
\end{Lemma}

\begin{proof}
Assume that $\var_{t'}^{N}\colon N(t')\rightarrow N(t)$ is represented
by an invertible matrix \[\begin{pmatrix}
\trans_{\sigma^{-1}} & 0\\
0 & \Id_{\fv}
\end{pmatrix}\begin{pmatrix}
\Id_{\ufv} & V_{\high}\\
0 & V_{\fv}
\end{pmatrix}\] as in \ref{eq:N_variation_triangular}. Then its inverse $\bigl(\var_{t'}^{N}\bigr)^{-1}$
is represented by \[\begin{pmatrix}
\Id_{\ufv} & -V_{\high}V_{\fv}^{-1}\\
0 & V_{\fv}^{-1}
\end{pmatrix}\begin{pmatrix}
\trans_{\sigma} & 0\\
0 & \Id_{\fv}
\end{pmatrix}=\begin{pmatrix}
\trans_{\sigma} & 0\\
0 & \Id_{\fv}
\end{pmatrix}\begin{pmatrix}
\Id_{\ufv} & -\trans_{\sigma}^{-1}V_{\high}V_{\fv}^{-1}\\
0 & V_{\fv}^{-1}
\end{pmatrix}.\] So $\bigl(\var_{t'}^{N}\bigr)^{-1}$ is also of the form in \eqref{eq:N_variation_triangular}.
The other conditions in Definition \ref{def:N_variation} can be checked
easily.
\end{proof}

\begin{Lemma}\label{lem:commute_variation_p}
Assume that $\var_{t'}^{N}\colon N(t')_{\Q}\rightarrow N(t)_{\Q}$ is a
linear map. Let $\var_{t}^{M}$denote its pullback. The following
diagram commutes if and only if $\var_{t'}^{N}$ preserves $\omega$
\begin{align}
\begin{array}{ccc}
\Mc(t)_{\Q} & \xrightarrow{\var_{t}^{M}} & \Mc(t')_{\Q}\\
\uparrow p^{*} & & \uparrow p^{*}\\
N(t)_{\Q} & \xleftarrow{\var_{t'}^{N}} & \,N(t')_{\Q}.
\end{array}\label{eq:commute_variation_p}
\end{align}
\end{Lemma}
Note that, by taking the matrix presentations, the commutative diagram \eqref{eq:commute_variation_p} is represented by a quadratic equation on the entries of $\var_{t'}^N$.

\begin{proof}
Take any $i,j\in I$. On the one hand, we have $p^{*}\ee_{i}'(\ee_{j}')=\omega(\ee'_{i},\ee'_{j})$.
On the other hand, we have
\begin{align*}
\var_{t}^{M} p^{*} \var_{t'}^{N}(\ee'_{i})(\ee'_{j}) & =p^{*}\var_{t'}^{N}(\ee'_{i})\big(\var_{t'}^{N}\ee'_{j}\big) =\omega\big(\var_{t'}^{N}\ee'_{i},\var_{t'}^{N}\ee'_{j}\big).
\end{align*}

Therefore, $p^{*}=\var_{t}^{M}p^{*}\var_{t'}^{N}$ if and only if
$\var_{t'}^{N}$ preserves $\omega$.
\end{proof}

\subsection{Change of seeds}

Let $t',t\in\Delta^{+}$ denote seeds similar up to a permutation
$\sigma$. For any $k\in I_{\ufv}$, consider the seeds~$s=\mu_{k}t$
and $s'=\mu_{\sigma k}t'$.

We define the linear map $\var_{s'}^{N}\colon N(s')\rightarrow N(s)$ to
be $\bigl(\psi_{t,s}^{N}\bigr)^{-1}\var_{t'}^{N}\psi_{t',s'}^{N}$ such that
the following diagram commutes:
\begin{align}
\begin{array}{@{}ccc}
N(t) & \xleftarrow{\var_{t'}^{N}} & N(t')\\
\uparrow\psi_{t,s}^{N} & & \uparrow\psi_{t',s'}^{N}\\
N(s) & \xleftarrow{\var_{s'}^{N}} & N(s').
\end{array}\label{eq:commute_N_variation}
\end{align}
By Proposition \ref{prop:matrix_properties}\,(3), the pullback $\big(\psi_{t,s}^{M}\big)^{*}$
is represented by the matrix $\smash{D^{-1}\bigl(\trans_{k,+}^{M}(t)\bigr)^{\rm T}D=}\allowbreak\smash{\trans_{k,+}^{N}(t)^{-1}}$.
So we have \smash{$\big(\psi_{t,s}^{M}\big)^{*}=\big(\psi_{t,s}^{N}\big)^{-1}$}. Diagram
\eqref{eq:commute_N_variation} should be compared with the pullback
of \eqref{eq:commute_M_variation}. We have the following result in
analogous of Proposition \ref{prop:change_seed_M_variation}.

\begin{Proposition}\label{prop:change_seed_N_variation}\quad
\begin{itemize}\itemsep=0pt
\item[$(1)$] If $\var_{t'}^{N}$ is a variation map, so is $\var_{s'}^{N}$.
\item[$(2)$] If $\var_{t'}^{N}$ is invertible, so is $\var_{s'}^{N}$.
\item[$(3)$] If $\var_{t'}^{N}$ is Poisson, so is $\var_{s'}^{N}$.
\end{itemize}
\end{Proposition}

\begin{proof}
(2) The statement is obvious.

(1)--(3) As in the proof of Proposition \ref{prop:change_seed_M_variation}\,(1),
straightforward computation shows that $\var_{s'}^{N}$ is a block
triangular matrix of the form in \eqref{eq:N_variation_triangular}.
By Lemma \ref{lem:bilinear_form_basis_change}, $\psi_{s,t}^{N}$
and $\psi_{s',t'}^{N}$ preserve the Poisson structure $\omega$.
So if \smash{$\omega\bigl(\var_{t'}^{N}\ee{}_{i}(t'),\var_{t'}^{N}\ee{}_{j}(t')\bigr)=\omega(\ee{}_{i}(t'),\ee{}_{j}(t'))$}
for some $i$, $j$, then $\omega\bigl(\var_{s'}^{N}\ee{}_{i}(s'),\var_{s'}^{N}\ee{}_{j}(s')\bigr)=\omega(\ee{}_{i}(s'),\ee{}_{j}(s'))$.
The desired claims follow.
\end{proof}

Let $\var_{t'}^{X}$ and $\var_{s'}^{X}$ be the monomial maps associated to the linear maps in Proposition \ref{prop:change_seed_N_variation}. Assume that $\var_{t'}^{X}$ is a variation map. Combining Proposition \ref{prop:change_seed_N_variation} and Remark \ref{rem:X-Poisson_quantum_variation}, we obtain that $\var_{s'}^{X}$ is also a variation map.

Let $T$ denote the subalgebra $\kk[X_k]_{k\in \Iufv}$ of the quantum torus algebra $\LP^X(t)$. Let $\cF_{\ufv}^X(t)$ denote the subalgebra of the fraction field $\cF^X(t)$ such that its elements take the form~$P*Q^{-1}$ with $P\in \LP^X(t)$, $Q\in T$. Then it is easy to check $\big(\mu_k^X\big)^* \cF_{\ufv}^X(t')=\cF_{\ufv}^X(t)$ by using the invertibility of $\mu_k^X$.

We naturally extend the variation map $\var_{t'}^X\colon\LP^X(t')\rightarrow\LP^X(t)$ to an algebra homomorphism $\var_{t'}^X$ from $\cF_{\ufv}^X(t')$ to $\cF_{\ufv}^X(t)$.

\begin{Proposition}\label{prop:change_seed_X_variation}
We have $\big(\mu_{k}^{X}\big)^{*}\var_{s'}^{X}=\var_{t'}^{X}\big(\mu_{\sigma k}^{X}\big)^{*}$,
i.e., the following diagram commutes:
\begin{align*}
 \begin{array}{@{}ccc}
 \cF^X_{\ufv}(t) & \xleftarrow{\var_{t'}^{X}} & \cF^X_{\ufv}(t')\\
 \big(\mu_{k}^{X}\big)^{*}\uparrow & & \uparrow\big(\mu_{\sigma k}^{X}\big)^{*}\\
 \cF^X_{\ufv}(s) & \xleftarrow{\var_{s'}^{X}} & \cF^X_{\ufv}(s').
 \end{array}
 \end{align*}
\end{Proposition}

\begin{proof}
 We will denote $\mu_{k}^{*}=\big(\mu_{k}^{X}\big)^{*}$ and $\psi=\psi^N$ below for simplicity. Relabeling the vertices for $s'$, $t'$ if necessary, we can
assume $\sigma=\Id$. It suffices to show $\var_{t'}^{X}\mu_{k}^{*}X_{i}(s')=\mu_{k}^{*}\var_{s'}^{X}X_{i}(s')$,~$\forall i$.

(i) Let us first consider classical cases. We have
\begin{align}
\mu_{k}^{*}X_{i}(s') & =\begin{cases}
X_{i}(t')\cdot X_{k}(t')^{[b_{ki}(t')]_{+}}\cdot(1+X_{k}(t'))^{-b_{ki}(t')}, & i\neq k,\\
X_{k}(t')^{-1}, & i=k
\end{cases} \nonumber\\
 & =X(t')^{\psi_{t',s'}\ee_{i}(s')}\cdot(1+X_{k}(t'))^{b_{ki}(s')}\nonumber \\
 & =X(t')^{\psi_{t',s'}\ee_{i}(s')}\cdot(1+X_{k}(t'))^{-\omega(\ee_{k}(s'),\ee_{i}(s'))d_{k}}. \label{eq:mutate_X_expression}
\end{align}
We set $\varepsilon=1$ if $-\omega(\ee_{k}(s'),\ee_{i}(s'))\geq 0$ and $\varepsilon=-1$ otherwise. Then we have
\begin{align}\label{eq:var_mu_X}
 \var_{t'}^{X}\mu_{k}^{*}X_{i}(s')^\varepsilon=X(t)^{\varepsilon\var_{t'}^{N}\psi_{t',s'}\ee_{i}(s')}\cdot(1+X_{k}(t))^{-\varepsilon\omega(\ee_{k}(s'),\ee_{i}(s'))d_{k}}.
\end{align}

Similar to \eqref{eq:mutate_X_expression}, we have
\begin{align}\label{eq:mu_var_X}
\mu_{k}^{*}\var_{s'}^{X}X_{i}(s')^{{\varepsilon'}} & =\mu_{k}^{*}X(s)^{{\varepsilon'}\var_{s'}^{N}\ee_{i}(s')} \nonumber\\
 & =X(t)^{{\varepsilon'}\psi_{t,s}\var_{s'}^{N}\ee_{i}(s')}\cdot (1+X_{k}(t))^{-{\varepsilon'}\omega(\ee_{k}(s),\var_{s'}^{N}\ee_{i}(s'))d_{k}},
\end{align}
where we set $\varepsilon'=1$ if $-\omega\big(\ee_{k}(s),\var_{s'}^{N}\ee_{i}(s')\big)\geq 0$ and $\varepsilon'=-1$ otherwise.

Note that $\var_{t'}^{N}\psi_{t',s'}\ee_{i}(s')=\psi_{t,s}\var_{s'}^{N}\ee_{i}(s')$
by the definition of $\var_{s'}^{N}$. Moreover, since $\var_{s'}^{N}$
is a variation map by Proposition \ref{prop:change_seed_N_variation},
we have \[\omega\big(\ee_{k}(s),\var_{s'}^{N}\ee_{i}(s')\big)=\omega\bigl(\var_{s'}^{N}\ee_{k}(s'),\var_{s'}^{N}\ee_{i}(s')\bigr)=\omega(\ee_{k}(s'),\ee_{i}(s')).\]
It follows that $\varepsilon=\varepsilon'$ and $\var_{t}^{X}\mu_{k}^{*}X_{i}(s')^{\varepsilon}=\mu_{k}^{*}\var_{s'}^{X}X_{i}(s')^{{\varepsilon'}}$.

(ii) In the quantum case, we can obtain $\var_{t}^{X}\mu_{k}^{*}X_{i}(s')^\varepsilon$ and $\mu_{k}^{*}\var_{s'}^{X}X_{i}(s')^{\varepsilon'}$ by replacing the binomial coefficients $\left(\begin{smallmatrix}
 a\\
 b\\
 \end{smallmatrix}\right)$ in polynomial expansions of the right hand sides of \eqref{eq:var_mu_X} and \eqref{eq:mu_var_X} by the quantum numbers $\left(\begin{smallmatrix}
 a\\
 b\\
 \end{smallmatrix}\right)_{v_k}$. Therefore, we still have \[\var_{t}^{X}\mu_{k}^{*}X_{i}(s')=\mu_{k}^{*}\var_{s'}^{X}X_{i}(s').\tag*{\qed} \]\renewcommand{\qed}{}
\end{proof}

The following statement follows from Proposition~\ref{prop:change_seed_X_variation}
by tracking the degrees of $X$-variables. It is an analog of Lemma~\ref{lem:deg_variation_change_A}.

\begin{Lemma}
The variation map $\var_{s'}^{N}\colon N(s')\rightarrow N(s)$ is the linear
map sending $\ee_{i}(s')$ to \smash{$\deg_{s}\big(\big(\mu_{k}^{X}\big)^{*}\big)^{-1}\var_{t'}^{X}\big(\mu_{\sigma k}^{X}\big)^{*}X_{i}(s')$}.
\end{Lemma}

Take any mutation sequence $\seq$ and denote $s=\seq t$, $s'=(\sigma\seq)t'$.
By using Proposition \ref{prop:change_seed_X_variation} recursively, $\var_{t'}^{X}$
and $\seq$ uniquely determine a variation map $\var_{s'}^{X}$ such
that the following diagram is commutative:
\begin{align}
\begin{array}{@{}ccc}
 \cF^X_{\ufv}(t) & \xleftarrow{\var_{t'}^{X}} & \cF^X_{\ufv}(t')\\
\uparrow\seq^{*} & & \uparrow(\sigma\seq)^{*}\\
\cF^X_{\ufv}(s) & \xleftarrow{\var_{s'}^{X}} & \cF^X_{\ufv}(s').
\end{array}\label{eq:commute_X_variation_seq}
\end{align}
By tracking the degrees of cluster variables, we see that linear the variation
map $\var_{s'}^{N}$ associated to $\var_{s'}^{X}$ is determined
by $\var_{t'}^{N}$ and $\seq$. It still makes \eqref{eq:commute_N_variation}
commutative.

Note that, by the commutative diagrams \eqref{eq:commute_X_variation_seq},
the variation maps $\var_{s'}^{X}$, $s'\in\Delta_{s'}^{+}$, defined
using $\var_{t'}^{X}$, are identified via the mutation maps.

 We have the following result in analogous to \cite[Lemma 4.2.2\,(ii)]{qin2017triangular}.
\begin{Lemma}\label{lem:variation_X_variables_similar}
For any $i\in I$, the variation map $\var_{t'}^X$ sends the cluster $X$-variable $(\sigma\seq)^* X_i(s')$ to a Laurent monomial of the form $(\seq)^* X(s)^n$, $n\in \Z^I$.
\end{Lemma}
\begin{proof}
By \eqref{eq:commute_X_variation_seq}, we have $\seq^*\var_{s'}^X =\var_{t'}^X(\sigma \seq)^*$, where $\var_{s'}^X$ is also a variation map. In the seed $s'$, the variation map $\var_{s'}^X$ sends the $X$-variable $X_i(s')$ to some Laurent monomial $X(s)^n$ by definition. The desired statement follows by applying mutations $(\sigma \seq)^*$ and $\seq^*$.
\end{proof}

We have the following result in analogous to Lemma \ref{lem:A_variation_restrict_cluster_alg}.
\begin{Lemma}\label{lem:X_variation_restrict_cluster_alg}
 A variation map $\var_{t'}^{X}$ restricts to an algebra homomorphism from $\XClAlg(t')$\linebreak to~$\XClAlg(t)$.
 \end{Lemma}
 \begin{proof}
The proof is similar to Lemma \ref{lem:A_variation_restrict_cluster_alg}, where we replace \cite[Lemma 4.2.2\,(ii)]{qin2017triangular} by Lem\-ma~\ref{lem:variation_X_variables_similar}.
 \end{proof}

 We have extended the variation map $\var^X_{t'}$ to a homomorphism from \smash{$\cF_{\ufv}^X(t')$} to \smash{$\cF_{\ufv}^X(t)$}. We can also extend it to a homomorphism from \smash{$\hLP^X(t')$} to \smash{$\hLP^X(t)$}. Let \smash{$\cF_{\ufv}^X(t')\cap \hLP^X(t')$} denote the set of formal Laurent series \smash{$Z\in \hLP^X(t')$} such that $Z*Q=P$ for some $P\in \kk[(X')_i]_{i\in I}$, $Q\in \kk[(X')_k]_{k\in \Iufv}$. Then $Z$ is identified with the element $P*Q^{-1}$ in $\cF_{\ufv}^X(t')$.
 \begin{Proposition}
 The above two extended maps restrict to the same homomorphism from \smash{$\cF_{\ufv}^X(t')\cap \hLP^X(t')$} to \smash{$\cF_{\ufv}^X(t)\cap \hLP^X(t)$}.
 \end{Proposition}
 \begin{proof}
We denote the above two extended maps by \smash{$f\colon\cF_{\ufv}^X(t')\rightarrow \cF_{\ufv}^X(t)$} and $\smash{g\colon\hLP^X(t')\rightarrow} \allowbreak\smash{\hLP^X(t)}$, respectively. Take any element \smash{$Z\in \cF_{\ufv}^X(t')\cap \hLP^X(t')$} described as above. Then we have $g(Z)*g(Q)=g(P)$. In addition, we have $f(P)=\var_{t'}^X(P)=g(P)\in \kk[X_i]_{i\in I}$ and
 \[f(Q)=\var_{t'}^X(Q)=g(Q)\in \kk[X_k]_{k\in \Iufv}.\]
 It follows that $g(Z)*f(Q)=f(P)$, i.e, $g(Z)$ is contained in \smash{$\cF_{\ufv}^X(t)\cap \hLP^X(t)$} and it is identified with $f\big(P*Q^{-1}\big)$.
 \end{proof}

\subsection{Twist endomorphisms\label{subsec:Twist-automorphisms-X}}

Let $t'=\seq t$ be two seeds similar up to $\sigma$ as before. Note
that the canonical Poisson structure on $\LP^{X}(t)$ naturally gives rise
to a Poisson structure on $\XClAlg(t)$.

\begin{Definition}

Let $\var_{t'}^{X}\colon\LP^X(t')\rightarrow \LP^X(t)$ be any variation map.\footnote{If we work with $\var_{t'}^{N}\colon N(t')_{\Q}\rightarrow N(t)_{\Q}$,
the corresponding twist map should be defined on an algebra such that
we allow roots of all the cluster $X$-variables, which are roots
of rational functions. We do not consider this situation here.} Then the composition $\tw_{t'}^{X}=\big(\seq^{-1}\big)^{*}\circ\var_{t'}^{X}$
is called a twist endomorphism on $\XClAlg(t')$ passing through $t$.

For classical cases, if $\var_{t'}^{X}$ is Poisson, then $\tw_{t'}^{X}$ is said to be Poisson.

\end{Definition}

Note that, by Proposition \ref{prop:change_seed_X_variation}, when
$\tw_{t'}^{X}$ is a twist endomorphism, it gives rise to twist endomorphisms
$\tw_{s'}^{X}$, $s'\in\Delta^{+}$. They are identified via mutations.
So we can simply denote them by $\tw^{X}$.

Assume that $\var_{t'}^{X}$ is an invertible variation map. We have
a twist automorphism $\tw_{t'}^{X}:=\big(\seq^{-1}\big)^{*}\var_{t'}^{X}$
on $\XClAlg(t')$. Lemma \ref{lem:inverse_N_variation} implies that
\smash{$\bigl(\var_{t'}^{X}\bigr)^{-1}$} is still a variation map. So $\eta:=\allowbreak\smash{\seq^{*}\bigl(\var_{t'}^{X}\bigr)^{-1}}$
is a twist automorphism on $\XClAlg(t)$.

\begin{Lemma}\label{lem:inverse_X_twist_automorphism}
We have $\seq^{*}\bigl(\tw_{t'}^{X}\bigr)^{-1}(\seq^{*})^{-1}=\eta$. Namely,
by identifying the fraction fields by mutations, the twist automorphisms
associated to invertible variation maps $\var_{t'}^{X}\colon\LP^{X}(t')\rightarrow\LP^{X}(t)$
and $(\var_{t'}^{X})^{-1}\colon\LP^{X}(t)\rightarrow\LP^{X}(t')$ are inverse
to each other.
\end{Lemma}

\begin{proof}
The proof is similar to that of Lemma \ref{lem:inverse_A_twist_automorphism}.
\end{proof}

\begin{Theorem}\label{thm:twist_X}
Consider the classical case $\kk=\Z$. Assume that we have a linear map $\var_{t'}^{N}\colon N(t')\rightarrow N(t)$.
Let $\var_{t}^{M}$ denote its pullback from $\Mc(t)_{\Q}$ to $\Mc(t')_{\Q}$ and $\var_{t'}^X$, $\var_t^A$ denote the associated monomial maps.
Define $\tw_{t'}^{X}=\big(\seq^{-1}\big)^{*}\circ\var_{t'}^{X}$ and $\tw_{t}^{A}=\seq^{*}\circ\var_{t}^{A}$.
\begin{itemize}\itemsep=0pt
\item[$(1)$] $\tw_{t'}^{X}$ is a twist endomorphism on $\XClAlg(t')$
if and only if $\tw_{t}^{A}$ is a twist endomorphism on~$\upClAlg(t)_{\fFacRoot}$.
\item[$(2)$] If $\tw_{t'}^{X}$ is a twist automorphism, then $\bigl(\tw_{t'}^{X}\bigr)^{-1}$
is also a twist automorphism by Lemma~{\rm \ref{lem:inverse_N_variation}}.
Let \smash{$\bigl(\tw_{t}^{X}\bigr)^{-1}:=\seq^{*}\bigl(\tw_{t'}^{X}\bigr)^{-1}(\seq^{*})^{-1}$}
denote the corresponding twist automorphism on $\XClAlg(t)$. Then
\smash{$\bigl(\tw_{t}^{X}\bigr)^{-1}$} preserves the Poisson structure if and only
if \smash{$p^{*}\bigl(\tw_{t}^{X}\bigr)^{-1}=\tw_{t}^{A}p^{*}$}.
\end{itemize}
\end{Theorem}
Because we can quantized a Poisson variation map to a quantum variation map Remark~\ref{rem:X-Poisson_quantum_variation}, by Theorem \ref{thm:twist_X}\,(2), the equality $p^{*}\bigl(\tw_{t}^{X}\bigr)^{-1}=\tw_{t}^{A}p^{*}$ provides a criterion for $\bigl(\tw_{t}^{X}\bigr)^\pm$ being quantum twist automorphisms for quantum cases. This equality can be presented by the quadratic equation in Lemma \ref{lem:commute_variation_p}.

\begin{proof}[Proof of Theorem \ref{thm:twist_X}]

(1) The statement follows from Proposition \ref{prop:linear_N_variation}.

(2) On the one hand, we have $\tw_{t}^{A}p^{*}=\seq^{*}\var_{t}^{A}p^{*}$.
On the other hand, we have \[p^{*}\bigl(\tw_{t}^{X}\bigr)^{-1}=p^{*}\seq^{*}\bigl(\tw_{t'}^{X}\bigr)^{-1}(\seq^{*})^{-1}=p^{*}\seq^{*}\bigl(\var_{t'}^{X}\bigr)^{-1}\]
by construction, which equals $\seq^{*}p^{*}\bigl(\var_{t'}^{X}\bigr)^{-1}$.
Moreover, Lemma \ref{lem:commute_variation_p} implies that $\var_{t}^{A}p^{*}=\allowbreak\smash{p^{*}\bigl(\var_{t'}^{X}\bigr)^{-1}}$
if and only if $\var_{t'}^{N}$ preserves $\omega$. The claim follows.
\end{proof}

\begin{Remark}[cluster Poisson algebras without coefficients]
Let $t'=\seq t$ be seeds similar up to $\sigma$. Let us restrict
to the unfrozen cluster Poisson algebra $\XClAlg_{\ufv}$, which consists of $P*Q^{-1}$ with~$P,Q\in \kk[X_k]_{k\in \Iufv}$. Then there is only one twist automorphism $\tw_{t}^{X}$ on $\XClAlg_{\ufv}(t)$
passing through $t'$. It sends $X_{k}(t)$ to $\seq^{*}X_{\sigma k}(t')$
and it preserves the Poisson structure when $\kk=\Z$.

Assume $I=\Iufv$. Take any twist endomorphism $\tw_{t}^{A}$ on $\upClAlg(t)$
passing through $t'$. Then we have $\tw_{t}^{A}p^{*}X_{k}(t)=\seq^{*}p^{*}X_{\sigma k}(t')$
for any $k\in I_{\ufv}$ by \eqref{eq:preserve_p_map}. Since $\seq^{*}p^{*}=p^{*}\seq^{*}$,
it follows that \[\tw_{t}^{A}p^{*}X_{k}(t)=\seq^{*}p^{*}X_{\sigma k}(t')=p^{*}\seq^{*}X_{\sigma k}(t')=p^{*}\tw_{t}^{X}X_{k}(t),\]
i.e., Theorem \ref{thm:twist_X}\,(2) always hold true.
\end{Remark}

\begin{Remark}[{existence of $\tw^X$}]\label{rem:exist_X_twist}
If all skew-symmetrizers $d_{i}$ are equal, we can work with lattices
$N$ and $\Mc$ instead of $\Q$-vector spaces in Lemma~\ref{lem:pullback_matrix}
and Proposition~\ref{prop:linear_N_variation}.

By using Proposition \ref{prop:linear_N_variation}, we can describe
the set of all linear variation maps $\var_{t'}^{N}\colon N(t')_{\Q}\rightarrow N(t)_{\Q}$
as in Proposition~\ref{prop:exist_linear_variation}\,(1)--(2).

Recall that $W$ denotes the matrix of $\omega$. If \smash{$\left(\begin{smallmatrix}
W_{\ufv}\\
W_{\low}
\end{smallmatrix}\right)$} is $\Z$-valued and an order $|I_{\ufv}|$ minor of it equals $\pm1$,
we can describe the set of all linear variation maps $\var_{t'}^{N}\colon N(t')\rightarrow N(t)$
as in Proposition~\ref{prop:exist_linear_variation}\,(3). In the classical case $\kk=\Z$, these linear variation
maps give rise to twist endomorphisms $\tw^{X}$ on $\XClAlg$ in
analogous to Theorem~\ref{thm:twist_A}\,(2).
\end{Remark}

\section{Twist automorphisms in special cases}\label{sec:Twist-automorphisms-special-cases}

\subsection{Twist automorphism of Donaldson--Thomas type}\label{subsec:Cluster-twist-automorphism-DT}

Assume that $t$ is injective-reachable such that $t'=t[1]$ is similar
to $t$ up to a permutation~$\sigma$. Then a twist automorphism on
$t$ passing through $t[1]$ is called a twist automorphism of Donaldson--Thomas
type\footnote{The mutation maps $\seq^{*}$ from $\cF^{A}(t[1])$ to $\cF^{A}(t)$
is often called a Donaldson--Thomas transformation in cluster theory.
See \cite{Nagao10} for the relation between Donaldson--Thomas theory
and cluster algebras.} (DT-type for short).

\begin{Proposition}\label{prop:exist_var_DT}
 There exists an invertible Poisson variation map $\var_{t}^{N}$ and an invertible variation map $\var_{t}^{M}$
passing through $t[1]$. Moreover, when $\cT^A(t)$ has a compatible Poisson structure, we can choose $\var_{t}^{M}$ to be Poisson.
\end{Proposition}

\begin{proof}

Relabel the vertices for $t[1]$ so that we can assume $\sigma=\Id$.
By \eqref{eq:injective_reachable_def}, we have $\smash{\psi_{t,t[1]}^{N}\uee'}=\allowbreak\smash{\uee\left(\begin{smallmatrix}
-\Id_{\ufv} & E_{\high}\\
0 & \Id_{\fv}
\end{smallmatrix}\right)=:\uee E}$ and, correspondingly, \smash{$\psi^M_{t,t[1]}\uff'=\uff\left(\begin{smallmatrix}
-\Id_{\ufv}&0\\
F_{\low} & \Id_{\fv}
\end{smallmatrix}\right)=:\uff F$}. Note that \smash{$\psi^M_{t,t[1]}$} is the identity on $\Mc_{\fv}$.

(1) Let us define the linear isomorphism $\var_{t}^{M}\colon\Mc(t)\rightarrow\Mc(t')$
such that
\begin{align*}
\var_{t}^{M}(\uff)=\uff'\begin{pmatrix}
\Id_{\ufv}&0\\
-F_{\low} & -\Id_{\fv}
\end{pmatrix}.
\end{align*}
Then we have $\psi_{t,t[1]}^{M}\var_{t}^{M}(\uff)=-\uff$.

On the one hand, for any $k\in I$, we have
\begin{gather*}
\psi_{t,t[1]}^{M}\var_{t}^{M}(p^{*}\ee_{k}) = \psi_{t,t[1]}^{M}\var_{t}^{M}\bigg(\sum_{i\in I}b_{ik}\ff_{i}\bigg) = -\sum_{i\in I}b_{ik}\ff_{i}= -p^{*}\ee_{k}.
\end{gather*}
On the other hand, by Lemma \ref{lem:bilinear_form_basis_change},
for any $k\in\Iufv$, we have
\begin{gather*}
\psi_{t,t[1]}^{M}p^{*}(\ee_{k}') = p^{*}\big(\psi_{t,t[1]}^{N}\ee_{k}'\big) = p^{*}(-\ee_{k}) = -p^{*}\ee_{k}.
\end{gather*}
It follows that $\var_{t}^{M}(p^{*}\ee_{k})=p^{*}(\ee_{k}')$ for
any $k\in\Iufv$. In particular, $\var_{t}^{M}$ is a variation map.

Assume there is a compatible Poisson structure $\lambda$. For any $i,j\in I$, we have
\begin{gather*}
\lambda\bigl(\var_{t}^{M}\ff_{i},\var_{t}^{M}\ff_{j}\bigr) = \lambda\big(\psi_{t,t[1]}^{M}\var_{t}^{M}\ff_{i},\psi_{t,t[1]}^{M}\var_{t}^{M}\ff_{j}\big)= \lambda(-\ff_{i},-\ff_{j}) = \lambda(\ff_{i},\ff_{j}).
\end{gather*}
Therefore, $\var_{t}^{M}$ preserves $\lambda$.

(2) Define $\var_{t}^{N}\colon N(t)\rightarrow N(t[1])$ such that \smash{$\var_{t}^{N}(\uee)=\uee'\left(\begin{smallmatrix}
\Id_{\ufv} & -E_{\high}\\
0 & -\Id_{\fv}
\end{smallmatrix}\right)$}. It follows that \[\psi_{t[1],t}^{N}\var_{t}^{N}(\uee)=-\uee.\] Then,
for any $i,j\in I$, we have
\begin{gather*}
\omega\bigl(\var_{t}^{N}\ee_{i},\var_{t}^{N}\ee_{j}\bigr) =\omega\big(\psi_{t,t[1]}^{N}\var_{t}^{N}\ee_{i},\psi_{t,t[1]}^{N}\var_{t}^{N}\ee_{j}\big) =\omega(-\ee_{i},-\ee_{j})=\omega(\ee_{i},\ee_{j}).
\end{gather*}
Therefore, $\var_{t}^{N}$ preserves $\omega$.
\end{proof}

\begin{Theorem}\label{thm:twist_DT}
Assume that $t$ is injective-reachable such that $t'=t[1]$ is similar
to $t$ up to a~permutation $\sigma$.
\begin{itemize}\itemsep=0pt
\item[$(1)$] There exists a twist automorphism $\tw_{t}^{A}$ on $\upClAlg(t)$
of Donaldson--Thomas type.
\item[$(2)$] There exists a twist automorphism $\tw_{t}^{X}$ on $\XClAlg(t)$
of Donaldson--Thomas type.
\item[$(3)$] We can choose $\tw_{t}^{A}$ and $\tw_{t}^{X}$ in $(1)$--$(2)$ such that $p^{*}\tw_{t}^{X}=\tw_{t}^{A}p^{*}$.
\end{itemize}
\end{Theorem}

\begin{proof}
(1)--(2) The statements follow from Proposition \ref{prop:exist_var_DT}.

(3) We claim that the two variation maps in the proof of Proposition
\ref{prop:exist_var_DT} are related by the pullback construction.
The variation map $\var_{t}^{M}$ is represented by \smash{$\left(\begin{smallmatrix}
\Id_{\ufv} & 0\\
-F_{\low} & -\Id_{\fv}
\end{smallmatrix}\right)$}. Then its pullback $\bigl(\var_{t}^{M}\bigr)^{*}$ satisfies \[\bigl(\var_{t}^{M}\bigr)^{*}(\uee')=\uee D^{-1}\begin{pmatrix}
\Id_{\ufv} & -F_{\low}^{\rm T}\\
0 & -\Id_{\fv}
\end{pmatrix}D\] by Lemma \ref{lem:pullback_matrix}. By \eqref{eq:E_D_F_relation},
$E=D^{-1}F^{\rm -T}D$. It follows that $E_{\high}=D_{\ufv}^{-1}F_{\low}^{\rm T}D_{\fv}$.
Define $\var_{t}^{N}\colon N(t)\rightarrow N(t[1])$ to be the inverse of
$\bigl(\var_{t}^{M}\bigr)^{*}$. We obtain
\begin{align*}
\var_{t}^{N}(\uee) & =\uee'D^{-1}\begin{pmatrix}
\Id_{\ufv} & -F_{\low}^{\rm T}\\
0 & -\Id_{\fv}
\end{pmatrix}D =\uee'\begin{pmatrix}
\Id_{\ufv} & -E_{\high}\\
0 & -\Id_{\fv}
\end{pmatrix}
\end{align*}
as in the previous proof.

The desired claim follows from Theorem \ref{thm:twist_X}\,(2).
\end{proof}

\begin{Example}\label{eg:A_1_variation_map}
Let us continue Example \ref{eg:A_1_eg}. Then we have $t'=t[1]$.

We can choose the variation map $\var_{t}^{M}\colon\Mc(t)\rightarrow\Mc(t[1])$
such that it is represented by the matrix $\left(\begin{smallmatrix}
\hphantom{-}1 & \hphantom{-}0\\
-1 & -1
\end{smallmatrix}\right)$. Then $\psi_{t,t[1]}^{M}\var_{t}^{M}$ is $-\Id$ on $\Mc(t)$.

We can choose the variation map $\var_{t}^{N}\colon N(t)\rightarrow N(t[1])$
such that it is represented by the matrix $\left(\begin{smallmatrix}
1 & -2\\
0 & -1
\end{smallmatrix}\right)$. Then $\psi_{t,t[1]}^{N}\var_{t}^{N}$ is $-\Id$ on $N(t)$.

\end{Example}

We refer the reader to Example \ref{eg:twist_auto_matrix} for an
alternative example.

\subsection{Twist automorphism for principal coefficients\label{subsec:Cluster-twist-automorphism-principal}}

Let $t_{0}$ denote an initial seed with principal coefficients, see
Section \ref{subsec:Principal-coefficients}. We denote its $B$-matrix
by \smash{$B=\left(\begin{smallmatrix}
B_{\ufv} & -\Id_{\ufv}\\
\Id_{\ufv} & 0
\end{smallmatrix}\right)$}. Recall that \smash{$D=\left(\begin{smallmatrix}
D_{\ufv}\\
 & D_{\ufv}
\end{smallmatrix}\right)$} in this case. We endow $t_{0}$ with the canonical Poisson structure whose $\Lambda$-matrix is
\[
\Lambda:= B^{\rm -T}D=\left(\begin{matrix}
0 & -\Id_{\ufv}\\
\Id_{\ufv} & B_{\ufv}^{\rm T}
\end{matrix}\right)D.
\] Then, for any seed $t'\in\Delta^{+}$, we can check
that $\Lambda(t')= B(t')^{\rm -T}D$ by using Proposition~\ref{prop:matrix_properties}
recursively.

Assume that $t\in\Delta_{t_0}^{+}$ is similar to the principal coefficients
seed $t_{0}$ up to a permutation $\sigma$.

\begin{Theorem}\label{thm:twist_principal_coeff}\quad
\begin{itemize}\itemsep=0pt
\item[$(1)$] There exists a twist automorphism $\tw_{t_{0}}^{A}$ on the $($quantum$)$ upper cluster algebra on
$\upClAlg(t_{0})$ passing through $t$.
\item[$(2)$] There exists a twist automorphism $\tw_{t_{0}}^{X}$ on
$\XClAlg(t_{0})$ passing through $t$.
\item[$(3)$] We can choose $\tw_{t_{0}}^{A}$ and $\tw_{t_{0}}^{X}$ such that
$p^{*}\tw_{t_{0}}^{X}=\tw_{t_{0}}^{A}p^{*}$.
\end{itemize}
\end{Theorem}

\begin{proof}
Relabeling the vertices in $t$ if necessary, we can assume $\sigma=\Id$.
Since $t_{0}$ has principal coefficients, we have \smash{$B(t)=\left(\begin{smallmatrix}
B_{\ufv} & -G^{-1}\\
C & 0
\end{smallmatrix}\right)$}, where $C$ and $G$ denote the $C$-matrix and the $G$-matrix of
the seed $t$ with respect to the initial seed $t_{0}$, respectively,
and we have $C^{\rm T}D_{\ufv}=D_{\ufv}G^{-1}$ (see Section~\ref{subsec:Principal-coefficients}).

(1) Define the linear map $\var_{t_{0}}^{M}\colon\Mc(t_{0})\rightarrow\Mc(t)$
such that{\samepage
\begin{align*}
\var_{t_{0}}^{M}(\uff(t_{0})) =\uff(t)\begin{pmatrix}
\Id_{\ufv} & 0\\
0 & C
\end{pmatrix}.
\end{align*}
It is clear that $\left(\begin{smallmatrix}
0 & C\end{smallmatrix}\right)\left(\begin{smallmatrix}
B_{\ufv}\\
\Id
\end{smallmatrix}\right)=C$. So $\var_{t}^{M}$ is a variation map by Lemma \ref{lem:matrix_eq_A_variation}.}

We have \smash{$\Lambda(t)= B(t)^{\rm -T}D=\left(\begin{smallmatrix}
0 & -G^{\rm T}\\
C^{\rm -T} & C^{\rm -T}B_{\ufv}^{\rm T}G^{\rm T}
\end{smallmatrix}\right)D$}. Recall that $\var_{t_{0}}^{M}$ is represented by $\left(\begin{smallmatrix}
\Id_{\ufv} & 0\\
0 & C
\end{smallmatrix}\right)$. We have
\begin{align*}
\begin{pmatrix}
\Id & 0\\
0 & C
\end{pmatrix}^{\rm T}\Lambda(t)\begin{pmatrix}
\Id & 0\\
0 & C
\end{pmatrix} & =\begin{pmatrix}
\Id & 0\\
0 & C^{\rm T}
\end{pmatrix}\begin{pmatrix}
0 & -G^{\rm T}\\
C^{\rm -T} & C^{\rm -T}B_{\ufv}^{\rm T}G^{\rm T}
\end{pmatrix}D\begin{pmatrix}
\Id & 0\\
0 & C
\end{pmatrix}\\
 & =\begin{pmatrix}
0 & -G^{\rm T}D_{\ufv}C\\
D_{\ufv} & B_{\ufv}^{\rm T}G^{\rm T}D_{\ufv}C
\end{pmatrix} =\begin{pmatrix}
0 & -D_{\ufv}\\
D_{\ufv} & B_{\ufv}^{\rm T}D_{\ufv}
\end{pmatrix} =\Lambda.
\end{align*}
Therefore, $\var_{t_{0}}^{M}$ preserves $\lambda$. The associated
twist automorphism provides a desired solution.

(2) Let $\bigl(\var_{t_{0}}^{M}\bigr)^{*}\colon N(t)\rightarrow N(t_{0})$ denote
the pullback of $\var_{t_{0}}^{M}$ in the proof of (1). Then $\bigl(\var_{t_{0}}^{M}\bigr)^{*}(\uee(t))=\uee(t_{0})D^{-1}\left(\begin{smallmatrix}
I & 0\\
0 & C^{\rm T}
\end{smallmatrix}\right)D$. Define $\var_{t_{0}}^{N}$ as the inverse:
\begin{align*}
\var_{t_{0}}^{N}(\uee(t_{0})) & =\uee(t)D^{-1}\begin{pmatrix}
I & 0\\
0 & C^{\rm T}
\end{pmatrix}^{-1}D =\uee(t)\begin{pmatrix}
\Id_{\ufv}\\
 & G
\end{pmatrix}.
\end{align*}

We have \smash{$W=B^{\rm T}D=\left(\begin{smallmatrix}
B_{\ufv}^{\rm T} & \Id\\
-\Id & 0
\end{smallmatrix}\right)D$} and \smash{$W(t)=B(t)^{\rm T}D=\left(\begin{smallmatrix}
B_{\ufv}^{\rm T} & C^{\rm T}\\
-G^{\rm -T} & 0
\end{smallmatrix}\right)D$}.

Using $D_{\ufv}G=C^{\rm -T}D_{\ufv}$, we obtain
\begin{align*}
\begin{pmatrix}
\Id_{\ufv} & 0\\
0 & G
\end{pmatrix}^{\rm T}W(t)\begin{pmatrix}
\Id_{\ufv}\\
 & G
\end{pmatrix} & =\begin{pmatrix}
\Id_{\ufv} & 0\\
0 & G^{\rm T}
\end{pmatrix}\begin{pmatrix}
B_{\ufv}^{\rm T} & C^{\rm T}\\
-G^{\rm -T} & 0
\end{pmatrix}D\begin{pmatrix}
\Id_{\ufv}\\
 & G
\end{pmatrix}\\
 & =\begin{pmatrix}
\Id_{\ufv} & 0\\
0 & G^{\rm T}
\end{pmatrix}\begin{pmatrix}
B_{\ufv}^{\rm T} & C^{\rm T}\\
-G^{\rm -T} & 0
\end{pmatrix}\begin{pmatrix}
\Id_{\ufv}\\
 & C^{\rm -T}
\end{pmatrix}D\\
 & =\begin{pmatrix}
B_{\ufv}^{\rm T} & \Id\\
-\Id & 0
\end{pmatrix}D =W.
\end{align*}
Therefore, $\var_{t_{0}}^{N}$ preserves $\omega$. The associated
twist automorphism provides a desired solution.

(3) We choose twist automorphisms as in the proof of (1)--(2). Then
the claim follows from Theorem \ref{thm:twist_X}\,(2).
\end{proof}

\begin{Example}
We can verify that $\var_{t_{0}}^{N}$ and $\var_{t_{0}}^{M}$ constructed
in the proof of Theorem~\ref{thm:twist_principal_coeff} make the
diagram~\eqref{eq:commute_variation_p} commutative, which is equivalent
to the following:
\begin{align}
\begin{array}{@{}ccc}
\Mc(t_{0}) & \xrightarrow{\var_{t_{0}}^{M}} & \Mc(t)\\
\uparrow p^{*} & & \uparrow p^{*}\\
N(t_{0}) & \xrightarrow{\var_{t_{0}}^{N}} & N(t).
\end{array}\label{eq:commute_variation_p_principal}
\end{align}

More precisely, $\var_{t_{0}}^{M}p^{*}\colon N(t_{0})\rightarrow\Mc(t)$
is represented by the matrix $\left(\begin{smallmatrix}
\Id_{\ufv} & 0\\
0 & C
\end{smallmatrix}\right)B(t_{0})=\left(\begin{smallmatrix}
B_{\ufv} & -\Id\\
C & 0
\end{smallmatrix}\right)$, and $p^{*}\var_{t_{0}}^{N}$ is also represented by the matrix \smash{$B(t)\left(\begin{smallmatrix}
\Id_{\ufv} & 0\\
0 & G
\end{smallmatrix}\right)=\left(\begin{smallmatrix}
B_{\ufv} & -\Id\\
C & 0
\end{smallmatrix}\right)$}.

Equivalently, $\var_{t_{0}}^{N}$ and $\var_{t_{0}}^{M}$ are determined
by each other such that the diagram \eqref{eq:commute_variation_p_principal}
is commutative.
\end{Example}

\section{Permutations on bases}\label{sec:Permutation-on-bases}

\subsection[Bases for U\^{}A]{Bases for $\boldsymbol{\upClAlg}$}\label{sec:permute_A_bases}

For any seed $t\in \Delta^+$, assume that we have chosen a collection of polynomials for $m\in \Z^I$: \[F_{m}(t)=\sum_{n\in\N^{\Iufv}}c_{n}(m;t)X(t)^{n}\]
with $c_{0}(m;t)=1$, $c_{n}(m;t)\in\kk$. Define the corresponding subset $\cS(t)=\big\{S_{m}(t)\mid m\in\Z^I\big\}$
of~$\LP^{A}(t)$, such that the elements $S_{m}$ are given by
\begin{align*}S_{m}(t)=A(t)^{m}\cdot F_{m}(t)|_{X(t)^n\mapsto A(t)^{p^{*}(t)n}}.\end{align*}
We often omit the symbol $t$ for simplicity. $F_m(t)$ are called the $F$-polynomials of $S_m(t)$. Note that we use the commutative product $\cdot$ here, so that the $F$-polynomials are defined as in standard literature.

Note that, if the full rank assumption (see Assumption \ref{fullrank}) holds, $S_m(t)$ are $m$-pointed and~$\cS(t)$ is
be $\kk$-linearly independent. Otherwise, $\cS(t)$ is not necessarily linearly independent and different choices of $F_m(t)$ might produce the same $S_m(t)$.

\begin{Assumption}\label{assumption:A_basis}
 We assume the following conditions hold.
\begin{enumerate}\itemsep=0pt
\item[(1)] The sets $\cS(t)$, $t\in\Delta^{+}$, are identified by mutations: $\forall t'=\seq t$, we have $\seq^{*}\cS(t')=\cS(t)$.

\item[(2)] The chosen functions $F_{m}(t)$ only depend on $\pr_{\Iufv}m$ and $B(t)_{\ufv}$:
\begin{itemize}\itemsep=0pt
\item If $m'=m+u$ for some $u\in\Mc_{\fv}$, then $F_{m'}(t)=F_{m}(t)$.
\item If $t'$ and $t$ are similar up to $\sigma$ and $\pr_{\Iufv}m'=\sigma\pr_{\Iufv}m$
for some $m,m'\in\Z^{I}$, we have $c_{\sigma n}(m';t')=c_{n}(m;t)$,
$\forall n$.
\end{itemize}
\end{enumerate}
\end{Assumption}
Note that Assumption \ref{assumption:A_basis} implies $\cS(t)\subset\upClAlg(t)$.

We have the following important examples of $\cS(t)$ satisfying
Assumption~\ref{assumption:A_basis}:
\begin{itemize}\itemsep=0pt
\item The set of generic cluster characters: When $B(t)_{\ufv}$ is skew-symmetric
and the corresponding cluster category satisfies some finiteness condition,
we can construct the set of generic cluster characters $\gen(t)=\{\gen_{m}(t)|m\in\Mc(t)\}$
in $\upClAlg$ using algebra representation theory. By~\cite{plamondon2013generic},
they satisfy the above assumptions.
\begin{itemize}\itemsep=0pt
\item By \cite{qin2019bases}, when $t$ is injective-reachable, $\gen(t)$
is a basis of $\upClAlg$, called the generic basis in the sense of
\cite{dupont2011generic}. This family of bases includes the dual
semicanonical basis \cite{lusztig2000semicanonical} of a unipotent
cell with symmetric Cartan datum \cite{GeissLeclercSchroeer10b}.
\end{itemize}
\item The common triangular basis \cite{qin2017triangular}: this basis
is known to satisfy the above assumptions~\cite{qin2017triangular}.
Its existence is known for some cluster algebras with a Lie theoretic
background. This family of bases includes the dual canonical basis
of a quantum unipotent cell \cite{kashiwara2019laurent,qin2017triangular,qin2020dual}.
\item The set of theta functions \cite{gross2018canonical}: the set of
the theta functions, under the assumption that they are Laurent polynomials,\footnote{In general, the theta functions are formal Laurent series, see \cite{qin2019bases}
for the formal completion. Our arguments for variation maps remain
valid for formal Laurent series. For simplicity, we only treat Laurent
polynomials in this paper.} satisfies the above assumption. They
are known to form a~basis of $\upClAlg(t)$ under appropriate conditions,
for example, when $t$ is injective-reachable.
\end{itemize}
\begin{Theorem}\label{thm:twist_A_basis}
Let $\tw^{A}$ denote any twist endomorphism on the $($quantum$)$ upper cluster algebra~$\upClAlg$. Assume
the sets $\cS(t)$, $t\in \Delta^+$, satisfy Assumption {\rm \ref{assumption:A_basis}}, then the following statements are true.
\begin{itemize}\itemsep=0pt
\item[$(1)$] We have $\tw^{A}\cS(t)\subset\cS(t)$.
\item[$(2)$] If $\tw^{A}$ is a twist automorphism, we have $\tw^{A}\cS(t)=\cS(t)$.
\end{itemize}
\end{Theorem}

\begin{proof}
Recall that $\tw^{A}=\seq^*\var_{t}^{A}$ for some $t'=\seq t$ and some variation map $\var_{t}^{A}$ from $\LP^A(t)$ to $\LP^A(t')$.

(1) Denote $\cS(t)=\big\{S_m\mid m\in \Z^I\big\}$ and $\cS(t')=\big\{S'_{m'}\mid m'\in \Z^I\big\}$ for simplicity. Take any \smash{$S_m=A^m \cdot F_m|_{X^n\mapsto A^{p^* n}}$} from $\cS(t)$, where $F_m =\sum c_n(m;t)X^n$. Then $\smash{\var_{t}^{A}(S_m)=A^{\var_t^M m} }\cdot \smash{F'|_{(X')^n\mapsto (A')^{p^*(t') n}}}$, where $F'=\sum c_{n}(m;t)(X')^{\sigma n}$. Note that $\pr_{\Iufv}\var_t^M m=\sigma \pr_{\Iufv}m$. By Assumption \ref{assumption:A_basis}\,(2), we have $F'=F_{\var_t^M m}(t')$. So $\var_{t}^{A}(S_m)=S'_{\var^M_t m}$.

By Assumption \ref{assumption:A_basis}\,(1), we have $\seq^* S'_{\var^M_t m}\in \cS(t)$. Therefore, $\tw^{A} S_m\in \cS(t)$.

(2) Recall that the inverse of $\tw_{t}^{A}$ is also a twist automorphism,
the claim follows.
\end{proof}

\subsection[Bases for U\^{}X]{Bases for $\boldsymbol{\XClAlg}$}

Similarly, for any $t$, assume that we have a subset $\cZ(t)=\big\{Z_{n}(t)\mid n\in \Z^I\big\}$
in $\LP^{X}(t)$, such that its elements $Z_{n}(t)$ are $n$-pointed.
Then we can denote $Z_{n}(t)=X(t)^{n}\cdot F_{n}(t)$, whose $F$-polynomials
take the form $F_{n}(t)=\sum_{e\in\N^{\Iufv}}c_{e}(n;t)X(t)^{e}$ with
$c_{0}(n;t)=1$, $c_{e}(n;t)\in\kk$.

Note that $\cZ(t)$ must be $\kk$-linearly independent.

\begin{Assumption}\label{assumption:X_basis}
 We assume the following conditions hold.
\begin{enumerate}\itemsep=0pt
\item[(1)] The sets $\cZ(t)$, $\forall t\in\Delta^{+}$, are identified by mutations: for any $t'=\seq t$, we have $\seq^{*}\cZ(t')=\cZ(t)$.
\item[(2)] The $F$-polynomials $F_{n}(t)$ only depend on $\pr_{\Iufv}p^{*}n$
and $B(t)_{\ufv}$:
\begin{itemize}\itemsep=0pt
\item If $p^{*}n'=p^{*}n+u$ for some $u\in\Mc_{\fv}$, then $F_{n'}(t)=F_{n}(t)$.
\item If $t'$ and $t$ are similar up to $\sigma$ and $\pr_{\Iufv}B(t')n'=\sigma\pr_{\Iufv}B(t)n$
for some $n,n'\in\Z^{I}$, we have $c_{\sigma e}(n';t')=c_{e}(n;t)$,
$\forall e$.
\end{itemize}
\end{enumerate}
\end{Assumption}

Assumption \ref{assumption:X_basis}\,(1) implies $\cZ(t)\subset\XClAlg(t)$.

\begin{Theorem}\label{thm:twist_X_basis}
Let $\tw^{X}$ denote any twist endomorphism on $\XClAlg$. Assume
the sets $\cZ(t)$, $\forall t\in\Delta^{+}$, satisfy Assumption
{\rm \ref{assumption:X_basis}}. Then the following statements are true.
\begin{itemize}\itemsep=0pt
\item[$(1)$] We have $\tw^{X}\cZ(t)\subset\cZ(t)$.
\item[$(2)$] If $\tw^{X}$ is a twist automorphism, we have $\tw^{X}\cZ(t)=\cZ(t)$.
\end{itemize}
\end{Theorem}

\begin{proof}
The proof is
similar to that of Theorem \ref{thm:twist_X_basis}.
\end{proof}

\subsection[Construction of bases for U\^{}X]{Construction of bases for $\boldsymbol{\XClAlg}$}\label{sec:construction_X_bases}

There is few literature on the bases for $\XClAlg$. Let us explain
how to construct a basis for $\XClAlg$ from a basis for $\upClAlg$
following the idea of \cite{gross2018canonical}.

Let $t^{\prin}$ denote the principal coefficient seed constructed
from $t$ by adding a copy ${I'\!=\!\{i'|i\in I\}}$ of $I$ as extra frozen
variables, such that \smash{$B\big(t^{\prin}\big)=\left(\begin{smallmatrix}
B(t) & -\Id\\
\Id & 0
\end{smallmatrix}\right)$}. Assume that a set $\cS\big(t^{\prin}\big)\subset \upClAlg\big(t^{\prin}\big)$ as in Section \ref{sec:permute_A_bases} has been given. Then we
can construct a set $\cZ(t)\subset\LP^{X}(t)$ such that $Z_{n}(t)=X(t)^{n}\cdot F_{n}$
where $F_{n}$ is the $F$-polynomial of $S_{m}\big(t^{\prin}\big)$ for $m=B\big(t^{\prin}\big)\left(\begin{smallmatrix}
 n\\
 0
 \end{smallmatrix}\right)=\left(\begin{smallmatrix}
B(t)n\\
n
\end{smallmatrix}\right)$.

\begin{Lemma}\label{lem:X_set_univ_Laurent}
If the sets $\big(\seq^A\big)^{*}\cS\big(t^{\prin}\big)\subset\upClAlg(\seq^{-1} t^{\prin})$,
for all mutation sequences $\seq$, satisfy Assumption {\rm \ref{assumption:A_basis}},
then the sets \smash{$\cZ\big(\seq^{-1}t\big):=\big(\seq^X\big)^{*}\cZ(t)$}, $\forall\seq$, are contained in $\XClAlg(t)$ and satisfy Assumption
{\rm \ref{assumption:X_basis}}.
\end{Lemma}
\begin{proof}
We claim that $\seq^* \cZ(t)$ is contained in $\LP^X\big(\seq^{-1} t\big)$ for any mutation sequence $\seq$. If so, we deduce that the \smash{$\cZ(t)\subset \cap_{\forall \seq}(\seq^*)^{-1}\LP^X\big(\seq^{-1}t\big) =\XClAlg(t)$}. Then the sets $\cZ\big(\seq^{-1}t\big)$, $\forall \seq$, satisfy Assumption \ref{assumption:X_basis}\,(1) by definition. In addition, they satisfy Assumption \ref{assumption:X_basis}\,(2) because the sets~$\cS\big(\seq^{-1} t^{\prin}\big)$, for all $\seq$, satisfy Assumption~\ref{assumption:A_basis}\,(2).

Denote $s=t^{\prin}$, $t'=(\seq)^{-1} t$, and $s'=(\seq)^{-1} s$. Then we have the natural inclusions $\LP^X(t)\subset \LP^X(s)$ and $\LP^X(t')\subset \LP^X(s')$.

Note that the monomial maps \smash{$p^*(s)\colon\LP^X(s)\rightarrow \LP^A(s)$} and \smash{$p^*(s')\colon\hLP^X(s)\rightarrow \hLP^A(s)$} are injective. We also have injective maps \smash{$\iota^\bullet\colon\LP^\bullet(s)\rightarrow \hLP^\bullet (s')$} taking the formal Laurent series after mutations (see Section \ref{sec:mutation}).

Take any $Z\in \cZ(t)$. Then we have $p^*(s)Z=S_m\in \cS(s)$ for some $m$. So we get $\big(\seq^A\big)^{*}p^{*}(s)Z=S_m'\in \cS(s')$ for some $m'$. Therefore, $\iota^A(s)p^*(s)Z=S_{m'}$. But we also have $\iota^A(s)p^*(s)Z=p^*(s')\iota^X(s)Z$. So the formal series $\iota^X(s)Z$ is sent to a Laurent polynomial $S_{m'}$ by the injective map $p^*(s')$. It follows that $\iota^X(s)Z$ is a Laurent polynomial. Therefore, $\big(\seq^X\big)^*Z=\iota^X Z$ is contained in $\LP^X(s')$.
\end{proof}

Using this construction, by taking $\cS\big(t^{\prin}\big)$ as the set of
generic cluster characters, the common triangular basis, or the set
of the theta functions which are assumed to be Laurent polynomials, we obtain the corresponding
set $\cZ(t)\subset\XClAlg(t)$.

Finally, let us discuss when $\cZ(t)$ becomes a basis of $\XClAlg(t)$.

\begin{Theorem}\label{thm:construct_X_bases}
Under the assumption of Lemma {\rm \ref{lem:X_set_univ_Laurent}}, if $p^{*}\XClAlg\big(t^{\prin}\big)$ is contained in the free $\kk$-module
$\oplus_{m}\kk S_{m}\big(t^{\prin}\big)\subset\upClAlg\big(t^{\prin}\big)$, then
$\cZ(t)$ is a $\kk$-basis of $\XClAlg$(t). In particular, when
$\cS\big(t^{\prin}\big)$ is a basis of $\upClAlg\big(t^{\prin}\big)$, $\cZ(t)$
is a basis of $\XClAlg(t)$.
\end{Theorem}
\begin{proof}
On the one hand, we know that $\bigl\{\cS_{m}\big(t^{\prin}\big)\mid \exists n,\ m=\left(\begin{smallmatrix}
B(t)n\\
n
\end{smallmatrix}\right)\bigr\}$
 is contained in $p^{*}\XClAlg\big(t^{\prin}\big)$.

On the other hand, let us endow $\LP\big(t^{\prin}\big)$ with the $\Mc(t)$-grading
such that, $\forall i\in I$, $\deg A_{i}(t)=\ff_{i}$ and $\deg A_{i'}(t)=-\sum_{j\in I}b_{ji}\ff_{j}$.
Then $S_{m}\big(t^{\prin}\big)$ are homogeneous elements, and they have degree
$0$ if and only if $m=\left(\begin{smallmatrix}
B(t)n\\
n
\end{smallmatrix}\right)$ for some $n$. So $\bigl\{\cS_{m}\big(t^{\prin}\big)\mid \exists n,\ m=\left(\begin{smallmatrix}
B(t)n\\
n
\end{smallmatrix}\right)\bigr\}$ is a basis for the degree-$0$ subspace of $\oplus_{m}\kk S_{m}\big(t^{\prin}\big)$.
Moreover, since $p^{*}X_{i}\big(t^{\prin}\big)$ has degree $0$ for any $i\in I$,
$p^{*}\XClAlg\big(t^{\prin}\big)$ consists of degree $0$ elements.

Therefore, $\bigl\{\cS_{m}\big(t^{\prin}\big)\mid \exists n,\ m=\left(\begin{smallmatrix}
B(t)n\\
n
\end{smallmatrix}\right)\bigr\}$ is a basis of $p^{*}\XClAlg\big(t^{\prin}\big)$. It implies the desired
claim.\looseness=-1
\end{proof}

\section{Examples}\label{sec:Examples}
We work at the classical level $\kk=\Z$ in this section for simplicity.
\subsection{Twist automorphisms on unipotent cells}
\cite{BerensteinFominZelevinsky96,berenstein1997total} introduced
an automorphism $\eta_{w}$ on the unipotent cell $\C[N_{-}^{w}]$,
which was called a twist automorphism. The quantized twist automorphism
$\eta_{w}$ was introduced by \cite{kimura2017twist}. We refer the
reader to loc. cit. for further details.

Moreover, it is known that (quantized) $\C[N_{-}^{w}]$ is a (quantum)
upper cluster algebra $\upClAlg$ by \cite{GeissLeclercSchroeer10,GeissLeclercSchroeer11,GY13,goodearl2020integral}. By \cite{qin2020dual}, the (quantum)
twist automorphism $\eta$ is a twist automorphism in our sense.

Let us illustrate the twist automorphism $\eta_{w}$ by an example.

\begin{Example}\label{eg:twist_auto_matrix}
Consider the algebraic group $G=SL_{3}(\kk)$ where $\kk$ is a field
of characteristic $0$. Let $B_{+}$ (resp.\ $B_{-}$) denote its Borel
subgroup consisting of the upper (resp.\ lower) triangular matrices
in $G$. Let $N_{\pm}$ denote the subgroups of $B_{\pm}$ whose diagonal
entries are $1$, respectively.

We take $w$ to be the longest element $w_{0}$ in the permutation
group of $3$ elements. It can be represented by the matrix{\samepage \[\trans_{w_{0}}=\begin{pmatrix}
0 & 0 & 1\\
0 & 1 & 0\\
1 & 0 & 0
\end{pmatrix}.\] The corresponding (negative) unipotent cell is $N_{-}^{w_{0}}:=B_{+}\trans_{w_{0}}B_{+}\cap N_{-}$.}

Let us denote an element $x\in N_{-}$ by \[x=\begin{pmatrix}
1 & 0 & 0\\
x_{1} & 1 & 0\\
x_{2} & x_{1}' & 1
\end{pmatrix}.\] Let $x_{3}$ denote the minor \smash{$\bigl|\left(\begin{smallmatrix}
x_{1} & 1\\
x_{2} & x_{1}'
\end{smallmatrix}\right)\bigr|=x_{1}x_{1}'-x_{2}$}. By \cite{berenstein1997total}, there exists a unique element $y$
in $N_{-}\cap B_{+}\trans_{w_{0}}x^{\rm T}$ when $x\in N_{-}^{w_{0}}$.
The twist automorphism $\tilde{\eta}_{w_{0}}$ on $N_{-}$ is defined
such that $\tilde{\eta}_{w_{0}}(x)=y$.

Let us compute $y$ explicitly. Take any element \[z=\begin{pmatrix}
a & b & c\\
 & d & e\\
 & & f
\end{pmatrix}\in B_{+}.\] Note that \[\trans_{w_{0}}x^{\rm T}=\begin{pmatrix}
0 & 0 & 1\\
0 & 1 & x_{1}'\\
1 & x_{1} & x_{2}
\end{pmatrix}.\] We have
\begin{align*}
z\trans_{w_{0}}x^{\rm T} & =\begin{pmatrix}
c & b+cx_{1} & a+bx_{1}'+cx_{2}\\
e & d+ex_{1} & dx_{1}'+ex_{2}\\
f & fx_{1} & fx_{2}
\end{pmatrix}.
\end{align*}
If $z\trans_{w_{0}}x^{\rm T}\in N_{-}$, then we must have $c=1$, $f=x_{2}^{-1}$,
$b=-cx_{1}$, $a=-bx_{1}'-cx_{2}$, and \[\begin{cases}
d+ex_{1} =1,\\
dx_{1}'+ex_{2} =0.
\end{cases} \] The last linear system gives us $d=\frac{x_{2}}{x_{2}-x_{1}x_{1}'}=\frac{x_{2}}{x_{3}}$,
$e=\frac{-x_{1}'}{x_{2}-x_{1}x_{1}'}=\frac{x_{1}'}{x_{3}}$. We deduce
that $N_{-}\cap B_{+}\trans_{w_{0}}x^{\rm T}$ consists of a single element
\begin{align*}
y & =\begin{pmatrix}
1 & 0 & 0\\
x_{1}'x_{3}^{-1} & 1 & 0\\
x_{2}^{-1} & x_{1}x_{2}^{-1} & 1
\end{pmatrix}
\end{align*}
if $x_{2},x_{3}\neq0$, and it is empty if not.

It is easy to check that $\big\{N_{-}\cap B_{+}\trans_{w_{0}}x^{\rm T}\mid \forall x\in N_{-}\big\}=N_{-}\cap B_{+}\trans_{w_{0}}B_{+}=:N_{-}^{w_{0}}$.
Then we can deduce from the above computation of $y$ that $N_{-}^{w_{0}}=\{x\in N_{-}\mid x_{2}\neq0,x_{3}\neq0\}$.

Define $A_{1}$ to be the function on $N_{-}^{w_{0}}$ sending a point
$x$ to $x_{1}$. Similarly define $A_{1}'$,~$A_{2}$,~$A_{3}$.
Then the coordinate ring $\kk[N_{-}^{w_{0}}]$ is an (upper) cluster
algebra with four cluster variables: $A_{1}$,~$A_{1}'$,~$A_{2}$,~$A_{3}$, see \cite[Example~2.2]{qin2021cluster}. Then $\tilde{\eta}_{w_{0}}$
gives rise to an automorphism $\eta_{w}$ such that
\begin{gather*}
\eta_{w_{0}}(A_{1}) =A_{1}'A_{3}^{-1}=\frac{A_{2}+A_{3}}{A_{1}}A_{3}^{-1},\\
\eta_{w_{0}}(A_{1}') =A_{1}A_{2}^{-1}=\frac{A_{2}+A_{3}}{A_{1}'}A_{2}^{-1},\\
\eta_{w_{0}}(X_{i}) =A_{i}^{-1},\qquad i=2,3.
\end{gather*}
\end{Example}

\subsection{Dehn twists for surface cases}

Let $\Sigma$ denote a triangulable surface $S$ with finitely many
marked point. For any of its (tagged) triangulation $\Delta$, one
can construct a seed $t_{\Delta}$, whose frozen vertices are contributed
from curves on the boundary $\partial S$. We refer the reader to
\cite{FominShapiroThurston08,Muller12} for details.

The punctures are the marked point in the interior $S^{\circ}$ of
$S$. By \cite{Muller12}, the seed $t_{\Delta}$ satisfies the full
rank condition when $\Sigma$ has no punctures.

Let $L$ denote any closed loop in the interior of $S$ which does
not pass a marked point or is contractible to a marked point. A Dehn
twist $\tw_{L}$ around $L$ produces a new (tagged) triangulation
$\tw_{L}\Delta$ from $\Delta$. By \cite{FominShapiroThurston08},
we have $t_{\tw_{L}\Delta}=\seq t_{\Delta}$ for some mutation sequence
$\seq$. By construction, $B(t_{\Delta})=B(t_{\tw_{L}\Delta})$. In
particular, $t_{\Delta}$ and $t_{\tw_{L}\Delta}$ are similar.

Then one can construct the variation maps $\var_{t_{\Delta}}^{M}\colon\Mc(t_{\Delta})\rightarrow\Mc(t_{\tw_{L}\Delta})$
and $\var_{t_{\Delta}}^{N}\colon N(t_{\Delta})\rightarrow N(t_{\tw_{L}\Delta})$
represented by the identity matrices, respectively. The corresponding
twist automorphisms are determined by
\begin{align*}
\var_{t_{\Delta}}^{A}(A_{i}(t_{\Delta})) & =\seq^{*}A_{i}(t_{\tw_{L}\Delta}),\qquad
\var_{t_{\Delta}}^{X}(X_{i}(t_{\Delta})) =\seq^{*}X_{i}(t_{\tw_{L}\Delta}).
\end{align*}
They give rise to automorphisms on $\upClAlg(t_{\Delta})$ and $\XClAlg(t_{\Delta})$
associated to the Dehn twist $\tw_{L}$.

\subsection{Once-punctured digon}

\begin{Example}
The following cluster Poisson algebra arises from ${\rm PGL}_{3}$-local
systems on a once-punctured digon \cite{schrader2019cluster}.

Take $I=\{1,2,3,4\}$, with $\Iufv=\{2,4\}$ and $\Ifv=\{1,3\}$.
Define the initial seed $t$ such that its $B$-matrix is \[B=\begin{pmatrix}
\hphantom{-}0 & -1 & \hphantom{-}0 & \hphantom{-}1\\
\hphantom{-}1 & \hphantom{-}0 & -1 & \hphantom{-}0\\
\hphantom{-}0 & \hphantom{-}1 & \hphantom{-}0 & -1\\
-1 & \hphantom{-}0 & \hphantom{-}1 & \hphantom{-}0
\end{pmatrix}.\] Note that $B_{\ufv}=0$. So we cannot endow a compatible Poisson
structure $\lambda$ for $\LP^{A}(t)$. We have $W=-B$.

Let $t'$ denote the seed $\seq t$ where $\seq=\mu_{2}\mu_{4}$.
Its $B$-matrix is $B'=-B$. We have $W'=-W$. We can check the following
mutation rule:
\begin{gather*}
\seq^{*}X_{2}' =X_{2}^{-1},\qquad
\seq^{*}X_{4}' =X_{4}^{-1},\qquad
\seq^{*}X_{1}' =X_{1}\cdot X_{2}\cdot(1+X_{2})^{-1}\cdot(1+X_{4}),\\
\seq^{*}X_{3}' =X_{3}\cdot X_{4}\cdot(1+X_{2})\cdot(1+X_{4})^{-1}.
\end{gather*}
Then $t'=t[1]$ with $\sigma=\Id$, $W'=-W$.

Let us construct Poisson automorphisms on $\XClAlg(t)$. Note that
we have \smash{$W=\left(\begin{smallmatrix}
0 & -W_{\low}^{\rm T}\\
W_{\low} & 0
\end{smallmatrix}\right)$} with row and column indices $(2,4,1,3)$, where \smash{$W_{\low}=\left(\begin{smallmatrix}
1 & -1\\
-1 & 1
\end{smallmatrix}\right)$} with row indices $\{1,3\}$ and column indices $\{2,4\}$. We want
to construct a linear isomorphism $\var_{t}^{N}:=\colon N(t)\rightarrow N(t')$
represented by an $\Z$-valued invertible matrix \smash{$\left(\begin{smallmatrix}
\Id_{\ufv} & V_{\high}\\
0 & V_{\fv}
\end{smallmatrix}\right)$}, such that
\begin{align*}
\begin{pmatrix}
\Id_{\ufv} & V_{\high}\\
0 & V_{\fv}
\end{pmatrix}^{\rm T}W'\begin{pmatrix}
\Id_{\ufv} & V_{\high}\\
0 & V_{\fv}
\end{pmatrix} & =W.
\end{align*}
The desired equation is equivalent to
\begin{align*}
-\begin{pmatrix}
\Id_{\ufv} & V_{\high}\\
0 & V_{\fv}
\end{pmatrix}^{\rm T}\begin{pmatrix}
0 & -W_{\low}^{\rm T}\\
W_{\low} & 0
\end{pmatrix}\begin{pmatrix}
\Id_{\ufv} & V_{\high}\\
0 & V_{\fv}
\end{pmatrix} & =\begin{pmatrix}
0 & -W_{\low}^{\rm T}\\
W_{\low} & 0
\end{pmatrix}.
\end{align*}
By computing its block submatrices, we reduce the equation to the
equations
\begin{align*}
-V_{\fv}^{\rm T}W_{\low} & =W_{\low},\qquad
V_{\high}^{\rm T}W_{\low}^{\rm T}V_{\fv}-V_{\fv}^{\rm T}W_{\low}V_{\high} =0.
\end{align*}

For the first equation (for $\var_{t}^{N}$ to be a variation map),
the solution takes the form $V_{\fv}=\smash{\left(\begin{smallmatrix}
\lambda-1 & \mu\\
\lambda & \mu-1
\end{smallmatrix}\right)}$ for any $\lambda$, $\mu$. Since we are looking for a bijection $\var_{t}^{N}$
between lattices, we must have $\lambda+\mu=0,2$.

The second equation (for $\var_{t}^{N}$ to further be Poisson) is
equivalent to \[-V_{\high}^{\rm T}W_{\low}^{\rm T}+W_{\low}V_{\high}=0\] by
the first equation. Then its solution takes the form \smash{$V_{\high}=\left(\begin{smallmatrix}
\alpha & \beta\\
\alpha & \beta
\end{smallmatrix}\right)$} for any $\alpha$, $\beta$.

The resulting twist automorphism on the fraction field takes the following
form:
\begin{gather*}
\tw X_{2} =X_{2}^{-1},\qquad
\tw X_{4} =X_{4}^{-1},\\
\tw X_{1} =\seq^{*}((X')^{(\lambda-1,\alpha,\lambda,\alpha)})=X^{(\lambda-1,\lambda-1-\alpha,\lambda,\lambda-\alpha)}(1+X_{2})(1+X_{4})^{-1},\\
\tw X_{3} =\seq^{*}((X')^{(\mu,\beta,\mu-1,\beta)})=X^{(\mu,\mu-\beta,\mu-1,\mu-1-\beta)}(1+X_{2})^{-1}(1+X_{4}).
\end{gather*}

Note that $\tw(X_{1}X_{3})=(X_{1}X_{3})^{\lambda+\mu-1}(X_{2}X_{4})^{\lambda+\mu-1-\alpha-\beta}$.

Let us consider elements $E=X_{1}\cdot(1+X_{4})$, $F=X_{3}\cdot(1+X_{2})$,
$K=X_{1}\cdot X_{3}\cdot X_{4}$ and $K'=X_{1}\cdot X_{2}\cdot X_{3}$
(the Chevalley generators, see \cite{schrader2019cluster}). It is
straightforward to check that
\begin{gather*}
\big(\seq^{-1}\big)^{*}E=X_{1}'\cdot(1+X_{2}'),\qquad \big(\seq^{-1}\big)^{*}F=X_{3}'\cdot(1+X_{2}'), \qquad\big(\seq^{-1}\big)^{*}K=X_{1}'X_{2}'X_{3}',\\
\big(\seq^{-1}\big)^{*}K'=X_{1}'X_{3}'X_{4}',
\end{gather*}
 so they are elements in the
cluster Poisson algebra $\XClAlg$.

We check that
\begin{align*}
&\tw E =X^{(\lambda-1,\lambda-1-\alpha,\lambda,\lambda-\alpha)}(1+X_{2})(1+X_{4})^{-1}\big(1+X_{4}^{-1}\big)=(X_{1}X_{3})^{\lambda-1}(X_{2}X_{4})^{\lambda-1-\alpha}F,\\
&\tw F =X^{(\mu,\mu-\beta,\mu-1,\mu-1-\beta)}(1+X_{4})(1+X_{2})^{-1}\big(1+X_{2}^{-1}\big)=(X_{1}X_{3})^{\mu-1}(X_{2}X_{4})^{\mu-1-\beta}E,\\
&\tw K =(X_{1}X_{3})^{\lambda+\mu-1}(X_{2}X_{4})^{\lambda+\mu-1-\alpha-\beta}X_{4}^{-1},\\
&\tw K' =(X_{1}X_{3})^{\lambda+\mu-1}(X_{2}X_{4})^{\lambda+\mu-1-\alpha\beta}X_{2}^{-1}.
\end{align*}
For the special case $\lambda=\mu=1$, $\alpha=\beta=0$, we recover
the following automorphism in (Drinfeld double) quantum groups:
\begin{align*}
\tw E & =F,\qquad
\tw F =E,\qquad
\tw K =K',\qquad
\tw K' =K.
\end{align*}
\end{Example}

\subsection*{Acknowledgements}

The authors express their sincere thanks to the anonymous referees for their insightful comments and constructive suggestions. The first author is supported by JSPS KAKENHI Grant Number 17K14168 and 22H01114. The second author is supported by NSFC Grant Number 12271347.

\pdfbookmark[1]{References}{ref}
\LastPageEnding

\end{document}